\documentclass[a4paper,12pt]{amsart}

\numberwithin{equation}{section}
\theoremstyle{plain}
\newtheorem{theorem}{\sc \bf Theorem}[section]
\newtheorem{lemma}[theorem]{\sc \bf Lemma}
\newtheorem{corollary}[theorem]{\sc \bf Corollary}
\newtheorem{proposition}[theorem]{\sc \bf Proposition}
\newtheorem{claim}{\sc \bf Claim}
\theoremstyle{definition} 
\newtheorem{definition}[theorem]{\sc \bf Definition}
\newtheorem{remark}[theorem]{\sc \bf Remark}
\newtheorem{example}[theorem]{\sc \bf Example}

\usepackage{amssymb}
\usepackage{amsfonts}
\usepackage{amsmath}
\usepackage{mathrsfs}   
\usepackage{graphicx}
\usepackage{latexsym}
\usepackage{accents}
\usepackage{enumerate}
\usepackage{ulem}
\usepackage{fancybox}
\usepackage{wasysym}
\usepackage{color}
\usepackage{cancel}
\usepackage{arydshln}

\newif\ifINITPACKNOTE
\INITPACKNOTEtrue

\newif\ifINITPACKKOUGIFUNC
\INITPACKKOUGIFUNCfalse

\makeatletter
\newenvironment{proofex}[1][\proofname]{\par
  \normalfont
  \topsep6\p@\@plus6\p@ \trivlist
  \item[\hskip\labelsep{\bfseries #1}\@addpunct{\bfseries}]\ignorespaces
}{
   \popQED\endtrivlist\@endpefalse
}
\makeatother

\newcounter{constants}
\makeatletter
\def\addconst{
\addtocounter{constants}{1}			
\def\@currentlabel{\arabic{constants}}	
\@currentlabel							
}
\makeatother
\newcommand{\adl}[1]{\addconst\label{c:#1}}
\newcommand{\adr}[1]{\ref{c:#1}}

\newcommand{\K}{{\mathbb{K}}}

\let\R=\relax
\let\C=\relax
\let\Z=\relax

\newcommand{\R}{{\mathbb{R}}}
\newcommand{\C}{{\mathbb{C}}}
\newcommand{\Z}{{\mathbb{Z}}}

\definecolor{COLgrey}{rgb}{0.65, 0.65, 0.65}
\definecolor{COLblue}{rgb}{0.65, 0.25, 0.25}
\definecolor{brickred}{cmyk}{0,0.89,0.94,0.28}

\definecolor{skyblue}{rgb}{0.82,0.94,1.0} 
\definecolor{blackblue}{rgb}{0.2,0.2,0.6}

\newcommand{\Si}{\Sigma}
\newcommand{\SiAp}{\Sigma_{A}^{+}}

\newcommand{\si}{\sigma}

\newcommand{\ph}{\varphi}

\newcommand{\phe}{\varphi(\epsilon,\cdot)}

\ifINITPACKNOTE
\let\lam=\relax
\let\e=\relax
\newcommand{\lam}{\lambda}
\newcommand{\e}{\epsilon}
\else
\fi

\let\p=\relax
\newcommand{\p}{{\prime}}
\newcommand{\pp}{{\prime\prime}}

\newcommand{\AR}{\mathcal{A}}
\newcommand{\BR}{\mathcal{B}}
\newcommand{\CR}{\mathcal{C}}
\newcommand{\DR}{\mathcal{D}}

\newcommand{\IR}{\mathcal{I}}

\newcommand{\LR}{\mathcal{L}}
\let\MR=\relax
\newcommand{\MR}{\mathcal{M}}
\newcommand{\NR}{\mathcal{N}}
\newcommand{\OR}{\mathcal{O}}

\newcommand{\RR}{\mathcal{R}}

\newcommand{\XR}{\mathcal{X}}

\newcommand{\tLR}{\tilde{\mathcal{L}}}

\newcommand{\ali}[1]{\begin{align*}#1\end{align*}}
\newcommand{\alil}[1]{\begin{align}#1\end{align}}
\newcommand{\ite}[1]{\begin{enumerate}[(1)]#1\end{enumerate}}
\newcommand{\ites}{\begin{enumerate}}
\newcommand{\itee}{\end{enumerate}}
\let\prop=\relax
\let\cor=\relax
\newcommand{\prop}[1]{\begin{proposition}#1\end{proposition}}
\newcommand{\props}{\begin{proposition}}
\newcommand{\prope}{\end{proposition}}
\newcommand{\cor}[1]{\begin{corollary}#1\end{corollary}}
\newcommand{\cors}{\begin{corollary}}
\newcommand{\core}{\end{corollary}}
\newcommand{\thm}[1]{\begin{theorem}#1\end{theorem}}
\newcommand{\thms}{\begin{theorem}}
\newcommand{\thme}{\end{theorem}}
\newcommand{\lem}[1]{\begin{lemma}#1\end{lemma}}
\newcommand{\lems}{\begin{lemma}}
\newcommand{\leme}{\end{lemma}}

\newcommand{\defis}{\begin{definition}}
\newcommand{\defie}{\end{definition}}

\newcommand{\exams}{\begin{example}}
\newcommand{\exame}{\end{example}}
\newcommand{\rem}[1]{\begin{remark}\normalfont #1\end{remark}}

\newcommand{\pros}{\begin{proof}}
\newcommand{\proe}{\end{proof}}
\newcommand{\pross}{\begin{proofex}}
\newcommand{\prossb}{\begin{proofex}\bf}
\newcommand{\smilemk}{\ifmmode\else\leavevmode\unskip\penalty9999\hbox{}\nobreak\hfill\fi\quad\hbox{\smileyex}}
\newcommand{\proes}{\smilemk\end{proofex}}
\newcommand{\prossn}{\begin{proofex}}
\newcommand{\proesn}{\end{proofex}}

\newcommand{\prossq}{\begin{proof}
\begin{enumerate}
\def\theenumi{(\arabic{enumi})}
\def\labelenumi{\theenumi}
\setlength{\leftmargin}{5pt}	
\setlength{\itemsep}{3pt}		
\setlength{\parskip}{0.0pt}		
\setlength{\itemindent}{15pt}	
\setlength{\leftskip}{-35pt}		
\setlength{\labelsep}{3pt}		
}
\newcommand{\case}[1]{\begin{cases}#1\end{cases}}
\newcommand{\cd}{\cdot}
\newcommand{\up}{\upsilon}
\newcommand{\om}{\omega}

\newcommand{\ti}[1]{\tilde{#1}}

\newcommand{\mat}[1]{\begin{matrix}#1\end{matrix}}

\newcommand{\matII}[1]{\left(\begin{array}{cc}#1\end{array}\right)}

\newcommand{\MatII}[1]{\left[\begin{array}{cc}#1\end{array}\right]}

\newcommand{\supp}{\mathrm{supp}}
\newcommand{\di}{\displaystyle}
\ifINITPACKKOUGIFUNC

\else
\ifINITPACKNOTE

\else

\fi
\fi

\let\MARU=\relax
\newcommand{\MARU}[1]{{\ooalign{\hfil#1\/\hfil\crcr\raise.167ex\hbox{\mathhexbox20D}}}}

\let\to=\relax
\newcommand{\to}{\ \rightarrow\ }

\newcommand{\qqquad}{\quad\qquad}

\newcommand{\graphexp}[4]{\mbox{\raisebox{-#3\zw}{\includegraphics[scale=#2,angle=0,page=#4]{#1}}}}

\newcommand{\zw}{em}

\begin{document}
\title[Markov systems with holes]{Perturbed infinite-state Markov systems with holes and its application}
\author[H. Tanaka]{Haruyoshi Tanaka}
\address{
{\rm Haruyoshi Tanaka}\\
Course of Mathematics Education\\
Education Graduate School of Education, Naruto University of Education\\
748, Nakajima, Takashima, Naruto-cho, Naruto-City, Tokushima, 772-8502, Japan
}
\email{htanaka@naruto-u.ac.jp}
\begin{abstract}
We consider a perturbed system $(X,\phe)$, where $X$ is a topological Markov shift with a countably infinite state space, and $\phe$ is a real-valued potential on X depending on a small parameter $\e\in (0,1)$. We assume that for each $\e>0$, the system has a unique transitive component and a unique Gibbs measure (or more generally, a Ruelle-Perron-Frobenius (RPF) measure) $\mu_{\e}$, while the unperturbed system possesses multiple transitive components and Gibbs measures on these components. We investigate the convergence of the measure $\mu_{\e}$ as $\e\to 0$ and the representation of the limiting measure, if it exists. In previous work [T. 2020], we considered the finite state case. Our approach relies on a development of the Schur–Frobenius factorization theorem, which we apply to demonstrate a spectral gap property for Perron complements of Ruelle operators in the infinite-state case. As an application, we examine perturbed piecewise expanding Markov maps with holes, defined over countably infinite partitions. We investigate the splitting behavior of the associated Gibbs measure under perturbation.
\end{abstract}
\maketitle
\tableofcontents
\section{Introduction and main results}\label{sec:intro}
Denote by $X=X_{A}$ the topological Markov shift with s countable set of states $S$ and with a transition matrix $A=(A(ij))_{S\times S}$ (see Section \ref{sec:TMS} for terminology). A word $w_{1}w_{2}\dots w_{n}\in S^{n}$ is said to be $A$-admissible if $A(w_{i}w_{i+1})=1$ for all $1\leq i<n$.
For an $A$-admissible word $w\in S^{n}$, we denote by $[w]$ the {\it cylinder set} $\{\om\in X\,:\,\om_{0}\cdots \om_{n-1}=w\}$. In what follows, we assume that $A$ is {\it irreducible}, that is, for any two states $a,b\in S$, there exists a word $w_{1}\cdots w_{n}\in S^{n}$ with $w_{1}=a$ and $w_{n}=b$ such that $w$ is $A$-admissible. We consider the following three conditions for functions $\ph(\e,\cd)\,:\,X\to \R$ with a small parameter $\epsilon\in (0,1)$:
\ite
{
\item[(A.1)] $[\phe]_{2}:=\sup_{\omega\neq \up\,:\,\om\in [\up_{0}\up_{1}]}\left|\ph(\e,\om)-\ph(\epsilon,\up)\right|/d_{\theta}(\omega,\up)$ is bounded uniformly in $\e$.
\item[(A.2)] $\sum_{a\in S}\exp(\sup_{\e>0}\sup_{\om\in [a]}\ph(\e,\om))$ is finite.
\item[(A.3)] There exists a function $\ph\,:\,X\to \R\cup\{-\infty\}$ such that for any $a\in S$, $\sup_{\om\in [a]}|e^{\ph(\e,\om)}-e^{\ph(\om)}|\to 0$ as $\e\to 0$, and $P(\ph)$ is finite, where $e^{-\infty}$ means $0$.
}
Here $d_{\theta}$ is defined in Section \ref{sec:Ruelleop}, and $P(\ph)$ denotes the topological pressure of $\ph$ which is given by (\ref{eq:toppres}). It follows from condition (A.1) that for any $a,b\in S$, $\ph(\om)$ is finite for some $\om\in [ab]$ if and only if $\ph(\om)$ is finite on $[ab]$. We define a zero-one matrix $B=(B(ij))$ such that $B(ij)=1$ if $\ph$ is finite on $[ij]$ and $0$ otherwise. For $a,b\in S$, we write $a\leftrightarrow b$ when either $a=b$ or there exist words $w_{1},w_{2}\in S^{*}:=\bigcup_{n=1}^{\infty}S^{n}$ such that $a\cd w_{1}\cd b$ and $b\cd w_{2}\cd a$ are $B$-admissible. Consider the quotient space $S/\!\!\leftrightarrow$. 
For a non-$1\times 1$ zero matrix $U\in S/\!\!\leftrightarrow$, we denote by $B[U]$ the irreducible submatrix of $B$ indexed by $U\times U$. Therefore we obtain countably many transitive components $X_{B[U]}$ $(U\in S/\!\!\leftrightarrow)$ of $X_{B}$, where $X_{B[U]}$ denotes the subsystem of $X$ with the transition matrix $B[U]$ and the state space $U$. Consider the {\it maximal pressure components}:
\ali
{
\mathscr{S}_{0}:=\{U\in S/\!\!\leftrightarrow\,:\,P(\ph)=P(\ph|X_{B[U]})\}:=\{U(1),U(2),\dots,U(m_{0})\}.
}
According to \cite{T2024_sum}, this set is non-empty and finite. Since $A$ is irreducible, there exists a unique {\it Ruelle-Perron-Frobenius (RPF) measure} $\mu_{\e}$ for $\phe$ which is given by $\mu_{\e}=h_{\e}\nu_{\e}$, $h_{\e}$ is the Perron eigenfunction of the Ruelle operator $\LR_{\phe}$ for $\phe$ and $\nu_{\e}$ is the Perron eigenvector of the dual of $\LR_{\phe}$ (see Section \ref{sec:Ruelleop} for the existence). Note that this measure is a generalization of a Gibbs measure (see also Theorem \ref{th:ex_Gibbs}). In contrast, the limiting potential $\ph$ has finitely many RPF measures $\mu(U,\cd)$ for $\ph|X_{B[U]}$ ($U\in \mathscr{S}_{0}$). For subset $U\subset S$, $\Si_{U}$ denote the set $\bigcup_{a\in U}[a]$, and $\chi_{U}$ denote the indicator of the set $\Si_{U}$. Denote by $F_{b}(X)$ the set of all complex-valued functions $f$ on $X_{A}$ that are continuous, have finite Lipschitz norm $\|f\|_{1}<+\infty$ (see Section \ref{sec:Ruelleop}), and by $F_{b,0}(X)$ the set of $f\in F_{b}(X)$ that {\it vanish at infinity}, that is, for any $\eta>0$ there exists a finite set $Q\subset S$ such that $|f|<\eta$ on $\Si_{S\setminus Q}$.

Now we are in a position to state one of our main results:
\thm
{\label{th:main_ac}
Assume that conditions (A.1)-(A.3) are satisfied. Then for any positive sequence $\e(n)$ with $\e(n)\to 0$, there exists a subsequence $\e^\p(n)$ of $\e(n)$ such that $\mu_{\e^\p(n)}$ converges weakly on $F_{b,0}(X)$ as $n\to \infty$, and the limit point has the form $\sum_{k=1}^{m_{0}}\delta(k)\mu(U(k),\cd)$ for some numbers $\{\delta(k)\}$ with $\sum_{k=1}^{m_{0}}\delta(k)\leq 1$.
}
This is a corollary of Theorem \ref{th:main} and will be proved in Section \ref{sec:proof_main1}.

Next we state how the coefficients $\delta(k)$ is determined. To do this, we introduce transfer-type operators associated with induced shifts. For distinct subsets $U, V\subset S$ with $U\neq \emptyset$, put
$W(U:V)=U\times\bigcup_{n=0}^{\infty}V^{n}$, where the set $\bigcup_{n=0}^{\infty}V^{n}$ includes the empty word. Let $\lam_{\e}:=\exp(P(\phe))$. We define a bounded linear operator $\LR_{\phe}[U,V,\lam_{\e}]$ acting on $F_{b}(X)$ by
\alil
{
\LR_{\phe}[U,V,\lam_{\e}] f(\omega):=\sum_{w\in W(U\,:\,V)\,:\,w\cdot \om_{0}:A\text{-admissible}}\lam_{\e}^{-|w|+1}e^{S_{|w|}\ph(\e,w\cdot \om)}f(w\cdot \om)\label{eq:LRekkp=...}
}
for $f\in F_{b}(X)$ and $\om\in X$, where $|w|$ denotes the length of the word $w$ and $S_{n}\ph(\e,\upsilon)=\sum_{j=0}^{n-1}\ph(\e,\si^{j}\upsilon)$. Here, the summation in (\ref{eq:LRekkp=...}) is taken over all words $w$ so that the itinerary $\up=w\cd\om$ starts at a state in $U$ and reaches the state $\om_{0}$ passing through states in $V$ on the way. In other words, this operator is a transfer operator with respect to an induced shift. Note also that this operator can be interpreted as a {\it Perron complement of Ruelle operators}, as first introduced in \cite{T2020} (see Section \ref{sec:proof} for the general setting).

For any finite subset $Q\subset S$ with $Q\cap U(k)\neq \emptyset$ for all $1\leq k\leq m_{0}$, and distinct integers $1\leq i,j\leq m_{0}$, 
we denote by $c_{\e}(Q,i,j)$ the spectral radius of the operator $\LR_{\phe}[Q\cap U(i),S\setminus ((U(i)\cup U(j))\cap Q),\lam_{\e}]$. For $1\leq i\leq m_{0}$, we define
\alil
{
\delta_{\e}(Q,i)=\frac{1}{\di 1+\sum_{1\leq j\leq m_{0}\,:\,j\neq i}\frac{\lam_{\e}-c_{\e}(Q,i,j)}{\lam_{\e}-c_{\e}(Q,j,i)}}.\label{eq:deQi=}
}
Now we state our second result.
\thm
{\label{th:main}
Assume that conditions (A.1)-(A.3) are satisfied, and $m_{0}\geq 2$. Then the following are equivalent:
\ite
{
\item There exists a Borel measure $\mu$ on $X$ such that for any $f\in F_{b,0}(X)$, $\mu_{\e}(f)$ converges to $\mu(f)$.
\item For any finite subset $Q\subset S$ with $Q\cap U(k)\neq \emptyset$ for each $1\leq k\leq m_{0}$, $\mu_{\e}(\chi_{Q})$ converges to a number $a(Q)$ and if $a(Q)>0$ then $\delta_{\e}(Q,k)$ converges to a number $\delta(Q,k)$ for all $1\leq k\leq m_{0}$ as $\e\to 0$.
\item There exists a sequence of finite subsets $Q_{n}\subset S$ with $Q_{n}\cap U(k)\neq \emptyset$ for all $1\leq k\leq m_{0}$, $Q_{n}\subset Q_{n+1}$ and $\bigcup_{n}Q_{n}=S$ such that $a(n):=\lim_{\e\to 0}\mu_{\e}(\chi_{Q_{n}})$ exists and if $a(n)>0$ then $\delta_{\e}(Q_{n},k)$ converges to a number $\delta(Q_{n},k)$ as $\e\to 0$ for all $k$.
\item There exists a sequence of finite subsets $Q_{n}\subset S$ with $Q_{n}\cap U(k)\neq \emptyset$ for all $k$, $Q_{n}\subset Q_{n+1}$ and $\bigcup_{n}Q_{n}=S$ such that $\delta(n,k):=\lim_{\e\to 0}\mu_{\e}(\chi_{Q_{n}})\delta_{\e}(Q_{n},k)$ exists for all $n\geq 1$ and for all $1\leq k\leq m_{0}$.
}
In this case, the limiting measure $\mu$ has the form $\sum_{k=1}^{m_{0}}\delta(k)\mu(U(k),\cd)$ with $\delta(k):=\lim_{n\to \infty}\lim_{\e\to 0}\mu_{\e}(\chi_{Q_{n}})\delta_{\e}(Q_{n},k)$.
}
We will show this in Section \ref{sec:proof_main2}.

We introduce an additional condition:
\ite
{
\item[(A.4)] For any $\eta\in (0,1)$ there exist a finite set $Q\subset S$ and $\e_{0}\in (0,1)$ such that $\mu_{\e}(\Si_{S\setminus Q})\leq \eta$ for any $0<\e<\e_{0}$.
}
\cor
{\label{cor:main}
Assume that conditions (A.1)-(A.4) are satisfied and $m_{0}\geq 2$. Then the following are equivalent:
\ite
{
\item There exists a Borel measure $\mu$ on $X$ such that for any $f\in F_{b}(X)$, $\mu_{\e}(f)$ converges to $\mu(f)$.
\item There exists a sequence of finite subsets $Q_{n}\subset S$ with $Q_{n}\cap U(k)\neq \emptyset$ for all $k$, $Q_{n}\subset Q_{n+1}$ and $\bigcup_{n}Q_{n}=S$ such that $\delta(n,k):=\lim_{\e\to 0}\delta_{\e}(Q_{n},k)$ exists for all $n\geq 1$ and for all $1\leq k\leq m_{0}$.
}
In this case, $\mu$ has the form $\sum_{k=1}^{m_{0}}\delta(k)\mu(U(k),\cd)$ with $\delta(k):=\lim_{n\to \infty}\lim_{\e\to 0}\delta_{\e}(Q_{n},k)$ and $\mu$ is a probability measure. Consequently, any accumulation point of $\{\mu_{\e}\}$ is a convex combination of $\{\mu(U(k),\cd)\}_{k}$. 
}
This will be proved in Section \ref{sec:proof_main2_cor}
\rem
{
\ite
{
\item Note that if $\# \mathscr{S}_{0}=1$ then $\mu_{\e}$ on $F_{b}(X)$ always converges weakly to $\mu(1,\cd)$ (see \cite{T2024_sum}).
\item If $S$ is finite then $F_{b,0}$ coincides with the set of bounded Lipschitz continuous functions on $X$. Moreover, we may take $Q=Q_{n}=S$ in (2)(3). Thus, our result is a generalization from the finite-state case \cite{T2020} to the infinite-state case.
\item We expect that conditions (A.1)-(A.3) imply condition (A.4). However, this has not been proven yet, and it will be the subject of future work.
}
}
As an application, we apply our results to piecewise expanding Markov maps with countably infinite partitions. We will show a split of the perturbed Gibbs measure for the perturbed physical potential as a convex combination of the unperturbed Gibbs measures (see Theorem \ref{th:Gibbsiff_PEMM}). Our results play an important role in studying the stable and unstable behavior of metastable systems. In our future work, we shall apply our results to problem of detecting metastable states. In particular, we aim to study specific specific metastable states via perturbation methods.

In the case when the state set $S$ is finite, we gave a necessary and sufficient condition for the convergence of $\mu_{\e}$ using the notion of Perron complements of Ruelle operators such as (\ref{eq:LRekkp=...}) \cite{T2020}. Moreover, we also proved that any limit point is a convex combination of $\{\mu(i,\cd)\}$. In the present paper, we investigate the corresponding problem in the case where $S$ is infinite. The technique of \cite{T2020} is used for the proof of our main results; however, many of the auxiliary results necessary for the proof cannot be derived using the same method as in \cite{T2020}. For example, it is difficult to obtain a spectral gap property for the transfer operators (\ref{eq:LRekkp=...}) in our main result using techniques from the finite-state case. We establish such a property by using a new method, namely a modified Schur-Frobenius factorization (Theorem \ref{th:PC_FS}). Furthermore, the proofs of the convergence of restricted Perron eigenfunctions and Perron eigenvectors are substantially modified (Section \ref{sec:TSC}). 
\smallskip
\par
In Section \ref{sec:TMS}, we recall the notion of symbolic dynamics with infinitely many states. The concept of thermodynamic formalism and Ruelle transfer operators, which drives thermodynamic features, is described in Section \ref{sec:Ruelleop}. In Section \ref{sec:subRuop_SFfactor}, we develop the Perron complements of transfer operators from the finite state case in \cite{T2020} to the infinite-state setting. In particular, we establish a connection between these operators and the so-called Schur-Frobenius factorization theorem in Theorem \ref{th:PC_FS}. Section \ref{sec:TSC} is devoted to convergence of eigenvalues and eigenfunctions of transfer operators, as well as eigenvectors of their dual operators. Using these results and the techniques from \cite{T2020}, the main results are proved in Section \ref{sec:proof}. In the final section \ref{sec:piecewiseCI}, we apply our results to piecewise expanding Markov maps with countably infinite partitions.
\medskip
\\
\noindent
{\it Acknowledgment.}\ 
This study was supported by JSPS KAKENHI Grant Number 20K03636.
\section{Auxiliary results}\label{sec:prelim}
\subsection{Topological Markov shifts}\label{sec:TMS}
In this section, we recall the notion of topological Markov shifts (TMS) introduced by \cite{Sar09,Sar99}, which is needed to state our main results. 

Let $S$ be a countable set equipped with the distinct topology, and $A=(A(ij))_{S\times S}$ a zero-one matrix, i.e. $A(ij)=0$ or $1$ for all $i,j\in S$.
Consider the set
\ali
{
\textstyle{ X=X_{A}=\{\om=\om_{0}\om_{1}\cdots\in\prod_{n=0}^{\infty}S\,:\,A(\om_{k}\om_{k+1})=1 \text{ for any }k\geq 0\}}
}
endowed with the subspace topology induced by the product topology on $\prod_{n=0}^{\infty}S$, and with the shift transformation $\si\,:\,X\to X$ defined by $(\si\om)_{n}=\om_{n+1}$ for any $n\geq 0$. This is called a {\it topological Markov shift} with state space $S$ and with transition matrix $A$.
For $\theta\in (0,1)$, a metric $d_{\theta}\,:\,X\times X\to \R$ on $X$ is given by $d_{\theta}(\om,\up)=\theta^{\min\{n\geq 0\,:\,\om_{n}\neq \up_{n}\}}$ if $\om\neq \up$ and $d_{\theta}(\om,\up)=0$ if $\om=\up$.
It is known that $(X,d_{\theta})$ is a complete and separable metric space; and however, it is not compact in general.

Next, we introduce subsystems on $X$. Fix a nonempty subset $S_{0}\subset S$ and an $S_{0}\times S_{0}$ zero-one matrix $M=(M(ij))$ with $M(ij)\leq A(ij)$. Consider the set
\alil
{
\textstyle{X_{M}=\{\om\in \prod_{n=0}^{\infty}S_{0}\,:\,M(\om_{k}\om_{k+1})=1\text{ for all }k\geq 0\}}.\label{eq:XM=...}
}
That is, $X_{M}$ is a subset of $X_{A}$ and is the topological Markov shift with the state space $S_{0}$ and the transition matrix $M$. For convenience, we set $M(ij)=0$ whenever $i\in S\setminus S_{0}$ or $j\in S\setminus S_{0}$.

A transition matrix $M$ is said to be {\it finitely irreducible} if there exists a finite subset $F$ of $\bigcup_{n=1}^{\infty}S^{n}$ such that for any $a,b\in S$, there is $w\in F$ so that $a\cd w\cd b$ is $M$-admissible. 
The matrix $M$ is said to have the {\it big images and pre-images} (BIP) property if there exists a finite set $\{a_{1},\cdots,a_{N}\}$ of $S$ such that for any $b\in S$, there exist $1\leq i,j\leq N$ such that $M(a_{i}b)=M(ba_{j})=1$. Note that $M$ is finitely irreducible if and only if $X_{M}$ is topologically transitive and $M$ has the BIP property.

Finally, we show a useful result:
\prop
{\label{prop:irresubmat}
Let $S$ be a countable set and $A$ an $S\times S$ zero-one irreducible matrix. Then for any non-empty finite subset $T\subset S$ there exists a finite set $S_{0}$ with $T\subset S_{0}\subset S$ such that $A|_{S_{0}\times S_{0}}$ is irreducible.
}
\pros
For any $i,j\in T$, there exist an integer $m(i,j)\geq 2$ and $w(i,j)=(w_{k}(i,j))_{k}\in S^{m(i,j)}$ such that $w_{1}(i,j)=i$, $w_{m(i,j)}(i,j)=j$ and $w_{1}(i,j)\cdots w_{m(i,j)}(i,j)\cd j$ is $A$-admissible. Let
\ali
{
S_{0}:=\{w_{k}(i,j)\,:\,i,j\in T,\ k=1,\dots, m(i,j)\}.
}
We will show that $B:=A|_{S_{0}\times S_{0}}$ is irreducible. For any $s,t\in S_{0}$, there exist $i_{1},i_{2},j_{1},j_{2}\in T$, $1\leq k_{1}\leq m(i_{1},j_{1})$ and $1\leq k_{2}\leq m(i_{2},j_{2})$ such that $s=w_{k_{1}}(i_{1},j_{1})$ and $t=w_{k_{2}}(i_{2},j_{2})$. Then
\ali
{
s\cd w_{k_{1}+1}(i_{1},j_{1})&\cdots w_{m(i_{1},j_{1})-1}(i_{1},j_{1})\cd j_{1}\cd w_{2}(j_{1},i_{2})\cdots w_{m(j_{1},i_{2})-1}(j_{1},i_{2})\cd i_{2}\cd\\
&\cd w_{2}(i_{2},j_{2})\cdots w_{k_{2}-1}(i_{2},j_{2})\cd t
}
is $B$-admissible. Hence $B=A|_{S_{0}\times S_{0}}$ is irreducible.
\proe
\subsection{Thermodynamic formalism and Ruelle transfer operators}\label{sec:Ruelleop}
We recall some basic facts about Ruelle transfer operators, which were mainly introduced by \cite{MU}.
Let $X$ be a topological Markov shift with countable state space $S$ and transition matrix $A$.
\smallskip
\par
Next, we study transfer operators $\LR_{M,\ph}$ associated with subsystems of $X$, as defined in (\ref{eq:transfer}). 
Let $S$ be a countable set equipped with the distinct topology, and let $A=(A(ij))_{S\times S}$ be a zero-one matrix.
Consider the topological Markov shift $X=X_{A}$ with state space $S$ and transition matrix $A$.
A function $\ph\,:\,X\to \R$ is called {\it summable} if
\alil
{
\sum_{s\in S\,:\,[s]\neq \emptyset}\exp(\sup_{\om\in [s]}\ph(\om))<\infty.\label{eq:sum}
}
Let $\K=\R$ or $\C$. For a function $f\,:\,X\to \K$, $k\geq 1$ and $\theta\in (0,1)$, we define
\ali
{
[f]_{k}:=\sup_{w\in S^{k}\,:\,[w]\neq \emptyset}\sup\left\{\frac{|f(\om)-f(\up)|}{d_{\theta}(\om,\up)}\,:\,\om,\up\in [w],\ \om\neq \up\right\}.
}
Note that $[f]_{k}\leq [f]_{k+1}$. When $[f]_{1}<\infty$, $f$ is called {\it locally $d_{\theta}$-Lipschitz continuous}, and when $[f]_{2}<\infty$, $f$ is called {\it weak $d_{\theta}$-Lipschitz continuous}. Let $C_{b}(X)$ denote the Banach space of continuous functions $f\,:\,X\to \C$ with $\|f\|_{\infty}:=\sup_{\om\in X}|f(\om)|<+\infty$.
The space $F_{b}(X)$ is a Banach space under the Lipschitz norm $\|f\|_{1}:=\|f\|_{\infty}+[f]_{1}$,

For a function $\ph\,:\,X\to \R\cup\{-\infty\}$, the {\it topological pressure} $P(\ph)$ is defined by
\alil
{
P(\ph)=\lim_{n\to \infty}\frac{1}{n}\log \sum_{w\in S^{n}\,:\,[w]\neq \emptyset}\exp(\sup_{\om\in [w]}S_{n}\ph(\om)),\label{eq:toppres}
}
where we put $S_{n}\ph(\om):=\sum_{k=0}^{n-1}\ph(\si^{k}\om)$. It is known that if $\ph$ is summable, then $P(\ph)$ exists in $[-\infty,+\infty)$.

A $\si$-invariant Borel probability measure $\mu$ on $X$ is called a {\it Gibbs measure} for the potential $\ph\,:\,X\to \R$ if there exist constants $c\geq 1$ and $P\in \R$ such that for any $\om\in X$ and $n\geq 1$
\ali
{
c^{-1}\leq \frac{\mu([\om_{0}\om_{1}\dots \om_{n-1}])}{\exp(-nP+S_{n}\ph(\om))}\leq c.
}
We recall a necessary and sufficient condition for the existence of a Gibbs measure:
\thms
[{\cite{MU}}]\label{th:ex_Gibbs}
Let $X$ be a topological Markov shift whose transition matrix is irreducible. Assume that $\ph\,:\,X\to \R$ is summable. Then $\ph$ admits a Gibbs measure if and only if $A$ is finitely irreducible.
\thme
The theorem tells us that if $A$ is not finitely irreducible, then $\ph$ does not admit a Gibbs measure. On the other hand, a Ruelle-Perron-Frobenius measure for $\ph$ still exists (see Theorem \ref{th:RPF_transitive_summ_anyA2}).

Let $S_{0}\subset S$ be a nonempty subset, and let $M=(M(ij))$ be a matrix indexed by $S_{0}$ such that $M(ij)\leq A(ij)$ for any $i,j\in S_{0}$. Let $\ph\,:\,X\to \R$ be a summable potential satisfying 
$[\ph]_{2}<+\infty$.
We define the transfer operator associated with the subsystem $X_{M}$ of $X$ as follows: The operator $\LR_{M,\ph}$ associated with $M$ and $\ph$ is defined by
\alil
{
\LR_{M,\ph} f(\om)=\sum_{a\in S\,:\,M(a\om_{0})=1}e^{\ph(a\cd \om)}f(a\cd\om), \label{eq:transfer}
}
where $a\cd \om=a\om_{0}\om_{1}\cdots$ is the concatenation of $a$ and $\om$. If $i\om_{0}\notin S_{0}\times S_{0}$, then $M(i\om_{0})$ is understood to be zero. Such an operator is considered by \cite{Demers_,T2020, T2019, T2009}, where perturbations of Ruelle operators with symbol deletions were studied.
These operators are also used to study systems with holes. Indeed, let $\Si=\bigcup_{ij\,:\,M(ij)=1}[ij]$. Then we regard the map $\si|_{\Si}\,:\,\Si\to X$ as an {\it open system} and the map $\si\,:\,X\to X$ as a {\it closed system}. Letting $\Si^{n}:=\bigcap_{k=0}^{n}\si^{-k}\Si$, the $n$-th iterate $(\si|_{\Si})^{n}$ is a map from $\Si^{n}$ to $X$. It os also known that $\bigcap_{n=0}^{\infty}\Si^{n}=X_{M}$. Thus, We may regard $\LR_{M,\ph}$ as a transfer operator associated with the open system $\si|_{\Si}$. When $M=A$, the operator $\LR_{A,\ph}$ is commonly referred to as the {\it Ruelle operator} for the potential $\ph$. 

Assume that conditions (A.1)-(A.3) are satisfied. 
Note that the operator $\LR_{M,\ph}$ acting on $F_{b}(X)$ and on $C_{b}(X)$ are both bounded.
When $M=A$, we write $\LR_{\ph}=\LR_{A,\ph}$ for simplicity.
We recall some useful results of the operator $\LR_{\ph}$.
For $c>0$, let
\alil
{
\Lambda_{c}=\Lambda_{c}(X)=\{f\in C(X)\,:\,0\leq f,\ f(\om)\leq e^{cd_{\theta}(\om,\up)}f(\up) \text{ if }\om\in [\nu_{0}]\}.\label{eq:Lambda_c=...}
}
Let $c_{\adl{Lamc}}=[\ph]_{1}\theta/(1-\theta)$. We recall the following spectral gap property of general Ruelle operators with holes:
\thms
[\cite{T2024_sum}]
\label{th:RPF_transitive_summ_anyA2}
Let $A$ be a zero-one matrix indexed by countable set of states $S$, and let $M$ be a zero-one matrix indexed by a subset $S_{0}\subset S$ such that $M(ij)=1$ yields $A(ij)=1$, and $X_{M}$ contains at least one periodic point. Let $\ph\,:\,X_{A}\to \R$ be a summable function with $[\ph]_{2}<\infty$. Then:
\ite
{
\item There exists a triplet $(\lam,g,\nu)$ such that (i) $\lam$ is the maximal positive eigenvalue of $\LR_{M,\ph}$ equal to $P(\ph|X_{M})$, (ii) $g\in \Lambda_{c_{\adr{Lamc}}}$ is the corresponding nonnegative eigenfunction with $\|g\|_{\infty}=1$, and (iii) $\nu$ is the corresponding positive eigenvector of the dual with $\nu(1)=1$.
\item if $M$ is irreducible, then $\lam$ is a simple eigenvalue, and $\nu(g)>0$. Letting $h:=g/\nu(g)$, the operator $h\otimes \nu$ is the projection onto the one-dimensional eigenspace of $\lam$. Moreover, $\supp\,h=\Si_{S_{0}}$, and $\nu$ is positive Borel probability measure supported on $X_{M}$. The mesure$\mu:=h\nu$ is a $\si$-invariant Borel probability measure, called the {\it Ruelle-Perron-Frobenius measure} for $\ph|X_{M}$.
}
\thme
When $M$ is irreducible, we call $\lam$ in Theorem \ref{th:RPF_transitive_summ_anyA2}(2) the {\it Perron eigenvalue} of $\LR_{M,\ph}$, $h$ the {\it Perron eigenfunction} of $\LR_{M,\ph}$, and $\nu$ the {\it Perron eigenvector} of the dual $\LR_{M,\ph}^{*}$. The triplet $(\lam,h,\nu)$ is then refereed to as the as the {\it Perron spectral triplet} of the operator $\LR_{M,\ph}$.
\cor
{\label{cor:RPF_transitive_summ_anyA2}
Assume that conditions (A.1)-(A.3) are satisfied. Let $U\in S/\!\!\leftrightarrow$ such that $X_{B[U]}$ contains a periodic point for the shift $\si|_{X_{B[U]}}$. Then the operator $\LR_{B[U],\ph}$ admits a Perron spectral triplet $(\lam(U),h(U,\cd),\nu(U,\cd))$.
}
\subsection{Perron complements of Ruelle operator and Schur-Frobenius factorization}\label{sec:subRuop_SFfactor}
Let $X$ be a topological Markov shift with a countable state space $S$ and a zero-one transition matrix $A$ indexed by $S\times S$. Recall the notation $\Si_{U}$ and $\chi_{U}$ in Section \ref{sec:intro} for a subset $U\subset S$.
Let $(\XR,\|\cd\|)$ be a Banach space satisfying $f\chi_{U}\in \XR$ for any $f\in \XR$ and $U\subset S$ (e.g. $\XR=L_{1}(X), C_{b}(X)$ or $F_{b}(X)$).
Denoted by $\LR(\XR)$ the set of all bounded linear operators acting on $\XR$.
For an operator $\MR\in \LR(\XR)$ and subsets $U_{1},U_{2}\subset S$, we define $\MR_{U_{1}U_{2}}\in \LR(\XR)$ by
\ali
{
\MR_{U_{1}U_{2}}f=\chi_{U_{1}}\MR(\chi_{U_{2}}f)
}
for $f\in \XR$. For the identity operator $\IR\in \LR(\XR)$, we write $\IR$ instead of $\IR_{UU}$, $\IR_{UV}$, $\IR_{VU}$ and $\IR_{VV}$ for simplicity. We represent the block operator matrix of an operator $\MR\in \LR(\XR)$ using index subsets $U,V\subset S$ with $U\cap V=\emptyset$. Let $W=U\cup V$. Then we write
\ali
{
\MR_{WW}=
\MatII{\MR_{UU}&\MR_{UV}\\
\MR_{VU}&\MR_{VV}\\
}.
}
That is, the $(i,j)$-element of the matrix is the restricted operator $\MR_{W(i)W(j)}$ of $\MR$ when $W(1)=U$ and $W(2)=V$. This block operator matrix acts as a standard matrix since $\chi_{U}\chi_{V}=\chi_{V}\chi_{U}=0$. We may also write $\MR_{WW}$ by $\MR_{UU}+\MR_{UV}+\MR_{VU}+\MR_{VV}$.

For an operator $\MR\in \LR(\XR)$, subsets $U,V\subset S$ and a number $\lam\in \C$ belonging in the resolvent set of $\MR_{VV}$,  we define a bounded linear operator $\MR[U,V,\lam]\,:\,\XR\to \XR$ by
\ali
{
\MR[U,V,\lam]=\MR_{UU}+\MR_{UV}(\lam\IR-\MR_{VV})^{-1}\MR_{VU}.
}
We begin with the following useful results.
\prop
{\label{prop:PC_evalprop}
Under the notation for $X,\XR,\LR(\XR)$ above, let $\MR\in \LR(\XR)$, $U,V\subset S$ with $U\neq \emptyset$ and $U\cap V=\emptyset$, and $\lam,\eta\in \C$ with $\lam\neq 0$, $\eta\neq 0$. Assume that $\lam$ is in the resolvent set of $\MR_{VV}$. Let $W=U\cup V$. Then
\ite
{
\item[(1)] A function $f\in \XR$ satisfies the equations $\MR[U,V,\lam](\chi_{U}f)=\eta (\chi_{U}f)$, $f\chi_{V}=((\lam/\eta)(\lam\IR-\MR_{VV})^{-1}\MR_{VU})(\chi_{U}f)$ and $f\chi_{S\setminus W}=0$ if and only if $(\MR_{WU}+(\eta/\lam)\MR_{WV})g=\eta g$. Moreover, $\eta$ is an eigenvalue of $\MR[U,V,\lam]$ if and only if $\eta$ is an eigenvalue of $\MR_{WU}+(\eta/\lam)\MR_{WV}$. In particular, in the case when $W=S$, $\eta$ is a simple eigenvalue of $\MR[U,V,\lam]$ if and only if $\eta$ is an eigenvalue of $\MR_{WU}+(\eta/\lam)\MR_{WV}$.
\item[(2)] A functional $\nu\in \XR^{*}$ satisfies $\MR[U,V,\lam]^{*}\nu(\chi_{U}f)=\eta \nu(\chi_{U}f)$, $\nu(\chi_{V}f)=((\lam/\eta)\MR_{UV}(\lam\IR-\MR_{VV})^{-1})^{*}\nu(f)$ and $\nu(\chi_{S\setminus W}f)=0$ for any $f\in \XR$ if and only if $(\MR_{UW}+(\eta/\lam)\MR_{VW})^{*}\nu=\eta \nu$. Furthermore, $\eta$ is an eigenvalue of $\MR[U,V,\lam]^{*}$ if and only if $\eta$ is an eigenvalue of $(\MR_{UW}+(\eta/\lam)\MR_{VW})^{*}$. In particular, in the case when $W=S$, $\eta$ is a simple eigenvalue of $\MR[U,V,\lam]^{*}$ if and only if $\eta$ is a simple eigenvalue of $(\MR_{UW}+(\eta/\lam)\MR_{VW})^{*}$.
}
}
\pros
(1) First we assume that $f\in \XR$ satisfies $\MR[U,V,\lam](\chi_{U}f)=\eta (\chi_{U}f), \chi_{V}f=(\lam/\eta)(\lam\IR-\MR_{VV})^{-1}\MR_{VU}f$ and $\chi_{S\setminus W}f=0$. Let $t=\eta/\lam$. The form of $f\chi_{V}$ implies 
\ali
{
(\lam\IR-\MR_{VV})(f\chi_{V})=&(\lam/\eta)\MR_{VU}f\\
\lam f\chi_{V}=&\MR_{VV}(f\chi_{V})+(\lam/\eta)\MR_{VU}(f\chi_{U})=(\lam/\eta)(t\MR_{VV}+\MR_{VU})f\\
\eta f\chi_{V}=&\chi_{V}(t\MR_{WV}+\MR_{WU})f.
}
Moreover,
\ali
{
\eta f\chi_{U}=&\MR[U,V,\eta](f\chi_{U})\\
=&\MR_{UU}f+t\MR_{UV}f\\
=&\chi_{U}(\MR_{WU}+t\MR_{WV})f
}
Thus we see the equation $\eta f=(\MR_{WU}+t\MR_{WV})f$.
In order to prove the converse, we next assume $(\MR_{WU}+t\MR_{WV})f=\eta f$. Then
\ali
{
\eta f\chi_{V}=&\MR_{VU}f+(\eta/\lam)\MR_{VV}f\\
\lam f\chi_{V}-\MR_{VV}f=&(\lam/\eta)\MR_{VU}f\\
f\chi_{V}=&(\lam/\eta)(\lam\IR-\MR_{VV})^{-1}\MR_{VU}f.
}
Furthermore
\ali
{
f\chi_{U}=\MR_{UU}f+(\eta/\lam)\MR_{UV}f=\MR[U,V,\lam]f.
}
and $f\chi_{S\setminus W}=\eta^{-1}\chi_{S\setminus W}(\MR_{WU}+t\MR_{WV})f=0$. Thus the converse assertion is valid.
Next, if $\eta$ is an eigenvalue of $\MR[U,V,\lam]$, i.e. $f\neq 0$ and $\MR[U,V,\lam]f=\eta f$, then $f\chi_{S\setminus U}=\eta^{-1}\chi_{S\setminus U}\MR[U,V,\lam]f=0$. Thus $f\chi_{U}\neq 0$. By letting $g=f\chi_{U}+(\lam/\eta)(\lam\IR-\MR_{VV})^{-1}\MR_{VU}f$, we get $(\MR_{WU}+(\lam/\eta)\MR_{WV})g=\eta g$. By $g\chi_{U}\neq 0$, $\eta$ is an eigenvalue of $\MR_{WU}+(\lam/\eta)\MR_{WV}$. Conversely, assume $f\neq 0$ and $(\MR_{WU}+t\MR_{WV})f=\eta f$. Since $f\chi_{U}=0$ implies $f\chi_{V}=(\lam/\eta)(\lam\IR-\MR_{VV})^{-1}\MR_{VU}f=0$ and then $f=0$, we get $f\chi_{U}\neq 0$. Thus $\eta$ is an eigenvalue of $\MR[U,V,\lam]$. Finally, under the case $W=S$, we assume that $\eta$ is a simple eigenvalue of $\MR[U,V,\lam]$ with an eigenvector $g\in \XR$. By the argument above, $(\MR_{WU}+t\MR_{WV})f=\eta f$ implies $\MR[U,V,\lam](f\chi_{U})=\eta f\chi_{U}$ and therefore $f\chi_{U}=c g$ for some $c\in \K$. Then $f\chi_{V}=(\lam/\eta)(\lam\IR-\MR_{VV})^{-1}\MR_{VU}f=c (\lam/\eta)(\lam\IR-\MR_{VV})^{-1}\MR_{VU}g$. By $f=\chi_{U}f+\chi_{V}f$ with $W=S$, $\eta$ is simple. On the other hand, if $\eta$ is a simple eigenvalue of $\MR_{WU}+t\MR_{WV}$ with eigenvector $g$, then $\MR[U,V,\lam]f=\eta f$ yields $f=c\chi_{U} g$ for some $c\in \K$. Hence $\eta$ is also simple for $\MR[U,V,\lam]$.
\item[(2)] Assume $\MR[U,V,\lam]^{*}\nu(f\chi_{U})=\eta \nu(f\chi_{U})$, $\nu(\chi_{V}f)=((\lam/\eta)\MR_{UV}(\lam\IR-\MR_{VV})^{-1})^{*}\nu(f)$ and $\nu(\chi_{S\setminus W}f)=0$. Put $t=\eta/\lam$. Then for $f\in \XR$
\ali
{
\nu(\chi_{V}(\lam \IR-\MR_{VV})f)=(\lam/\eta)\nu(\MR_{UV}f).
}
Therefore
\ali
{
\eta\nu(\chi_{V}f)=\nu(\MR_{UV}f)+(\eta/\lam)\nu(\MR_{VV}f)=\nu((\MR_{UW}+(\eta/\lam)\MR_{VW})(\chi_{V}f)).
}
Moreover, 
\ali
{
\eta\nu(\chi_{U}f)=&\nu(\MR_{UU}f+\MR_{UV}(\lam\IR-\MR_{VV})^{-1}\MR_{VU}f)\\
=&\nu(\MR_{UU}f)+(\eta/\lam)\nu(\MR_{VU}f)\\
=&\nu((\MR_{UW}+(\eta/\lam)\MR_{VW})(\chi_{U}f)).
}
Thus we obtain $\eta\nu(f)=\nu((\MR_{UW}+(\eta/\lam)\MR_{VW})f)$. Conversely, we assume $(\MR_{UW}+(\eta/\lam)\MR_{VW})^{*}\nu=\eta \nu$. We notice
\ali
{
\eta\nu(\chi_{V}g)=&\nu(\MR_{UV}g+(\eta/\lam)\MR_{VV}g)\\
\lam\nu(\chi_{V}g)-\nu(\MR_{VV}g)=&(\lam/\eta)\nu(\MR_{UV}g).
}
By letting $g:=(\lam\IR-\MR_{VV})^{-1}f$, we obtain $\nu(\chi_{V}f)=\nu((\lam/\eta)\MR_{UV}(\lam\IR-\MR_{VV})^{-1}f)$. Finally $\nu(\chi_{S\setminus W}f)=\eta^{-1}\nu((\MR_{UW}+(\eta/\lam)\MR_{VW})(\chi_{S\setminus W}f)=0$. Thus the converse assertion is fulfilled.

If $\MR[U,V,\lam]\nu=\eta\nu$ and $\nu\neq 0$, then $\nu(f)=\nu(\chi_{U}f)=\mu(\chi_{U}f)\neq 0$ for some $f\in \XR$. Thus $\mu$ becomes an eigenvector of $\eta$ of $(\MR_{UW}+(\eta/\lam)\MR_{VW})^{*}$. On the other hand, we assume $(\MR_{UW}+(\eta/\lam)\MR_{VW})^{*}\mu=\eta\mu$ and $\mu\neq 0$. By the form $\mu(\chi_{V}f)=\mu((\chi_{U}+(\lam/\eta)\MR_{UV}(\lam\IR-\MR_{VV})^{-1})(\chi_{U}f))$, $\mu(\chi_{U}f)=0$ implies $\mu(\chi_{V}f)=0$ and therefore $\mu(f)=0$. Thus $\mu(\chi_{U}f)\neq 0$ for some $f\in \XR$. This means that the functional $f\mapsto \mu(\chi_{U}f)$ is an eigenvector of $\eta$ of the operator $\MR[U,V,\lam]^{*}$. The proof for simplicity is implied by a similar argument above (1). Hence the assertion is valid.
\proe
\thm
{[A modified Schur-Frobenius factorization]
\label{th:PC_FS}
Under the notation $X,\XR$ and $\LR(\XR)$ in the above, let $\MR\in \LR(\XR)$, $U,V\subset S$ with $U\neq \emptyset$ and $U\cap V=\emptyset$, and $\lam,\eta\in \C$ with $\lam\neq 0$ and $\eta\neq 0$. Assume also that $\lam$ is in the resolvent set of $\MR_{VV}$. Put $W=U\cup V$. Then the operator $\MR_{\eta}:=\MatII{\MR_{UU}&\eta\lam^{-1}\MR_{UV}\\\MR_{VU}&\eta\lam^{-1}\MR_{VV}}$ satisfies
\alil
{
&\MR_{\eta}-\eta\IR\label{eq:M-eta=}\\
=&
\MatII
{
\IR&\MR_{UV}\RR_{\lam}\\
\OR&\IR
}
\MatII
{
\MR[U,V,\lam]-\eta\IR&\OR\\
\OR&\eta\lam^{-1}(\MR_{VV}-\lam\IR)
}
\MatII
{
\IR&\OR\\
\lam\eta^{-1}\RR_{\lam}\MR_{VU}&\IR\\
},\nonumber
}
where we put $\RR_{\lam}=(\MR_{VV}-\lam\IR)^{-1}$. Consequently, the number $\eta$ is in resolvent set of $\MR_{\eta}$ if and only if $\eta$ is in resolvent set of $\MR[U,V,\lam]$.
}
\pros
First we show the equation (\ref{eq:M-eta=}). Denoted by $\MR_{\eta}=\MatII{\AR&\BR\\\CR&\DR}$ with $\AR=\MR_{UU}$, $\BR=\eta\lam^{-1}\MR_{UV}$, $\CR=\MR_{VU}$ and $\DR=\eta\lam^{-1}\MR_{VV}$. By noticing $(\DR-\eta\IR)^{-1}=\lam\eta^{-1}\RR_{\lam}$, $\eta$ is in the resolvent set of $\DR$. The original Schur-Frobenius factorization \cite{AJS} says
\ali
{
&\MR_{\eta}-\eta\IR\\
=&\MatII{
\IR&\BR(\DR-\eta\IR)^{-1}\\\OR&\IR
}
\MatII{
\AR+\BR(\DR-\eta\IR)^{-1}\CR-\eta\IR&\OR\\
\OR&\DR-\eta\IR
}
\MatII{
\IR&\OR\\
(\DR-\eta\IR)^{-1}\CR&\IR\\
}.
}
By using the equation $(\DR-\eta\IR)^{-1}=\eta^{-1}\lam(\MR_{VV}-\lam\IR)^{-1}$, it is not hard to check that the above equation implies the equation (\ref{eq:M-eta=}).
\smallskip
\par
Next we prove the final assertion. We assume that $(\MR[U,V,\lam]-\eta\IR)^{-1}$ exists as bounded linear operator. Then the right hand side of the equation (\ref{eq:M-eta=}) has an inverse bounded linear operator and so does for $\MR_{\eta}-\eta\IR$. On the other hand, $(\MR_{\eta}-\eta\IR)^{-1}$ exists as bounded linear operator. We denote by 
\ali
{
\MatII{\NR_{11}&\NR_{12}\\\NR_{21}&\NR_{22}}:=
\MatII
{
\IR&-\MR_{UV}\RR_{\lam}\\
\OR&\IR
}(\MR_{\eta}-\eta\IR)
\MatII
{
\IR&\OR\\
-\lam\eta^{-1}\RR_{\lam}\MR_{VU}&\IR\\
}.
}
Then (\ref{eq:M-eta=}) implies that $(\MR[U,V,\lam]-\eta\IR)\NR_{11}=\NR_{11}(\MR[U,V,\lam]-\eta\IR)=\IR_{UU}$. Thus $(\MR[U,V,\lam]-\eta\IR)^{-1}$ exists as the operator $\NR_{11}-\eta^{-1}\IR_{VV}+\IR_{YY}$ with $Y=S\setminus W$. Hence the proof is complete.
\proe
We give an application of Theorem \ref{th:PC_FS}.
\thm
{\label{th:PC_specgap}
Let $\MR\in \LR(\XR)$, $U,V\subset S$ with $U\neq \emptyset$ and $U\cap V=\emptyset$, and $\lam,\eta\in \C$ with $\lam\neq 0$ and $\eta\neq 0$. Put $W=U\cup V$. Assume that $\MR_{WW}\,:\,\XR\to \XR$ has the spectral decomposition $\MR_{WW}=\lam (h\otimes \nu)+\RR$ satisfying that $\lam$ is a simple eigenvalue of this operator, $h\in \XR$ is the eigenfunction of $\lam$, $\nu\in \XR^{*}$ is a eigenvector of the dual $\MR_{WW}^{*}$ of $\MR_{WW}$ with $\nu(h)=1$, and the spectrum of $\RR$ is the outside of a neighborhood of $\lam$. Assume also $\nu(\MR_{WU}h)\neq 0$. Then $\MR[U,V,\lam]$ has a spectral gap at the simple eigenvalue $\lam$. 
}
\pros
By the assumption, we have the spectral decomposition $\MR_{WW}=\lam (h\otimes \nu)+\RR$, where $h\in \XR$ is the eigenfunction of $\lam$ and $\nu\in \XR^{*}$ is a eigenvector of the dual $\MR_{WW}^{*}$ of $\MR_{WW}$ with $\nu(h)=1$. there exists a small $\e_{0}>0$ such that for any $\eta\in B_{\C}(\lam,\e_{0})$ with $\eta\neq \lam$, $\eta$ is in the both resolvent sets of $\MR_{WW}$ and of $\MR_{VV}$. Proposition \ref{prop:PC_evalprop}(1) implies that $\lam$ is also a simple eigenvalue of $\MR[U,V,\lam]$. Notice that $\eta$ is an eigenvalue of $\eta\lam^{-1}\LR_{UU}=\eta\lam^{-1}\MR_{\lam}$ and the corresponding eigenfunction is $h$. Since the map 
\ali
{
\C\ni \eta\mapsto \MR_{\eta}=\MR_{WW}+
(\eta\lam^{-1}-1)\MR_{WV}
}
is analytic on $\LR(\XR)$, The analytic perturbation implies that (i) $\MR_{\eta}$ has the spectral decomposition $\MR_{\eta}=\lam(\eta)(h(\eta,\cd)\otimes \nu(\eta,\cd))+\RR_{\eta}$, (ii) $\lam(\eta)$ is a simple eigenvalue, $h(\eta,\cd)$ is the corresponding eigenfunction and $\nu(\eta,\cd)$ is the corresponding eigenvector of the dual of this operator with $\nu(\eta,h(\eta,\cd))=1$, (iii) the maps $\eta\mapsto \lam(\eta)$, $\eta\mapsto h(\eta,\cd)\in \XR$, and $\eta\mapsto \nu(\eta,\cd) \in \XR^{*}$ are analytic, (iv) $\lam(\lam)=\lam$, $h(\lam,\cd)=h$ and $\nu(\lam,\cd)=\nu$, (v) the spectrum of $\RR_{\eta}$ is the outside of a neighborhood $B_{\C}(\lam(\eta),\e_{1})$ of $\lam(\eta)$ for any $\eta\in B_{\C}(\lam,\e_{2})$ for some $0<\e_{1}$ and $0<\e_{2}$. Since $\nu(h(\eta,\cd))\to \nu(h)=1$ and $\nu(\eta,\MR_{WU}h)\to \nu(\MR_{WU}h)\neq 0$ as $\eta\to \lam$, we may assume $\nu(h(\eta,\cd))\neq 0$ and $\nu(\eta,\MR_{WU}h)\neq 0$ for any $\eta\in B_{\C}(\lam,\e_{2})$. Similarity, for any such a number $\eta$, the number $\eta$ is in $B_{\C}(\lam(\eta),\e_{1})$ from $\lam(\eta)\to \lam$ as $\eta\to \lam$.

Choose any $\eta\in B_{\C}(\lam,\min(\e_{0},\e_{1},\e_{2}))$ with $\eta\neq \lam$. Now we show $\eta$ is in the resolvent set of the operator $\MR_{\eta}$. By using $\MR_{\eta}^{*} \nu(\eta,\cd)=\lam(\eta)\nu(\eta,\cd)$ and $\eta\lam^{-1}\MR_{WW}h=\eta h$
\ali
{
(\lam(\eta)-\eta)\nu(\eta,h)=\nu(\eta,(\MR_{\eta}-\eta\lam^{-1}\MR_{WW})h)=(1-\eta\lam^{-1})\nu(\eta,\MR_{WU}h).
}
By $\eta\neq \lam$, $\nu(\eta,h)\neq 0$ and $\nu(\eta,\MR_{WU}h)\neq 0$, we get $\lam(\eta)-\eta\neq 0$. In addition to the fact $\eta\in B_{\C}(\lam(\eta),\e_{1})$, we obtain that $\eta$ is in the resolvent set of $\MR_{\eta}$. Hence $\eta$ is in the resolvent set of $\MR[U,V,\lam]$ and we see this operator has a spectral gap at $\lam$.
\proe
We mainly consider Perron complements of Ruelle type operators $\LR_{M,\ph}$ acting on $F_{b}(X)$ under conditions (A.1)-(A.3). Let $U,V\subset S$ with $U\cap V=\emptyset$ and $U\neq \emptyset$. Take $\eta\in \C$ so that $r(\LR_{VV})<|\eta|$. The resolvent $(\eta\IR-\LR_{VV})^{-1}$ has the series form
\ali
{
(\eta\IR-\LR_{VV})^{-1}=\sum_{n=0}^{\infty}\frac{\LR_{VV}^{n}}{\eta^{n+1}}.
}
Therefore, we have for $f\in C_{b}(X)$ and for $\om\in \Si_{U}$
\ali
{
\LR[U,V,\eta]f(\om)=\sum_{w\in W(U:V)\,:\,w\cd \om_{0}\in W(A)}\eta^{-|w|+1}e^{S_{|w|}\ph(w\cd\om)}f(w\cd\om),
}
where $|w|$ means the length of the word $w$ and
\ali
{
W(U:V):=U\times \bigcup_{n=0}^{\infty}V^{n}
}
and $\bigcup_{n=0}^{\infty}V^{n}$ includes an empty word.
\thm
{\label{th:spec_PC}
Let $A$ be a zero-one matrix indexed by countable states $S$ and $M$ a zero-one matrix indexed by $S_{0}\subset S$ such that $M$ is irreducible,  $M(ij)=1$ yields $A(ij)=1$ and $X_{M}$ has at least one periodic point. Let $\ph\,:\,X_{A}\to \R$ be a function with $[\ph]_{2}<\infty$ and summable. Assume also that $M$ is irreducible. Choose any $U\subset S_{0}$ with $U\neq S_{0}$ and let $V=S_{0}\setminus U$. Then $\LR_{M,\ph}[U,V,\lam]\,:\,F_{b}(X)\to F_{b}(X)$ has the spectral gap at the spectral radius $\lam$ of $\LR_{M,\ph}$.
}
\pros
Take $(\eta,h,\nu)$ the Perron spectral triplet of the operator $\LR_{M,\ph}$. By the positively of $h$ and $\nu$, and by the support of $\supp\,\nu=X_{M}$ and $\supp\,h=\Si_{S_{0}}$, we see $\nu((\LR_{M,\ph})_{S_{0}U}h)=\nu(\LR_{M,\ph}(\chi_{U}h))=\lam\nu(\chi_{U}h)>0$. Hence Theorem \ref{th:PC_specgap} yields the assertion.
\proe
\subsection{Convergence of the spectral triplets of Ruelle type operators}\label{sec:TSC}
In this section, we will give convergence of Perron spectral triplet of transfer operators.
Let $X$ be a topological Markov shift with transition matrix $A$ and countable state space $S$. Assume that conditions (A.1)-(A.3) are satisfied. We use all notation in Section \ref{sec:intro} and in Section \ref{sec:prelim}. For Let $U\in S/\!\!\leftrightarrow$ satisfying that $X_{B[U]}$ has a periodic point for $\si|_{X_{B[U]}}$. denoted by 
\ali
{
(\lam(U),h(U,\cd),\nu(U,\cd))
}
the Perron spectral triplet of $\LR_{B[U],\ph}$ which appears in Corollary \ref{cor:RPF_transitive_summ_anyA2}. Note that if there is no periodic point of $X_{B[U]}$, the spectral radius of $\LR_{B[U],\ph}$ is equal to $0$ (see \cite[Proposition 4.2]{T2024_sum}) and therefore we may put $\lam(U):=0$. For $U\in \mathscr{S}_{0}$ and a nonempty subset $Q\subset U$, we introduce restricted eigenfunctions and conditional eigenvectors as follows:
\ali
{
\nu^{P}(U,f):=\nu(U,f|\Si_{P})=\frac{\nu(U,f\chi_{P})}{\nu(U,\chi_{P})},\quad 
h^{P}(U,\cd):=\frac{h(U,\cd)}{\nu^{P}(U,h(U,\cd)\chi_{P})}.
}
By virtue of Proposition \ref{prop:PC_evalprop}(1), we notice that $\LR[P,U\setminus P,\lam]h^{P}(U,\cd)=\lam h^{P}(U,\cd)$ and $\LR[P,U\setminus P,\lam]^{*}\nu^{P}(U,\cd)=\lam\nu^{P}(U,\cd)$. Therefore, we may call the triplet $(\lam,h^{P}(U,\cd),\nu^{P}(U,\cd))$ the Perron spectral triplet of $\LR[P,U\setminus P,\lam]$. 

For each $\e\in (0,1)$, we write by
\ali
{
(\lam_{\e},h_{\e},\nu_{\e})
}
the Perron spectral triplet of the operator $\LR_{A,\phe}$ and
\ali
{
\nu_{\e}^{P}(f):=\nu_{\e}(f|\Si_{P}),\quad 
h_{\e}^{P}:=\frac{h_{\e}\chi_{P}}{\|h_{\e}\chi_{P}\|_{\infty}}.
}
We see that $(\lam_{\e},h_{\e}^{P}(U,\cd),\nu_{\e}^{P}(U,\cd))$ becomes the Perron spectral triplet of $\LR_{\e}[P,S\setminus P,\lam_{\e}]$. Further, we put
\ali
{
&g(U,\cd):=h(U,\cd)/\|h(U,\cd)\|_{\infty},\quad g_{\e}:=h_{\e}/\|h_{\e}\|_{\infty}\\
&g^{P}(U,\cd):=h^{P}(U,\cd)/\|h^{P}(U,\cd)\|_{\infty},\quad g_{\e}^{P}:=h_{\e}^{P}/\|h_{\e}^{P}\|_{\infty}\\
&\mu(U,\cd):=h(U,\cd)\nu(U,\cd),\quad \mu_{\e}:=h_{\e}\nu_{\e}.
}
The notation $\mu^{P}(U,\cd)$ and $\mu_{\e}^{P}$ are similar defined. Put $c_{\adl{cL}}=\max\{[\ph]_{2},\sup_{\e>0}[\phe]_{2}\}\theta/(1-\theta)$. We note that all eigenfunctions given in the above are in $\Lambda_{c_{\adr{cL}}}$. For simplicity, we denote
\ali
{
\lam:=\max_{U\in \mathscr{S}_{0}}\lam(U),\quad \LR:=\LR_{B,\ph},\quad \LR_{\e}:=\LR_{A,\phe}.
}
We begin with the following basic properties of $\Lambda_{c}$:
\prop
{\label{prop:prop_Lamc0}
\ite
{
\item For any $a\in S$, $\chi_{[a]}$ is in $\Lambda_{c}$ for all $c\geq 0$.
\item For $f,g\in \Lambda_{c}$, $f+g\in \Lambda_{c}$.
\item For $f\in \Lambda_{c}$ and $a\geq 0$, $af\in \Lambda_{c}$.
\item For $f\in \Lambda_{c}$ with $\LR_{A_{0},\ph}f\in C(X)$ and $c\geq [\ph]_{2}\theta/(1-\theta)$, $\LR_{A_{0},\ph}f \in \Lambda_{c}$.
}
}
\pros
We will only show the assertion (4). For $f\in \Lambda_{c}$ and $\om_{0}=\up_{0}$
\ali
{
\LR_{A_{0},\ph} f(\om)=&\sum_{a\in S\,:\,A_{0}(a\om_{0})=1}e^{\ph(a\cd \om)}f(a\cd\om)\\
\leq&\sum_{a\in S\,:\,A_{0}(a\up_{0})=1}e^{\ph(a\cd\up)+[\ph]_{2}d_{\theta}(a\cd \om,a\cd\up)}f(a\cd\up)e^{cd_{\theta}(a\cd\om,a\cd\up)}\\
=&\LR_{A_{0},\ph}f(\up) e^{\theta([\ph]_{2}+c)d_{\theta}(\om,\up)}.
}
Here we note $\theta([\ph]_{2}+c)\leq c$ $\iff$ $\theta [\ph]_{2}/(1-\theta)\leq c$. Thus we get the inequality $\LR_{A_{0},\ph} f(\om)\leq \LR_{A_{0},\ph}f(\up) e^{c d_{\theta}(\om,\up)}$. Hence $\LR_{A_{0},\ph} f\in \Lambda_{c}$.
\proe
\prop
{\label{prop:B2=>sum}
Assume that conditions (A.1)-(A.3) are satisfied. Then there exists a constant $c_{\adl{cbd}}\geq 1$ such that for any $n\geq 1$ and $\om,\up\in X$ with $\om\in [\up_{0}\cdots \up_{n-1}]$, we have
$e^{S_{n}\ph(\e,\om)}\leq c_{\adr{cbd}}e^{S_{n}\ph(\e,\up)}$.
}
\pros
\item By $c_{\adr{cL}}=\sup_{\e>0}[\phe]_{2}<+\infty$, we have
\ali
{
e^{S_{n}\ph(\e,\om)}=&e^{\sum_{i=0}^{n-1}\ph(\e,\si^{i}\om)}\\
\leq&e^{\sum_{i=0}^{n-1}c_{\adr{cL}}d_{\theta}(\si^{i}\om,\si^{i}\up)+\sum_{i=0}^{n-1}\ph(\e,\si^{i}\up)}\\
\leq&e^{c_{\adr{cL}}\sum_{i=0}^{n-1}\theta^{n+1-i}}e^{S_{n}\ph(\e,\up)}\\
\leq&e^{c_{\adr{cL}}\theta^{2}/(1-\theta)}e^{S_{n}\ph(\e,\up)}.
}
Hence the assertion is valid by putting $c_{\adr{cbd}}:=e^{c_{\adr{cL}}\theta^{2}/(1-\theta)}$.
\proe
\prop
{[\cite{T2020}]\label{prop:conv_op}
Assume that conditions (A.1)-(A.3) are satisfied. Then
\ite
{
\item $\lam_{\e}$ converges to $\lam$ as $\e\to 0$;
\item $\|\LR_{\e}-\LR\|_{\infty}\to 0$ as $\e\to 0$. Therefore, we obtain $\|(\LR_{\e})_{UV}-\LR_{UV}\|_{\infty}\to 0$ for any $U,V\subset S$.
}
}
\pros
Note that if (A.1)-(A.3) are satisfied, then conditions (A.1)-(A.3) and (B.1)-(B.3) given in \cite{T2024_sum} are satisfied. Then the assertion (1) is due to \cite[Proposition 6.4]{T2024_sum}. The assertion (2) follows from \cite[Proposition 4.3]{T2024_sum}.
\proe
\prop
{\label{prop:conv_PC_Rue1_v2}
Assume that conditions (A.1)-(A.3) are satisfied. Let $U\subset S$ be a nonempty subset and $V\subset S\setminus U$ with $r(\LR_{VV})<\lam$. Then $\|\LR_{\e}[U,V,\lam_{\e}]-\LR[U,V,\lam]\|_{\infty}\to 0$ as $\e\to 0$. 
}
\pros
Recall equations
\ali
{
\LR[U,V,\lam]=&\LR_{UU}+\LR_{UV}(\lam-\LR_{VV})^{-1}\LR_{VU}\\
\LR_{\e}[U,V,\lam_{\e}]=&(\LR_{\e})_{UU}+(\LR_{\e})_{UV}(\lam_{\e}-(\LR_{\e})_{VV})^{-1}(\LR_{\e})_{VU}.
}
By Proposition \ref{prop:conv_op}(2), we see $\|(\LR_{\e})_{UV}-\LR_{UV}\|_{\infty}\to 0$ for any $U,V\subset S$. Therefore, it suffices to show that $\|(\lam-\LR_{VV})^{-1}-(\lam_{\e}-(\LR_{\e})_{VV})^{-1}\|_{\infty}\to 0$. By the resolvent equation
\ali
{
&(\lam_{\e}-(\LR_{\e})_{VV})^{-1}-(\lam-\LR_{VV})^{-1}\\
=&(\lam_{\e}-(\LR_{\e})_{VV})^{-1}(\lam-\LR_{VV}-\lam_{\e}+(\LR_{\e})_{VV})(\lam-\LR_{VV})^{-1}.
}
Since $\|(\LR_{\e})_{VV}-\LR_{VV}\|_{\infty}\to 0$, $r((\LR_{\e})_{VV})\to r(\LR_{VV})<\lam$ and $\lam_{\e}\to \lam$, we obtain that for any small $\eta_{2}>\eta_{1}>0$ with $r(\LR_{VV})<\eta_{1}<\eta_{2}<\lam$, there exists $\e_{0}>0$ such that for any $0<\e<\e_{0}$, $r((\LR_{\e})_{VV})<\eta_{1}<\eta_{2}<\lam_{\e}$. Moreover, there exists a constant $c_{\adl{cop}}>0$ for any $n\geq 1$, 
\ali
{
\|(\LR_{\e})_{VV}^{n}\|_{\infty}\leq c_{\adr{cop}}\eta_{1}^{n}.
}
Then
\ali
{
\|(\lam_{\e}-(\LR_{\e})_{VV})^{-1}\|_{\infty}=\left\|\sum_{n=0}^{\infty}\frac{(\LR_{\e})_{VV}^{n}}{\lam_{\e}^{n+1}}\right\|_{\infty}\leq\sum_{n=0}^{\infty}\frac{\|(\LR_{\e})_{VV}^{n}\|_{\infty}}{|\lam_{\e}|^{n+1}}\leq\sum_{n=0}^{\infty}\frac{c_{\adr{cop}}}{\eta_{2}}\left(\frac{\eta_{1}}{\eta_{2}}\right)^{n}=c_{\adr{cod}}<\infty
}
letting $c_{\adl{cod}}=c_{\adr{cop}}/(\eta_{2}-\eta_{1})$.
Hence we obtain
\ali
{
&\|\LR_{\e}[U,V,\lam_{\e}]-\LR[U,V,\lam]\|_{\infty}\\
\leq&\|(\LR_{\e})_{UU}-\LR_{UU}\|_{\infty}+\|(\LR_{\e})_{UV}-\LR_{UV}\|_{\infty}\|(\lam_{\e}-(\LR_{\e})_{VV})^{-1}(\LR_{\e})_{VU}\|_{\infty}\\
&+\|\LR_{UV}\|_{\infty}\|(\lam_{\e}-(\LR_{\e})_{VV})^{-1}-(\lam-\LR_{VV})^{-1}\|_{\infty}\|(\LR_{\e})_{VU}\|_{\infty}\\
&+\|\LR_{UV}(\lam-\LR_{VV})^{-1}\|_{\infty}\|(\LR_{\e})_{VU}-\LR_{VU}\|_{\infty}\\
\leq&\|(\LR_{\e})_{UU}-\LR_{UU}\|_{\infty}+c_{\adr{cod}}\|(\LR_{\e})_{UV}-\LR_{UV}\|_{\infty}\\
&+\|\LR_{UV}\|_{\infty}c_{\adr{cod}}(|\lam_{\e}-\lam|+\|(\LR_{\e})_{VV}-\LR_{VV}\|_{\infty})\|(\lam-\LR_{VV})^{-1}\|_{\infty}\\
&+\|\LR_{UV}(\lam-\LR_{VV})^{-1}\|_{\infty}\|(\LR_{\e})_{VU}-\LR_{VU}\|_{\infty}\\
\to&0.
}
\proe
\prop
{\label{prop:gP_nonzero}
Assume that conditions (A.1)-(A.3) are satisfied. Take $U\in S/\!\!\leftrightarrow$. Assume also that a nonempty subset $P\subset U$ satisfies that there exists a finite subset $P_{0}\subset P$ such that $\{i\in S\,:\,A(ij)=1 \text{ for some }j\in P\}=\{i\in S\,:\,A(ij)=1 \text{ for some }j\in P_{0}\}$. Then $\liminf_{\e\to 0}g_{\e}^{P}\neq 0$. 
}
\rem
{
If either $P$ is finite or a transition matrix given in application in directed-graph (see the proof of Proposition \ref{prop:positivebdd_ge^E}), then the above additional condition is satisfied.
}
\pros
For any $\om\in \Si_{P}$, we have
\ali
{
g_{\e}^{P}(\om)=&\lam_{\e}^{-1}\LR_{\e}[P,S\setminus P,\lam_{\e}]g_{\e}^{P}(\om)\\
=&\lam_{\e}^{-1}\sum_{w\in W(P,S\setminus P)\,:\,w\cd\om_{0}\in W(A)}\lam_{\e}^{-|w|+1}\exp(S_{|w|}\ph(\e,w\cd\om))g_{\e}^{P}(w\cd\om)\\
\leq &c_{\adr{cbd}}e^{c_{\adr{cL}}\theta}\lam_{\e}^{-1}\sum_{j\in P_{0}}\sum_{w\in W(P,S\setminus P)\,:\,w\cd\om_{0}\in W(A)}\lam_{\e}^{-|w|+1}\exp(S_{|w|}\ph(\e,w\cd\om^{j}))g_{\e}^{P}(w\cd\om^{j})\\
=&c_{\adr{cbd}}e^{c_{\adr{cL}}\theta}\sum_{j\in P_{0}}\lam_{\e}^{-1}\LR_{\e}[P,S\setminus P,\lam_{\e}]g_{\e}^{P}(\om^{j})
=c_{\adr{cbd}}e^{c_{\adr{cL}}\theta}\sum_{j\in P_{0}}g_{\e}^{P}(\om^{j}).
}
Since $\|g_{\e}^{P}\|_{\infty}=1$, we see
\ali
{
c_{\adr{cbd}}^{-1}e^{-c_{\adr{cL}}\theta}\leq \sum_{j\in P_{0}}g_{\e}^{P}(\om^{j})
}
for any $\e\in (0,1)$. Now we take a positive sequence $(\e(n))$ with $\e(n)\to 0$ as $n\to \infty$ so that for any $j\in P_{0}$, $\liminf_{\e\to 0}g_{\e}^{P}(\om^{j})=\lim_{n\to \infty}g_{\e}^{P}(\om_{j})$. By Ascoli Theorem, there exist a subsequence $\{\e^\p(n)\}$ of $\{\e(n)\}$ and $g_{0}\in \Lambda_{c_{\adr{cL}}}$ such that $g_{\e^\p(n)}^{P}(\om)\to g_{0}(\om)$ for each $\om\in X$. Since $c_{\adr{cbd}}^{-1}e^{-c_{\adr{cL}}\theta}\leq \sum_{j\in P_{0}}g_{0}(\om^{j})$, we have $g_{0}(\om^{j})>0$ for some $j\in P_{0}$. Hence we obtain $\liminf_{\e\to 0}g_{\e}^{P}(\om^{j})\neq 0$ for some $j\in P_{0}$.
\proe
\prop
{\label{prop:gP_nonzero_inf}
Assume that conditions (A.1)-(A.3) are satisfied. Take $U\in S/\!\!\leftrightarrow$ and a nonempty subset $P\subset U$ . Assume also that $\liminf_{\e\to 0}g_{\e}^{P}\neq 0$. Then $g_{\e}^{P}(\om)/\nu_{\e}^{P}(U,g_{\e}^{P})\to h^{P}(U,\om)$ for each $\om\in \Si_{P}\cap X_{U}$.
}
\pros
Choose any positive sequence $(\e(n))$ with $\inf_{n}\e(n)=0$. 
Recall that $g_{\e}^{P}$ is the corresponding eigenfunction of the eigenvalue $\lam_{\e}$ of the operator $\LR_{\e}[P,S\setminus P,\lam_{\e}]$ with $\|g_{\e}^{P}\|_{\infty}=1$. Since $g_{\e}^{P}$ is in $\Lambda_{c_{\adr{cL}}}$ for all $\e>0$, Ascoli Theorem implies that there exist a subsequence $(\e^\p(n))$ of $(\e^\p(n))$ and $g_{0}\in \Lambda_{c_{\adr{cL}}}$ such that $g^{P}_{\e^\p(n)}(\om)\to g_{0}(\om)$ as $n\to \infty$ for each $\om\in X$.
By the assumption, $g_{0}\neq 0$. 
\alil
{
\lam_{\e}g_{\e}^{P}=\LR_{\e}[P,S\setminus P,\lam_{\e}]g_{\e}^{P}\geq \LR_{\e}[P,U\setminus P,\lam_{\e}]g_{\e}^{P}.\label{eq:legep<=Le[PUP...]}
}
We show convergence $\LR_{\e}[P,U\setminus P,\lam_{\e}]g_{\e}^{P}(\om)\to \LR[P,U\setminus P,\lam_{\e}]g^{P}(\om)$ for each $\om\in X$. We have
\alil
{
&|\LR_{\e}[P,U\setminus P,\lam_{\e}]g_{\e}^{P}(\om)-\LR[P,U\setminus P,\lam]g_{0}(\om)|\nonumber\\
\leq &\|\LR_{\e}[P,U\setminus P,\lam_{\e}]-\LR[P,U\setminus P,\lam]\|_{\infty}+\LR[P,U\setminus P,\lam]|g_{\e}^{P}-g_{0}|(\om).\label{eq:Lge-Lg0}
}
Now we show $\LR[P,U\setminus P,\lam]|g_{\e}^{P}-g_{0}|(\om)\to 0$ for each $\om$. Fix $\om\in X$. Choose any $\eta>0$. Put $V:=U\setminus P$. By $r(\LR_{VV})<\lam$, there exists $n_{0}\geq 1$ such that
\ali
{
\sum_{n=n_{0}}^{\infty}\frac{\|\LR_{VV}^{n}\|_{\infty}}{\lam^{n+1}}\leq \frac{\eta}{1+\|\LR_{UV}\|_{\infty}\|\LR_{VU}\|_{\infty}}.
}
Then we see
\ali
{
\LR[U,V,\lam]|g_{\e}^{P}-g_{0}|(\om)=&(\LR_{UU}+\LR_{UV}(\lam-\LR_{VV})^{-1}\LR_{VU})|g_{\e}^{P}-g_{0}|(\om)\\
\leq &(\LR_{UU}+\sum_{n=0}^{n_{0}}\lam^{-n-1}\LR_{UV}\LR_{VV}^{n}\LR_{VU})|g_{\e}^{P}-g_{0}|(\om)+2\eta\\
\leq &\sum_{n=1}^{n_{0}+2}\lam^{-n+1}\LR^{n}|g_{\e}^{P}-g_{0}|(\om)+2\eta.
}
Let $c_{\adl{psm}}:=\sum_{a\in S}\sup_{\up\in [a]}e^{\ph(\up)}$. Since $\ph$ is summable, for any $1\leq n\leq n_{0}+1$, there exists a finite subset $F_{n}\subset S$ such that 
\ali
{
\sum_{a\in S\setminus F_{n}}\sup_{\up\in [a]}e^{\ph(\up)}\leq \eta\lam^{n-1}/((n_{0}+2)c_{\adr{psm}}^{n-1}).
}
In addition to the finiteness of $\bigcup_{n=1}^{n_{2}+2}F_{n}$, there exists $m_{0}\geq 1$ such that for any $1\leq n\leq n_{0}+2$, $m\geq m_{0}$, and $w\in F_{n}$ with $w\cd\om_{0}\in W_{n+1}(A)$, 
\ali
{
e^{S_{n}\ph(w\cd\om)}|g_{\e^\p(m)}^{P}(w\cd\om)-g_{0}(w\cd\om)|<\eta\lam^{n-1}/(n_{0}+2).
}
Thus, we have
\ali
{
\sum_{n=1}^{n_{0}+2}\lam^{-n+1}\LR^{n}|g_{\e}^{P}-g_{0}|(\om)\leq&\sum_{n=1}^{n_{0}+2}\lam^{-n+1}\sum_{w\in F_{n}\,:\,w\cd\om_{0}\in W_{n+1}(A)}e^{S_{n}\ph(w\cd\om)}|g_{\e}^{P}(w\cd\om)-g_{0}(w\cd\om)|\\
&+\sum_{n=1}^{n_{0}+2}2\lam^{-n+1}c_{\adr{psm}}^{n-1}\sum_{a\in S\setminus F_{n}}\sup_{\up\in [a]}e^{\ph(\up)}\\
\leq&\eta+2\eta.
}
Hence
\ali
{
|\LR[U,V,\lam]g_{\e^\p(m)}^{P}(\om)-\LR[U,V,\lam]g_{0}(\om)|\leq
\LR[U,V,\lam]|g_{\e(m)}^{P}-g_{0}|(\om)\leq 5\eta
}
for any $m\geq m_{0}$. This means that $\LR[U,V,\lam]g_{\e^\p(m)}^{P}(\om)$ converges to $\LR[U,V,\lam]g_{0}(\om)$.

Together with Proposition \ref{prop:conv_PC_Rue1_v2}, (\ref{eq:Lge-Lg0}) yields convergence of $\LR_{\e}[P,U\setminus P,\lam_{\e}]g_{\e}^{P}(\om)$ to $\LR_{\e}[P,U\setminus P,\lam]g_{0}(\om)$. Letting $\e=\e^\p(m)\to 0$,  (\ref{eq:legep<=Le[PUP...]}) implies
\ali
{
\lam g_{0}(\om)\geq \LR[P,U\setminus P,\lam]g_{0}(\om)
}
for each $\om\in X$. Let $f:=\lam g_{0}-\LR[P,U\setminus P,\lam]g_{0}$. We notice $\nu^{P}(U,f)=0$. Since $f$ is nonnegative and continuous, and since $\nu^{P}(U,\cd)$ is a Borel probability measure, we see $\lam g_{0}-\LR[P,U\setminus P,\lam]g_{0}=f\equiv 0$ on $\supp\,\nu^{P}(U,\cd)=\Si_{P}\cap X_{B(U)}$. By the simplicity of the eigenvalue $\lam$ of the operator $\LR[P,U\setminus P,\lam]$ restricted to the space $F_{b}(X_{U})$, we have $g_{0}=c h^{P}(U,\cd)$ on $X_{U}\cap \Si_{P}$ for some positive $c>0$. 

Now we prove the assertion. By Lebesgue Dominated convergence theorem. we obtain $\nu^{P}(U,h_{\e(n^\p)}^{P})\to c \nu^{P}(U,h^{P}(U,\cd))=c>0$ as $n^\p\to \infty$. Thus $g_{\e}^{P}(\om)/\nu^{P}(U,g_{\e}^{P})\to h^{P}(U,\om)/\nu^{P}(U,h^{P}(U,\cd))=h^{P}(U,\om)$ for each $\om$ runs though the subsequence $(\e^\p(n))$ of $(\e(n))$. This fact does not depend on choosing $(\e(n))$. Hence the proof is complete.
\proe
\prop
{\label{prop:gP_positive_inf}
Assume that conditions (A.1)-(A.3) are satisfied. Take $U\in S/\!\!\leftrightarrow$ and a nonempty subset $P\subset U$. Assume also that 
\ite
{
\item[(i)] $\liminf_{\e\to 0}g_{\e}^{P}\neq 0$.
\item[(ii)] There exists a finite set $Q_{0}\subset U$ such that for any $b\in P$, $B(ab)=1$ for some $a\in Q_{0}$.
\item[(iii)] If $\# P=+\infty$, then $\sup_{\e>0}[\phe]_{1}<+\infty$.
}
Then there exist a constant $c_{\adr{ubge}}>0$ and $\e_{0}>0$ such that $g_{\e}^{P}\geq c_{\adr{ubge}}\chi_{P}$ for any $0<\e<\e_{0}$.
}
\pros
First we consider the case that $B_{UU}$ is $1\times 1$ zero matrix. In this case, we have $P=U$ and $\# P=1$. For any $\om\in \Si_{P}$, we have $\om_{0}=\up_{0}$, $g_{\e}^{P}(\om)\leq e^{c_{\adr{cL}}\theta}g_{\e}^{P}(\up)$. Therefore the assertions holds.

Next we assume that $B_{UU}$ is not zero matrix. Since $B_{UU}$ is irreducible, for any $a,b\in P$, there exist an integer $m(a,b)\geq 1$ and $w(a,b)=w_{1}(a,b)\cd\cdots\cd w_{m(a,b)}(a,b)\in W(P,U\setminus P)^{m(a,b)}$, $w(a,b)\cd b$ is $B$-admissible and $(w_{1}(a,b))_{1}=a$. For any $a,b\in P$, fix $\tau(a,b)\in [w(a,b)\cd b]$. Moreover, since $\lam_{\e}\to \lam$ and $\ph(\e,\om)\to \ph(\om)$ for each $\om\in X$ with $B(\om_{0}\om_{1})=1$, there exists $0<\e_{0}<1$ such that for any $0<\e<\e_{0}$
\ali
{
\frac{\lam}{2}\leq \lam_{\e}\leq 2\lam,\quad 
e^{\ph(\e,\tau(a,b))}\geq \frac{1}{2}e^{\ph(\tau(a,b))}.
}
We also note that $\lam_{\e}^{m}\geq \underline{\lam}^{|m|}$ for any integer $m$ by putting
$\underline{\lam}=\min(2^{-1}\lam^{-1},2^{-1}\lam)$.
\smallskip
\\
Choose any $\om,\up\in \Si_{P}$. Put $a=\up_{0}$ and $b=\om_{0}$. Write $m=m(a,b)$ and $w=w(a,b)$ for simplicity. Then we have
\ali
{
g^{P}_{\e}(\om)=&\lam_{\e}^{-m}\LR_{\e}^{m}[P,S\setminus P,\lam_{\e}]g^{P}_{\e}(\om)\\
\geq &\lam_{\e}^{-m}\LR_{\e}^{m}[P,U\setminus P,\lam_{\e}]g^{P}_{\e}(\om)\\
=&\lam_{\e}^{-m+|w|-1}\exp(S_{|w|}\ph(\e,w\cd \om))g^{P}_{\e}(w\cd\om)\\
\geq &(\underline{\lam})^{|-m+|w|-1|}c_{\adr{cbd}}^{-1}\exp(S_{|w|}\ph(\e,\tau(a,b)))e^{-c_{\adr{cL}}\theta}g^{P}_{\e}(\up)\\
\geq &(\underline{\lam})^{|-m+|w|-1|}c_{\adr{cbd}}^{-1}2^{-|w|}e^{-c_{\adr{cL}}\theta}(e^{S_{|w|}\ph(\tau(a,b))})g^{P}_{\e}(\up).
}
using Proposition \ref{prop:B2=>sum} and the fact $g^{P}_{\e}\in \Lambda_{c_{\adr{cL}}}$.
By $\liminf_{\e\to 0}g_{\e}^{P}\neq 0$, there exists a sequence $(\e(n))$ such that $\lim_{n\to \infty}g_{\e(n)}^{P}(\up)>0$ for some $\up\in \Si_{P}$. Put $a:=\up_{0}$. If $P$ is finite, then by letting $c_{\adl{ingPe}}=\min_{b\in P}((\underline{\lam})^{|-m(a,b)+|w(a,b)|-1|}2^{-|w(a,b)|}e^{S_{|w|}\ph(\tau(a,b))})$, we obtain the inequality $g^{P}_{\e}\geq c_{\adr{ingPe}}^{-1}\chi_{P}$. 

Next we consider the case $\# P=+\infty$. When we put
\ali
{
c_{\adr{ingPe}}=\min_{b\in Q_{0}}((\underline{\lam})^{|-m(a,b)+|w(a,b)|-1|}2^{-|w(a,b)|}e^{S_{|w|}\ph(\tau(a,b))}),
}
we obtain $g^{P}_{\e}(\om)\geq c_{\adr{cbd}}^{-1}g^{P}_{\e}(\up)>0$ for any $\om\in \Si_{Q_{0}}$ and for any $0<\e<\e_{0}$. Fix $\om(b)\in [b]$ for each $b\in Q_{0}$. 
By the assumption, for any $\om\in \Si_{P}$, there exists $b\in Q_{0}$ such that $B(b\om_{0})=1$.
We have 
\ali
{
g^{P}_{\e}(\om)=\frac{g_{\e}(\om)}{\|g_{\e}\chi_{P}\|_{\infty}}=&\lam_{\e}^{-1}\LR_{\e}\frac{g_{\e}(\om)}{\|g_{\e}\chi_{P}\|_{\infty}}\\
=&\lam_{\e}^{-1}\sum_{s\in S\,:\,A(s\om_{0})=1}e^{\ph(\e,s\cd \om)}\frac{g_{\e}(s\cd\om)}{\|g_{\e}\chi_{P}\|_{\infty}}\\
\geq &\lam_{\e}^{-1}\sum_{s\in Q_{0}\,:\,B(s\om_{0})=1}e^{\ph(\e,s\cd \om)}g^{P}_{\e}(s\cd \om)\\
\geq& \underline{\lam}^{-1}e^{-c_{\adr{cL}}\theta-\sup_{\e>0}[\phe]_{1}\theta+\ph(\e,\om(b))}\min_{b\in Q_{0}}g^{P}_{\e}(\om(b))\\
\geq& \underline{\lam}^{-1}e^{-c_{\adr{cL}}\theta-\sup_{\e>0}[\phe]_{1}\theta+\min_{b\in Q_{0}}\ph(\e,\om(b))}c_{\adr{cbd}}^{-1}g^{P}_{\e}(\up)=:c_{\adr{ubge}}>0
}
for any $0<\e<\e_{0}$. Hence the assertion holds.
\proe
\prop
{\label{prop:conv_PC_Rue2}
Assume that conditions (A.1)-(A.3) are satisfied. Let $U\in \mathscr{S}_{0}$ and $P\subset S$ a nonempty subset. We also assume that $g_{\e}^{P}\geq c_{\adr{ubge}}\chi_{P}$ for some $c_{\adl{ubge}}>0$. Take $V\subset S$ so that $P\cap V=\emptyset$ and $U\setminus P\subset V$. Then $\|\LR_{\e}[P,V,\lam_{\e}]-\LR[P,U\setminus P,\lam]\|_{\infty}\to 0$ as $\e\to 0$.
}
\pros
We write
\ali
{
\LR_{\e}[P,S\setminus P,\lam_{\e}]=&(\LR_{\e})_{PP}+\mathscr{M}_{\e}\\
\LR_{\e}[P,V,\lam_{\e}]=&(\LR_{\e})_{PP}+\mathscr{P}_{\e}\\
\LR_{\e}[P,U\setminus P,\lam_{\e}]=&(\LR_{\e})_{PP}+\mathscr{N}_{\e}.
}
Notice that $\mathscr{P}_{\e}-\mathscr{N}_{\e}$ and $\mathscr{M}_{\e}-\mathscr{N}_{\e}$ are positive operators by $U\setminus P\subset V\subset S\setminus P$. In order to show the assertion, we remark $\|\mathscr{P}_{\e}-\mathscr{N}_{\e}\|_{\infty}=\|(\mathscr{P}_{\e}-\mathscr{N}_{\e})1\|_{\infty}\leq \|(\mathscr{M}_{\e}-\mathscr{N}_{\e})1\|_{\infty}$. Therefore we may show $\|(\mathscr{M}_{\e}-\mathscr{N}_{\e})1\|_{\infty}\to 0$.
Put $f_{\e}^{P}=g_{\e}^{P}/\nu^{P}(U,g_{\e}^{P})$. By the assumption, we see $f_{\e}^{P}\geq c_{\adr{ubge}}\chi_{P}$. 
We have
\ali
{
\lam_{\e}f_{\e}^{P}=&\LR_{\e}[P,S\setminus P,\lam_{\e}]f_{\e}^{P}\\
=&\LR_{\e}[P,U\setminus P,\lam_{\e}]f_{\e}^{P}+(\mathscr{M}_{\e}-\mathscr{N}_{\e})f_{\e}^{P}\\
\geq& \LR_{\e}[P,U\setminus P,\lam_{\e}]f_{\e}^{P}+(\inf_{\Si_{P}}f_{\e}^{P})(\mathscr{M}_{\e}-\mathscr{N}_{\e})1.
}
Then for $\om\in \Si_{P}$
\ali
{
|(\mathscr{M}_{\e}-\mathscr{N}_{\e})1(\om)|\leq c_{\adr{ubge}}^{-1}|\lam_{\e}f_{\e}^{P}(\om)-\LR_{\e}[P,U\setminus P,\lam_{\e}]f_{\e}^{P}(\om)|.
}
For $a\in P$, fix $\om(a)\in [a]$. 
\ali
{
\|(\mathscr{M}_{\e}-\mathscr{N}_{\e})1\|_{\infty}=&\max_{a\in P}\sup_{\om\in [a]}\sum_{w\in W(P:S\setminus P)\setminus W(P:U\setminus P)\,:\,|w|\geq 2,\ \atop{w\cd a\in W(A)}}\lam_{\e}^{-|w|+1}e^{S_{|w|}\ph(\e,w\cd\om)}\\
\leq &\max_{a\in P}\sum_{w\in W(P:S\setminus P)\setminus W(P:U\setminus P)\,:\,|w|\geq 2,\ \atop{w\cd a\in W(A)}}\lam_{\e}^{-|w|+1}c_{\adr{cbd}}e^{S_{|w|}\ph(\e,w\cd\om(a))}\\
\leq &\max_{a\in P} (\mathscr{M}_{\e}-\mathscr{N}_{\e})1(\om(a))\\
\leq &c_{\adr{ubge}}^{-1}\max_{a\in P}|\lam_{\e}f_{\e}^{P}(\om(a))-\LR_{\e}[P,U\setminus P,\lam_{\e}]f_{\e}^{P}(\om(a))|\\
\to &0
}
as $\e\to 0$, where the last convergence uses a similar argument of (\ref{eq:Lge-Lg0}). Hence the assertion is valid.
\proe
\prop
{\label{prop:spec_PC_Rue}
Assume that conditions (A.1)-(A.3) are satisfied. Let $U\in \mathscr{S}_{0}$ and $P\subset U$ a nonempty subset. We also assume that $h^{P}(U,\cd)\geq c_{\adr{ubh}}\chi_{P}$ for some $c_{\adl{ubh}}>0$. Then $\LR[P,U\setminus P,\lam]\,:\,F_{b}(X)\to F_{b}(X)$ has the spectral decomposition
\alil
{
\LR[P,U\setminus P,\lam]=\lam\mathscr{P}+\mathscr{R}\label{eq:LPUPl=}
}
satisfying that $\lam$ is the spectral radius of the operator $\LR[P,U\setminus P,\lam]$ and becomes the simple eigenvalue, 
$\mathscr{P}$ is the projection onto the one-dimensional eigenspace of the $\lam$ which has the form $\mathscr{P}=h^{P}(U,\cd)\otimes \nu^{P}(U,\cd)$, and the spectrum of $\mathscr{R}$ is in the outside of a neighborhood of the $\lam$. 
}
\pros
Let $V=U\setminus P$.
By using Theorem \ref{th:spec_PC}, $\LR[P,V,\lam]$ has the spectral gap at the simple eigenvalue of $\lam$. It suffices to check that $\lam$ is the spectral radius of this operator. Let $f\in F_{b}(X)$ and $\om,\up\in X$ with $\om_{0}=\up_{0}$. By the basic inequality
\alil
{
&|\LR[P,V,\lam]^{n}f(\om)-\LR[P,V,\lam]^{n}f(\up)|\label{eq:spec_PC_Rue1}\\
\leq &\sum_{w\in W(U:V)^{n}\,:\,\atop{w\cd\om_{0}\in W(B[U])}}\{\lam^{-|w|+1}e^{S_{|w|}\ph(w\cd\om)}|f(w\cd\om)-f(w\cd\up)|+\nonumber\\
&\quad \lam^{-|w|+1}e^{S_{|w|}\ph(w\cd\up)}|e^{S_{|w|}\ph(w\cd\om)-S_{|w|}\ph(w\cd\up)}-1||f(w\cd\up)|\}\nonumber\\
\leq&\|\LR[P,V,\lam]^{n}1\|_{\infty}([f]_{1}\theta+c_{\adr{bi}}\|f\|_{\infty})d_{\theta}(\om,\up)\nonumber
}
with $c_{\adl{bi}}=e^{\theta^{2}[\ph]_{2}/(1-\theta)}\theta[\ph]_{2}/(1-\theta)$. Therefore $\|\LR[P,V,\lam]^{n}\|_{1}\leq \|\LR[P,V,\lam]^{n}1\|_{\infty}(\theta+c_{\adr{bi}})$. Thus the spectral radius of $\LR[P,V,\lam]$ is equal to $\lim_{n\to \infty}\|\LR[P,V,\lam]^{n}1\|_{\infty}^{1/n}$. Moreover, we see
\alil
{
\|\LR[P,V,\lam]^{n}1\|_{\infty}\leq &\frac{\|\LR[P,V,\lam]^{n}h^{P}(U,\cd)\|_{\infty}}{\inf_{\Si_{P}}h^{P}(U,\cd)}\leq c_{\adr{ubh}}^{-1}\lam^{n}\|h^{P}(U,\cd)\|_{\infty}.\label{eq:spec_PC_Rue2}
}
We get $\lim_{n\to \infty}\|\LR[P,V,\lam]^{n}1\|_{\infty}^{1/n}\leq \lam$. Hence we see the spectral radius of $\LR[P,V,\lam]$ equals $\lam$.
\proe
\prop
{\label{prop:conv_PC_nue}
Assume that conditions (A.1)-(A.3) are satisfied. Let $U\in \mathscr{S}_{0}$ and $P\subset S$ a nonempty subset. We also assume that $g_{\e}^{P}\geq c_{\adr{ubge}}\chi_{P}$ for some $c_{\adr{ubge}}>0$. Then $\nu_{\e}^{P}$ converges to $\nu^{P}(U,\cd)$ in the functional on $F_{b}(X)$ as $\e\to 0$.
}
\pros
Together with Proposition \ref{prop:gP_nonzero_inf}, the assumption $g_{\e}^{P}\geq c_{\adr{ubge}}\chi_{P}$ implies $h(U,\cd)\geq c_{\adr{ubge}}\chi_{P}$. Therefore conditions of Proposition \ref{prop:conv_PC_Rue2} is satisfied. Let $Q=U\setminus P$. We use the spectral decomposition (\ref{eq:LPUPl=}) in Proposition \ref{prop:spec_PC_Rue}. By using \cite[Theorem 2.2]{T2022pre} (or \cite[Apeendix B]{T2024_sum}), the measure $\kappa_{\e}:=\nu_{\e}^{P}(f)/\nu_{\e}^{P}(h^{P}(U,\cd))$ has the formula
\alil
{
\kappa_{\e}(f)-\nu^{P}(U,f)=\kappa_{\e}((\LR_{\e}[P,S\setminus P,\lam_{\e}]-\LR[P,U\setminus P,\lam])(h\otimes \kappa_{\e}-\IR)(\mathscr{R}-\lam\IR)^{-1}f)\label{eq:kapef-nUf=}
}
for each $f\in F_{b}(X)$. 
To show the assertion, we note that $\kappa_{\e}$ is a finite measure with the upper bound by $c_{\adl{ubk}}:=(\min_{a\in P}\inf_{[a]}h^{P}(U,\cd))^{-1}$. Then we see that for $f\in F_{b}(X)$ with $\|f\|_{1}\leq 1$,
\ali
{
&|\kappa_{\e}(f)-\nu^{P}(U,f)|\\
\leq &c_{\adr{ubk}}\|\LR_{\e}[P,S\setminus P,\lam_{\e}]-\LR[P,U\setminus P,\lam]\|_{\infty}\|h\otimes \kappa_{\e}-\IR\|_{\infty}\|(\mathscr{R}-\lam\IR)^{-1}f\|_{\infty}\\
\leq& c_{\adr{ubk}}\|\LR_{\e}[P,S\setminus P,\lam_{\e}]-\LR[P,U\setminus P,\lam]\|_{\infty}(c_{\adr{ubk}}\|h\|_{\infty}+1)\|(\mathscr{R}-\lam\IR)^{-1}\|_{1}\to 0
}
as $\e\to 0$ by Proposition \ref{prop:conv_PC_Rue2}. Thus we obtain $\kappa_{\e}(f)\to \nu^{P}(U,f)$. Finally letting $f:=1$, this convergence implies $1/\nu_{\e}^{P}(h^{P}(U,\cd))\to 1$. Hence we get convergence of the assertion.
\proe
\prop
{\label{prop:conv_PC_mue}
Assume that conditions (A.1)-(A.3) are satisfied. Let $U\in \mathscr{S}_{0}$ and $P\subset S$ a nonempty subset. We also assume that $g_{\e}^{P}\geq c_{\adr{ubge}}\chi_{P}$ for some $c_{\adr{ubge}}>0$. Then $\mu_{\e}^{P}$ converges weakly to $\mu^{P}(U,\cd)$ as $\e\to 0$.
}
\pros
Put $f_{\e}^{P}=g_{\e}^{P}/\nu^{P}(U,g_{\e}^{P})$. Note the equation
\ali
{
\mu_{\e}^{P}(f)=\frac{\nu_{\e}^{P}(g_{\e}^{P}f)}{\nu_{\e}^{P}(g_{\e}^{P})}=\frac{\nu_{\e}^{P}(f_{\e}^{P}f)}{\nu_{\e}^{P}(f_{\e}^{P})}\qquad \mu^{P}(U,f)=\nu^{P}(U,h^{P}(U,\cd)f).
}
By Proposition \ref{prop:gP_nonzero_inf} and Proposition \ref{prop:conv_PC_nue}, we have that for each $f\in F_{b}(X)$
\ali
{
&|\nu_{\e}^{P}(f_{\e}^{P}f)-\mu^{P}(U,f)|\\
\leq&\|\nu_{\e}^{P}-\nu^{P}(U,\cd)\|_{F_{b}(X)\to \R}\|f_{\e}^{P}f\|_{1}+|\nu^{P}(U,f_{\e}^{P} f-h^{P}(U,\cd)f)|\to 0
}
as $\e\to 0$ in addition to $\sup_{\e>0}\|f_{\e}^{P}\|_{1}<+\infty$. In particular, the case $f\equiv 1$ yields $\nu_{\e}^{P}(f_{\e}^{P})\to \mu^{P}(U,1)=1$. Hence the assertion is valid.
\proe
For an operator $\MR\in \LR(\XR)$, denoted by $r(\MR$ the spectral radius. Next we show the existence of Perron spectral triplet of operators $\LR_{\e}[U,V,\eta]$ which is an important role for proof our main result.
\prop
{\label{prop:ex_PSC_PC}
Assume that conditions (A.1)-(A.3) are satisfied. Let $U\in \mathscr{S}_{0}$ and $P\subset S$ a nonempty subset with $r(\LR_{PP})>0$. Take $V\subset S$ so that $U\setminus P\subset V$. Put $W=P\cup V$. Let $\lam$ be a number with $\lam\geq r(\LR_{WW})$ and $\lam>r(\LR_{VV})$. Then there exists a triplet $(\hat{\lam},\hat{h},\hat{\nu})$ such that
\ite
{
\item $\hat{\lam}$ is a positive eigenvalue of the operator $\LR[P,V,\lam]$. If $P$ is finite, then $\hat{\lam}$ becomes the spectral radius of this operator;
\item $\hat{\nu}$ is the Borel probability measure with $\hat{\nu}(\Si_{P})>0$ such that $\LR[P,V,\lam]^{*}\hat{\nu}=\hat{\lam}\hat{\nu}$;
\item $\hat{h}\in \Lambda_{c_{\adr{cL}}}$ satisfies $\hat{h}\neq 0$ on $\Si_{P}$ and $\LR[P,V,\lam]\hat{h}=\hat{\lam}\hat{h}$ with $\hat{\nu}(\hat{h})=1$.
}
}
\pros
For $t>0$, we define
\ali
{
\ph_{t}:=
\case
{
\ph, &\text{ on }\Si_{U}\\
\ph+\log t, &\text{ on }X\setminus \Si_{U}\\
}
}
Observe the form
\ali
{
(\LR_{\ph_{t}})_{WW}=\LR_{WP}+t\LR_{WV}=
\MatII
{
\LR_{PP}&t\LR_{PV}\\
\LR_{VP}&t\LR_{VV}\\
}
}
It is not hard to check that $[\ph_{t}]_{2}<+\infty$ and $\ph_{t}$ is summable for each $t>0$. Then Theorem \ref{th:RPF_transitive_summ_anyA2}(1) implies that there exists the Perron spectral triplet $(\lam_{t},g_{t},\nu_{t})$ for the operator $(\LR_{\ph_{t}})_{WW}$, where $g_{t}$ is the corresponding eigenfunction with $\|g_{t}\|_{\infty}=1$. In order to show the assertion, we remark the following inequality: for $0<t\leq 1$
\ali
{
0<r(\LR_{PP})\leq r((\LR_{\ph_{t}})_{WW})=\lam_{t}\leq r(\LR_{WW})\leq \lam.
}
By the continuity of $t\mapsto \lam_{t}/t$, there exists $0<s\leq 1$ such that $\lam_{s}/s=\lam$. In addition to the form $s=\lam_{s}/\lam$, Proposition \ref{prop:PC_evalprop}(1) implies $\LR[P,V,\lam](\chi_{P}g_{s})=\lam_{s}(\chi_{P}g_{s})$, $\LR[P,V,\lam]^{*}\nu_{s}(\chi_{P}f)=\lam_{s}\nu_{s}(\chi_{P}f)$ for $f\in F_{b}(X)$, $\chi_{P}g_{s}\neq 0$ and $\nu_{s}(\chi_{P})>0$. Now we show $g_{s}>0$ on $\Si_{P}$. Choose any $\om,\up\in \Si_{P}$ with $g_{s}(\up)>0$. Since $U\setminus P\subset V$ and $B[U]$ is irreducible, there exists $w:=w_{1}\cdots w_{m}\in W(P:V)^{m}$, $w\cd \om_{0}$ is $B[U]$-admissible and $w_{0}=\up_{0}$. We have
\ali
{
g_{s}(\om)=&\lam_{s}^{-m}\LR[P,V,\lam]^{m}g_{s}(\om)\\
\geq& \lam_{s}^{-m}\lam^{-|w|+1}e^{S_{|w|}\ph(w\cd\om)}g_{s}(w\cd\om)\\
\geq&  \lam_{s}^{-m}\lam^{-|w|+1}e^{S_{|w|}\ph(w\cd\om)}c_{\adr{cL}}^{-1}g_{s}(\up)>0.
}
Thus $g_{s}>0$ on $\Si_{P}$. 
Letting $\hat{\lam}:=\lam_{s}$, $\hat{\nu}:=\nu_{s}(\cd|\Si_{P})$ and $\hat{h}:=\chi_{P}f_{s}/\hat{\nu}(\chi_{P}f_{s})$, we obtain the assertions (1)-(3) excluding a part of (1). Finally we assume $P$ is finite and prove that $\hat{\lam}$ is the spectral radius of $\LR[P,V,\lam]$. By a similar argument in (\ref{eq:spec_PC_Rue1}), we see that $\hat{\lam}$ is equal to the spectral radius of $\LR[P,V,\lam]\,:\,C_{b}(X)\to C_{b}(X)$. Moreover, the facts $\hat{h}\in \Lambda_{c_{\adr{cL}}}$ and $\hat{h}>0$ on $\Si_{P}$ imply $\inf_{\Si_{P}}\hat{h}>0$. Thus a similar argument in (\ref{eq:spec_PC_Rue2}) says $\hat{\lam}=r(\LR[P,V,\lam])$. Hence the proof is complete.
\proe
\prop
{\label{prop:ex_PSC_PC_per}
Assume that conditions (A.1)-(A.3) are satisfied. Let $U\in \mathscr{S}_{0}$ and $P\subset S$ a nonempty subset with $r((\LR_{\e})_{PP})>0$ for some $\e>0$. Assume also that $g_{\e}^{P}\geq c_{\adr{ubge}}\chi_{P}$ for some $c_{\adr{ubge}}>0$. Take $V\subset S$ so that $U\setminus P\subset V\subset S\setminus P$. Put $W=P\cup V$. Take the Perron spectral triplet $(\hat{\lam}_{\e},\hat{h}_{\e},\hat{\nu}_{\e})$ from Proposition \ref{prop:ex_PSC_PC} for the operator $\LR_{\e}[P,V,\lam_{\e}]$. Then
\ite
{
\item $\hat{\lam}_{\e}\to \lam$ as $\e\to 0$;
\item $\hat{h}_{\e}(\om)\to h^{P}(U,\om)$ for each $\om\in X$;
\item $\hat{\nu}_{\e}\to \nu^{P}(U,\cd)$ in $(F_{b}(X))^{*}$.
}
}
\pros
(1) Note that since the symbolic structure in $(\LR_{\e})_{PP}$ does not chenge for any positive $\e>0$, the assumption $r((\LR_{\e})_{PP})>0$ for some $\e>0$ yields $r((\LR_{\e})_{PP})>0$ for any $\e>0$. 
In addition to $r(\LR_{\e}[P,U\setminus P,\lam_{\e}])\to \lam$ and $\lam_{\e}\to \lam$, the inequality
\ali
{
r(\LR_{\e}[P,U\setminus P,\lam_{\e}])\leq r(\LR_{\e}[P,V,\lam_{\e}])=\hat{\lam}_{\e}\leq r(\LR_{\e}[P,S\setminus P,\lam_{\e}])=\lam_{\e}
}
implies $\hat{\lam}_{\e}\to \lam$.
\smallskip
\\
(2) By a similar argument in the proof of Proposition \ref{prop:gP_nonzero_inf}, we see the assertion.
\smallskip
\\
(3) By using the formula (\ref{eq:kapef-nUf=}) and using the fact $\|\LR_{\e}[P,V,\lam_{\e}]-\LR[P,U\setminus P,\lam]\|_{\infty}\to 0$, we get $\hat{\nu}_{\e}\to \nu^{P}(U,\cd)$. 
\proe
\subsection{Decomposition of the measure $\mu_{\e}$}\label{sec:coupling}
In this section, we show decomposition of the measure $\mu_{\e}$ and related results via Perron complements of Ruelle operators. We will develop techniques in \cite{T2020} from finite state case to infinite state case.

Assume that conditions (A.1)-(A.3) are satisfied. 
Choose any $W\subset S$. We decompose $W$ into a finite number of non-empty subsets:
\ali
{
W=W(1)\cup W(2)\cup\cdots \cup W(m_{0}).
}
We will use mainly techniques in \cite[Section 3]{T2020}.
Denote by $(\lam_{\e},h_{\e}^{W(i)},\nu_{\e}^{W(i)})$ the Perron spectral triplet of $\LR_{\e}[W(i),S\setminus W(i),\lam_{\e}]$, namely $h_{\e}^{W(i)}$ is the corresponding eigenfunction of the eigenvalue $\lam_{\e}$ of this operator, and $\nu_{\e}^{W(i)}$ is the corresponding eigenvector of $\lam_{\e}$ of the dual of this operator with $\nu_{\e}^{W(i)}(1)=\nu_{\e}^{W(i)}(h_{\e}^{W(i)})=1$. Recall the equations
\ali
{
h_{\e}^{W(i)}=\frac{h_{\e}\chi_{W(i)}}{\nu_{\e}^{W(i)}(h)}=\frac{h_{\e}^{W}\chi_{W(i)}}{\nu_{\e}^{W(i)}(h_{\e}^{W})},\qquad \nu_{\e}^{W(i)}(f)=\nu_{\e}(f | \Si_{W(i)})=\frac{\nu_{\e}(f\chi_{W(i)})}{\nu_{\e}(\Si_{W(i)})}=\frac{\nu_{\e}^{W}(f\chi_{W(i)})}{\nu_{\e}^{W}(\Si_{W(i)})}.
}
Recall the conditional probability measure $\mu_{\e}^{W}$ of $\mu_{\e}$ which is defined by $\mu_{\e}^{W}=h_{\e}^{W}\nu_{\e}^{W}$. We shall give a necessary and sufficient condition for convergence of $\mu_{\e}^{W}$.
We define matrices $D_{\e}^{W}=(D_{\e}^{W}(ij))$ and $E_{\e}^{W}=(E_{\e}^{W}(ij))$ indexed as $T_{0}:=\{1,2,\dots, m_{0}\}$ by
\alil
{
D_{\e}^{W}(ij):=&\nu_{\e}^{W(i)}(\LR_{\e}[W,S\setminus W,\lam_{\e}]h_{\e}^{W(j)})\nonumber\\
E_{\e}^{W}(ij):=&\nu_{\e}^{W(i)}(\LR_{\e}[W,S\setminus W,\lam_{\e}]\chi_{W(j)}).\label{eq:D(ij)=}
}
Note $\supp\,h_{\e}^{W(i)}=\Si_{W(i)}$ and $\supp\,\nu_{\e}^{W(i)}=X_{W(i)}$.
We set the column vector $\beta_{\e}^{W}=(\beta_{\e}^{W}(i))_{i\in T_{0}}$ and the low vector $\gamma_{\e}^{W}=(\gamma_{\e}^{W}(i))_{i\in T_{0}}$ by
\ali
{
\beta_{\e}^{W}(i)=\nu_{\e}^{W(i)}(h_{\e}^{W}),\quad \gamma_{\e}^{W}(i)=\nu_{\e}^{W}(\chi_{W(i)}).
}
\prop
{\label{prop:DeEe_irre2_inf}
Assume that conditions (A.1)-(A.3) are satisfied. Choose any $W\subset S$ and any decomposition $W(1), W(2),\dots, W(m_{0})$ with $W(k)\neq \emptyset$ $(k\in T_{0})$. Then $D_{\e}^{W}$ and $E_{\e}^{W}$ defined the above notation are irreducible. Moreover, the right positive eigenvector of $D_{\e}^{W}$ is equal to $\beta_{\e}^{W}$, and the left positive eigenvector of $D_{\e}^{W}$ equals $\gamma_{\e}^{W}$ with $\sum_{k\in T_{0}}\beta_{\e}^{W}(k)\gamma_{\e}^{W}(k)=\sum_{k\in T_{0}}\gamma_{\e}^{W}(k)=1$. Moreover, let $\beta_{\e}=(\beta_{\e}(k))$ be the Perron right eigenvector of $D_{\e}^{W}$ and $\gamma_{\e}={}^{t}(\gamma_{\e}(k))$ the Perron left eigenvector of $E_{\e}^{W}$ with $\sum_{k\in T_{0}}\gamma_{\e}^{W}(k)=\sum_{k\in T_{0}}\beta_{\e}^{W}(k)\gamma_{\e}^{W}(k)=1$. Then the function $\sum_{k\in T_{0}}\beta_{\e}(k)h_{\e}^{W(k)}$ coincides with the Perron eigenfunction $h_{\e}^{W}$ of the operator $\LR_{\e}[W,S\setminus W,\lam_{\e}]$ and the measure $\sum_{k\in T_{0}}\gamma_{\e}^{W}(k)\nu_{\e}^{W}$ is equal to the Perron eigenvector $\nu_{\e}^{W}$ of this dual operator.
}
\pros
Throughout this proof, we fix the parameter $\e>0$. Then for simplicity, we remote the symbol `$\e$' from the notation (e.g. $\LR=\LR_{\e}$, $\lam=\lam_{\e}$, $h=h_{\e}$, $\nu=\nu_{\e}$ etc.). 
By the definitions $D^{W}$ and $E^{W}$, $D^{W}(ij)>0$ iff $E^{W}(ij)>0$ iff $\LR[W,S\setminus W,\lam]_{W(i)W(j)}\neq O$. Since $A$ is irreducible, so are for $D^{W}$ and $E^{W}$.
We have that
\ali
{
(D^{W} \beta^{W})(i)=&\sum_{j=1}^{m_{0}}D^{W}(ij)\beta^{W}(j)=\sum_{j=1}^{m_{0}}\nu^{W(i)}(\LR[W,S\setminus W,\lam]h^{W(j)})\nu^{W(j)}(h^{W})\\
=&\sum_{j=1}^{m_{0}}\nu^{W(i)}(\LR[W,S\setminus W,\lam]\frac{h^{W}\chi_{W(j)}}{\nu^{W(j)}(h^{W})})\nu^{W(j)}(h^{W})\\
=&\nu^{W(i)}(\LR[W,S\setminus W,\lam]h^{W})
=\nu^{W(i)}(\lam h^{W})
=\lam\beta^{W}
}
\ali
{
(\gamma^{W}E^{W})(i)=&\sum_{j=1}^{m_{0}}\gamma^{W}(j)E^{W}(ji)
=\sum_{j=1}^{m_{0}}\nu^{W}(\chi_{W(j)})\frac{\nu^{W}(\chi_{W(j)}\LR[W,S\setminus W,\lam]\chi_{W(i)})}{\nu^{W}(\Si_{W(j)})}\\
=&\nu^{W}(\LR[W,S\setminus W,\lam]\chi_{W(i)})
=\lam\nu^{W}(\chi_{W(i)})
=\lam\gamma^{W}(i)
}
for each $i$. Moreover 
$\sum_{i=1}^{m_{0}}\gamma^{W}(i)=\nu^{W}(\Si_{W})=1$
and
\ali
{
\sum_{i=1}^{m_{0}}\beta^{W}(i)\gamma^{W}(i)=\sum_{i=1}^{m_{0}}\frac{\nu^{W}(\chi_{W(i)}h^{W})}{\nu^{W}(\chi_{W(i)})}\nu^{W}(\chi_{W(i)})=\nu^{W}(\chi_{W}h^{W})=1
}
are satisfied. Thus $\beta^{W}$ and $\gamma^{W}$ are the Perron right eigenvalue and the Perron left eigenvalue of $D^{W}$ and $E^{W}$, respectively. We also obtain $\beta=\beta^{W}$ and $\gamma=\gamma^{W}$ by the simplicity of the eigenvalue $\lam$. Finally
\ali
{
\sum_{i=1}^{m_{0}}\beta(i)h^{W(i)}=\sum_{i=1}^{m_{0}}\nu^{W(i)}(h^{W(i)})h^{W(i)}=\sum_{i=1}^{m_{0}}\nu^{W(i)}(h^{W})\frac{h^{W}\chi_{W(i)}}{\nu^{W(i)}(h^{W})}=h^{W}
}
and
\ali
{
\sum_{i=1}^{m_{0}}\gamma(i)\nu^{W(i)}(f)=\sum_{i=1}^{m_{0}}\nu^{W}(\chi_{W(i)})\frac{\nu^{W}(f\chi_{W(i)})}{\nu^{W}(\Si_{W(i)})}=\nu^{W}(f).
}
Hence the proof is complete.
\proe
\lem
{\label{lem:mueQ=v1_2_inf}
Assume that conditions (A.1)-(A.3) are satisfied. Choose any $W\subset S$ and any decomposition $W(1), W(2),\dots, W(m_{0})$ with $W(k)\neq \emptyset$ $(k\in T_{0})$. Then $\mu_{\e}^{W}=\sum_{i\in T_{0}}\beta_{\e}^{W}(i)\gamma_{\e}^{W}(i)\mu_{\e}^{W(i)}$.
}
\pros
We omit the symbol `$\e$'. By Proposition \ref{prop:DeEe_irre2_inf}, we have
\ali
{
\mu^{W}(f)=\nu^{W}(h^{W} f)=&\sum_{i=1}^{m_{0}}\sum_{j=1}^{m_{0}}\gamma^{W}(i)\beta^{W}(j)\nu^{W(i)}(h^{W(j)}f)\\
=&\sum_{i=1}^{m_{0}}\sum_{j=1}^{m_{0}}\gamma^{W}(i)\beta^{W}(j)\nu^{W(i)}(\chi_{W(i)}\chi_{W(j)}h^{W(j)}f)\\
=&\sum_{i=1}^{m_{0}}\gamma^{W}(i)\beta^{W}(i)\nu^{W(i)}(h^{W(i)}f)\\
=&\sum_{i=1}^{m_{0}}\gamma^{W}(i)\beta^{W}(i)\mu^{W(i)}(f).
}
\proe
Finally we mention the case where $\LR_{\e}$ is a normalized Ruelle operator.
\cor{\label{cor:coup2_inf}
Assume that conditions (A.1)-(A.3) are satisfied. Choose any $W\subset S$ and any decomposition $W(1), W(2),\dots, W(m_{0})$ with $W(k)\neq \emptyset$ $(k\in T_{0})$. Assume also the equation $(\LR_{\e})_{WW} \chi_{W}=\lam_{\e} \chi_{W}$. Then the two matrices $D_{\e}^{W}$ and $E_{\e}^{W}$ defined by (\ref{eq:D(ij)=}) are identical, and these right eigenvector are a vector whose entries are all $1$. In particular, the conditional probability measure $\mu_{\e}^{W}$ of $\mu_{\e}$ has the form $\mu_{\e}^{W}=\sum_{i=1}^{m_{0}}\gamma_{\e}(i)\mu_{\e}^{W(i)}$ with the Perron left eigenvector $\gamma_{\e}(i)$ of $D_{\e}^{W}=E_{\e}^{W}$ with $\sum_{i\in T_{0}}\gamma_{\e}(i)=1$.
}
\pros
We omit the symbol `$\e$'. By Proposition \ref{prop:PC_evalprop}(1), the equation $\LR_{WW} \chi_{W}=\lam \chi_{W}$ gives $\LR[W(i),S\setminus W(i),\lam]\chi_{W(i)}=\lam \chi_{W(i)}$. By Proposition \ref{prop:DeEe_irre2_inf}, we obtain the form $\chi_{W}=\sum_{i=1}^{m_{0}}\beta(i)h^{W(i)}$. This implies that $h^{W(i)}$ is a constant function on $\Si_{W(i)}$. By $\nu^{W(i)}(h^{W(i)})=\nu^{W(i)}(\chi_{W(i)})=1$, we have $h^{W(i)}=\chi_{W(i)}$ and therefore $\beta(i)=1$ for all $i$. Thus $D^{W}(ij)=\nu^{W(i)}(\LR[W,S\setminus W,\lam]\chi_{W(j)})=E^{W}(ij)$. Finally, Proposition \ref{prop:DeEe_irre2_inf} yields the last assertion $\sum_{i\in T_{0}}\gamma_{\e}(i)=1$. Hence the assertion is fulfilled.
\proe
We decompose $Q$ into a finite number of non-empty subsets $Q=Q(1)\cup Q(2)\cup\cdots \cup Q(m_{0})$. For $i_{1},\dots, i_{n}\in T_{0}$ mutually disjoint and a subset $P_{0}\subset T_{0}\setminus\{i_{1},\dots, i_{n}\}$, put
\ali
{
b_{\e}(i_{1}\cdots i_{n}:P_{0}):=\lam_{\e}-\frac{\nu_{\e}^{Q}(\LR[\bigcup_{l=1}^{n}Q(i_{l}),(S\setminus Q)\cup \bigcup_{k\in P_{0}}Q(k),\lam_{\e}]g_{\e}^{Q})}{\nu_{\e}^{Q}(g_{\e}^{Q}\sum_{l=1}^{n}\chi_{Q(i_{l})})}.
}
Observe that this value does not change even if $i_{1},\dots, i_{n}$ are replaced.
For $i,j\in T_{0}$ with $i\neq j$
\ali
{
B_{\e}(i,j):=\frac{b_{\e}(i : T_{0}\setminus \{i,j\})}{b_{\e}(j : T_{0}\setminus \{i,j\})}.
}
For $k\in T_{0}$, let
\ali
{
\ti{\delta}_{\e}(k):=\frac{1}{1+\sum_{j\in T_{0}\,:\,j\neq k}B_{\e}(k,j)}.
}
\prop
{\label{prop:mueQ=sumk_inf}
Assume that conditions (A.1)-(A.3) are satisfied. Choose any $Q\subset S$ and any decomposition $Q(1), Q(2),\dots, Q(m_{0})$ with $Q(k)\neq \emptyset$ $(k\in T_{0})$. Then we have the form $\mu_{\e}^{Q}=\sum_{k\in T_{0}}\ti{\delta}_{\e}(k)\mu_{\e}^{Q(k)}$. 
}
\pros
We omit the symbol `$\e$'. 
For pairwise disjoint elements $i_{1},\dots, i_{n}, j_{1}\dots, j_{u}\in T_{0}$, we let
\ali{
f(i_{1}\cdots i_{n},j_{1}\cdots j_{u})&=
\frac{b(i_{1}\cdots i_{n}:P_{0})}{b(i_{1}\cdots i_{n} : P_{0})+b(j_{1}\cdots j_{u} : P_{0})},
}
where $P_{0}$ denotes $T_{0}\setminus \{i_{1},\dots, i_{n}, j_{1}\dots, j_{u}\}$. Observe that the following is satisfied:
\alil
{
f(j_{1}\cdots j_{u},i_{1}\cdots i_{n})&=1-f(i_{1}\cdots i_{n},j_{1}\cdots j_{u}).\label{eq:f(j1...)=...}
}
We first show that the equation
\alil
{
f(j_{1}\cdots j_{n},j_{0})=\frac{f(j_{1},j_{2}\cdots j_{n})f(j_{2}\cdots j_{n},j_{0})}{f(j_{0},j_{2}\cdots j_{n})+f(j_{1},j_{2}\cdots j_{n})f(j_{2}\cdots j_{n},j_{0})}\label{eq:fe=0_2}
}
is satisfied for any disjoint elements $j_{0},\cdots, j_{n}\in T_{0}$.
We write 
\ali
{
T_{0}=\{k_{1},k_{2},\dots, k_{m_{0}}\}.
}
Let $\ti{Q}(k_{j})=\bigcup_{l=j+1}^{m_{0}}Q(k_{l})$ for $j=1,2,\dots, m_{0}-1$. Denote $\ti{\LR}$ by the normalized Ruelle operator of $\LR$, namely, $\ti{\LR}=\LR_{\ph-\log h\circ \si+\log h}$. By virtue of Proposition \ref{prop:DeEe_irre2_inf} with $W(1)=Q(k_{1})$ and $W(2)=\ti{Q}(k_{1})$, the Perron eigenvector $\mu^{Q}$ of $\tLR[Q,S\setminus Q,\lam]^{*}$ has the decomposition $\mu^{Q}=\delta^{1}(1)\mu^{Q(k_{1})}+\delta^{1}(2)\mu^{\ti{Q}(k_{1})}$,
where $\delta^{1}=(\delta^{1}(i))_{i=1}^{2}$ is the left Perron eigenvector with $\delta^{1}(1)+\delta^{1}(2)=1$ of the $2\times 2$ matrix $E^{1}=(E^{1}(ij))$ whose entries are 
\ali
{
E^{1}=\matII{\mu^{Q(k_{1})}(\tLR[Q(k_{1}),S\setminus Q,\lam]1)&\mu^{Q(k_{1})}(\tLR[Q,S\setminus Q,\lam]_{Q(k_{1})\ti{Q}(k_{1})}1)\\
\mu^{\ti{Q}(k_{1})}(\tLR[Q,S\setminus Q,\lam]_{\ti{Q}(k_{1})Q(k_{1})}1)&\mu^{\ti{Q}(k_{1})}(\tLR[\ti{Q}(k_{1}),S\setminus Q,\lam]1)
}.
}
Note that its Perron eigenvalue is equal to $\lam$ and the equations
\ali
{
E^{1}(11)=&\frac{\nu(\LR[Q(k_{1}),S\setminus Q,\lam]g)}{\nu(g\chi_{Q(k_{1})})}=\frac{\nu^{Q}(\LR[Q(k_{1}),S\setminus Q,\lam]g^{Q})}{\nu^{Q}(g^{Q}\chi_{Q(k_{1})})}\\
E^{1}(22)=&\frac{\nu(\LR[\ti{Q}(k_{1}),S\setminus Q,\lam]g)}{\nu(g\chi_{\ti{Q}(k_{1})})}=\frac{\nu^{Q}(\LR[\bigcup_{l=2}^{m_{0}}Q(k_{l}),S\setminus Q,\lam]g^{Q})}{\nu^{Q}(g^{Q}\sum_{l=2}^{m_{0}}\chi_{Q(k_{l})})}
}
are satisfied by using $\mu=g\nu/\nu(g)=g^{Q}\nu^{Q}/\nu^{Q}(g^{Q})$ with $g^{Q}=h^{Q}/\|h^{Q}\|_{\infty}$. Since 
In addition to the equation $\tLR 1=\lam 1$, it follows from Corollary \ref{cor:coup2_inf} that the right Perron eigenvector $\beta^{1}=(\beta^{1}(i))_{i=1}^{2}$ of $E^{1}$ becomes a vector whose entries are $1$.
The basic result of matrix theory  (e.g. Corollary 2, p 9 in \cite{Seneta}) yields
\ali
{
\delta^{1}(1)=&\gamma^{1}(1)\beta^{1}(1)=\frac{\lam-E^{1}(22)}{\lam-E^{1}(11)+\lam-E^{1}(22)}=f(k_{2}\cdots k_{m_{0}}, k_{1})\\
\delta^{1}(2)=&\gamma^{1}(2)\beta^{1}(2)=\frac{\lam-E^{1}(11)}{\lam-E^{1}(11)+\lam-E^{1}(22)}=f(k_{1}, k_{2}\cdots k_{m_{0}}).
}
Thus
\alil
{
\mu^{Q}=&f(k_{2}\cdots k_{m_{0}}, k_{1})\mu^{Q(k_{1})}+f(k_{1}, k_{2}\cdots k_{m_{0}})\mu^{\ti{Q}(k_{1})}.\label{eq:mQ=...}
}
Observe that the measure $\mu^{\ti{Q}(k_{1})}$ is the Perron eigenvector of $\tLR[\ti{Q}(k_{1}),S\setminus \ti{Q}(k_{1}),\lam]^{*}$.

Next we give a decomposition of $\mu^{\ti{Q}(k_{1})}$. 
By Proposition \ref{prop:DeEe_irre2_inf} with $W(1)=Q(k_{2})$ and $W(2)=\ti{Q}(k_{2})$, the Perron eigenvector $\mu^{\ti{Q}(k_{1})}$ of $\tLR[\mu^{\ti{Q}(k_{1})},S\setminus \mu^{\ti{Q}(k_{1})},\lam]^{*}$ has the decomposition $\mu^{\ti{Q}(k_{1})}=\delta^{2}(1)\mu^{Q(k_{2})}+\delta^{2}(2)\mu^{\ti{Q}(k_{2})}$,
where $\delta^{2}=(\delta^{2}(i))_{i=1}^{2}$ is the left Perron eigenvector with $\delta^{2}(1)+\delta^{2}(2)=1$ of the $2\times 2$ matrix $E^{2}=(E^{2}(ij))$ whose entries are 
\ali
{
E^{2}=\matII{\mu^{Q(k_{2})}(\tLR[Q(k_{2}),V,\lam]1)&\mu^{Q(k_{2})}(\tLR[\ti{Q}(k_{1}),V,\lam]_{Q(k_{2})\ti{Q}(k_{2})}1)\\
\mu^{\ti{Q}(k_{2})}(\tLR[\ti{Q}(k_{1}),V,\lam]_{\ti{Q}(k_{2})Q(k_{2})}1)&\mu^{\ti{Q}(k_{2})}(\tLR[\ti{Q}(k_{2}),V,\lam]1)
}.
}
with $V=S\setminus \ti{Q}(k_{1})$.
We also have
\ali
{
E^{2}(11)=&\frac{\nu(\LR[Q(k_{2}),S\setminus \ti{Q}(k_{1}),\lam]g)}{\nu(g\chi_{Q(k_{2})})}=\frac{\nu^{Q}(\LR[Q(k_{2}),S\setminus \ti{Q}(k_{1}),\lam]g^{Q})}{\nu^{Q}(g^{Q}\chi_{Q(k_{2})})}\\
E^{2}(22)=&\frac{\nu(\LR[\ti{Q}(k_{2}),S\setminus \ti{Q}(k_{1}),\lam]g)}{\nu(g\chi_{\ti{Q}(k_{2})})}=\frac{\nu^{Q}(\LR[\bigcup_{l=3}^{m_{0}}Q(k_{l}),S\setminus \ti{Q}(k_{1}),\lam]g^{Q})}{\nu^{Q}(g^{Q}\sum_{l=3}^{m_{0}}\chi_{Q(k_{l})})}.
}
The right Perron eigenvector $\beta^{2}=(\beta^{2}(i))_{i=1}^{2}$ of $E^{2}$ becomes a vector whose entries are $1$. We obtain
\ali
{
\delta^{2}(1)=&\gamma^{2}(1)\beta^{2}(1)=\frac{\lam-E^{2}(22)}{\lam-E^{2}(11)+\lam-E^{2}(22)}=f(k_{3}\cdots k_{m_{0}}, k_{2})\\
\delta^{2}(2)=&\gamma^{2}(2)\beta^{2}(2)=\frac{\lam-E^{2}(11)}{\lam-E^{2}(11)+\lam-E^{2}(22)}=f(k_{2}, k_{3}\cdots k_{m_{0}}).
}
Therefore
\alil
{
\mu^{\ti{Q}(k_{1})}=&f(k_{3}\cdots k_{m_{0}}, k_{2})\mu^{Q(k_{2})}+f(k_{2}, k_{3}\cdots k_{m_{0}})\mu^{\ti{Q}(k_{2})}.\label{eq:mutQk1=}
}
Thus (\ref{eq:mQ=...}) and (\ref{eq:mutQk1=}) imply
\ali
{
\mu^{Q}=&f(k_{2}\cdots k_{m_{0}}, k_{1})\mu^{Q(k_{1})}+f(k_{3}\cdots k_{m_{0}}, k_{2})f(k_{1}, k_{2}\cdots k_{m_{0}})\mu^{Q(k_{2})}\\
&+f(k_{1}, k_{2}\cdots k_{m_{0}})f(k_{2}, k_{3}\cdots k_{m_{0}})\mu^{\ti{Q}(k_{2})}.
}
By using Proposition \ref{prop:DeEe_irre2_inf} repeatedly, we have the decomposition
\ali
{
\mu^{\ti{S}(k_{j-1})}=\delta^{j}(1)\mu^{S(k_{j})}+\delta^{j}(2)\mu^{\ti{S}(k_{j})}
}
for $j=2,3,\dots, m_{0}-1$, and $\delta^{j}=(\delta^{j}(i))_{i=1}^{2}$ is the left Perron eigenvalue with $\sum_{i}\delta^{j}(i)=1$ of the matrix
\alil
{
E^{j}=
\matII
{
\mu^{Q(k_{j})}(\tLR[Q(k_{j}),V,\lam]1)&\mu^{Q(k_{j})}(\tLR[\ti{Q}(k_{j-1}),V ,\lam]_{Q(k_{j})\ti{Q}(k_{j})}1)\\
\mu^{\ti{Q}(k_{j})}(\tLR[\ti{Q}(k_{j-1}),V,\lam]_{\ti{Q}(k_{j})Q(k_{j})}1)&\mu^{\ti{Q}(k_{j})}(\tLR[\ti{Q}(k_{j}),V,\lam]1)
}\label{eq:Ej=...}
}
with $V=S\setminus \ti{Q}(k_{j-1})$. Note the decomposition $\ti{Q}(k_{j-1})=Q(k_{j})\cup \ti{Q}(k_{j})$. 
Note also that the elements of the right eigenvector of this matrix are all $1$.
By a similar argument above, we have
\ali
{
\delta^{j}(1)=&f(k_{j+1}\cdots k_{m_{0}},k_{j}),\quad \delta^{j}(2)=f(k_{j},k_{j+1}\cdots k_{m_{0}})
}
and then
\alil
{
\mu^{\ti{S}(k_{j-1})}=&f(k_{j+1}\cdots k_{m_{0}},k_{j})\mu^{S(k_{j})}+f(k_{j},k_{j+1}\cdots k_{m_{0}})\mu^{\ti{S}(k_{j})}.\label{eq:mutSkj-1=...}
}
Consequently, $\mu^{Q}$ satisfies the form
\alil
{
\mu^{Q}=&f(k_{2}\cdots k_{m_{0}},k_{1})\mu^{S(k_{1})}+
\sum_{j=2}^{m_{0}-1}\big(f(k_{j+1}\cdots k_{m_{0}},k_{j})\prod_{l=1}^{j-1}f(k_{l},k_{l+1}\cdots k_{m_{0}})\big)\mu^{S(k_{j})}\label{eq:mue=expa}\\
&+\prod_{l=1}^{m_{0}-1}f(k_{l},k_{l+1}\cdots k_{m_{0}})\mu^{S(k_{m_{0}})}\nonumber\\
=&f(k_{2}\cdots k_{m_{0}},k_{1})\mu^{S(k_{1})}\nonumber\\
&+f(k_{3}\cdots k_{m_{0}},k_{2})f(k_{1},k_{2}\cdots k_{m_{0}})\mu^{S(k_{2})}\nonumber\\
&\vdots\nonumber\\
&+f(k_{l+1}\cdots k_{m_{0}},k_{l})f(k_{1},k_{2}\cdots k_{m_{0}})\cdots f(k_{l-1},k_{l}\cdots k_{m_{0}})\mu^{S(k_{l})}\nonumber\\
&+f(k_{l+2}\cdots k_{m_{0}},k_{l+1})f(k_{1},k_{2}\cdots k_{m_{0}})\cdots f(k_{l},k_{l+1}\cdots k_{m_{0}})\mu^{S(k_{l+1})}\nonumber\\
&\vdots\nonumber\\
&+f(k_{m_{0}},k_{m_{0}-1})f(k_{1},k_{2}\cdots k_{m_{0}})\cdots f(k_{m_{0}-2},k_{m_{0}-1}k_{m_{0}})\mu^{S(k_{m_{0}-1})}\nonumber\\
&+f(k_{1},k_{2}\cdots k_{m_{0}})\cdots f(k_{m_{0}-1},k_{m_{0}})\mu^{S(k_{m_{0}})}.\nonumber
}
Now for each $l=1,2,\dots, m_{0}-1$, we exchange $k_{l}$ and $k_{l+1}$ in the above expansion.
Namely
\ali
{
&\mu^{Q}\\
=&f(k_{2}\cdots k_{m_{0}},k_{1})\mu^{S(k_{1})}\\
&+f(k_{3}\cdots k_{m_{0}},k_{2})f(k_{1},k_{2}\cdots k_{m_{0}})\mu^{S(k_{2})}\\
&\vdots\\
&+f(k_{l}k_{l+2}\cdots k_{m_{0}},k_{l+1})f(k_{1},k_{2}\cdots k_{m_{0}})\cdots f(k_{l-1},k_{l}\cdots k_{m_{0}})\mu^{S(k_{l+1})}\\
&+f(k_{l+2}\cdots k_{m_{0}},k_{l})f(k_{1},k_{2}\cdots k_{m_{0}})\cdots f(k_{l-1},k_{l+1}k_{l}k_{l+2}\cdots k_{m_{0}})f(k_{l+1},k_{l}k_{l+2}\cdots k_{m_{0}})\mu^{S(k_{l})}\\
&\vdots\\
&+f(k_{m_{0}},k_{m_{0}-1})f(k_{1},k_{2}\cdots k_{m_{0}})\cdots f(k_{m_{0}-2},k_{m_{0}-1}k_{m_{0}})\mu^{S(k_{m_{0}-1})}\\
&+f(k_{1},k_{2}\cdots k_{m_{0}})\cdots f(k_{m_{0}-1},k_{m_{0}})\mu^{S(k_{m_{0}})}.
}
Compare the two coefficients of $\mu^{S(k_{l})}$ and the following equation is obtained:
\alil
{
f(k_{l+1}\cdots k_{m_{0}},k_{l})=&f(k_{l+2}\cdots k_{m_{0}},k_{l})f(k_{l+1},k_{l}k_{l+2}\cdots k_{m_{0}})\label{eq:fe=1_2}.
}
Similarity, 
by comparing the two coefficients of $\mu^{S(k_{l+1})}$
\ali
{
f(k_{l+2}\cdots k_{m_{0}},k_{l+1})f(k_{l},k_{l+1}\cdots k_{m_{0}})=f(k_{l}k_{l+2}\cdots k_{m_{0}},k_{l+1}).
}
It follows from the fact (\ref{eq:f(j1...)=...}) that the above equation is equivalent to
\alil
{
(1-f(k_{l+1},k_{l+2}\cdots k_{m_{0}}))(1-f(k_{l+1}\cdots k_{m_{0}},k_{l}))=1-f(k_{l+1},k_{l}k_{l+2}\cdots k_{m_{0}}).\label{eq:fe=2_2}
}
By (\ref{eq:fe=1_2}) and (\ref{eq:fe=2_2}), we have
\ali
{
\case
{
f(k_{l+1}\cdots k_{m_{0}},k_{l})-f(k_{l+2}\cdots k_{m_{0}},k_{l})f(k_{l+1},k_{l}k_{l+2}\cdots k_{m_{0}})=0\\
(f(k_{l+1},k_{l+2}\cdots k_{m_{0}})-1)f(k_{l+1}\cdots k_{m_{0}},k_{l})+f(k_{l+1},k_{l}k_{l+2}\cdots k_{m_{0}})=f(k_{l+1},k_{l+2}\cdots k_{m_{0}}).\\
}
}
By observing the equation (\ref{eq:mutSkj-1=...}), we have $0<f(k_{l+2}\cdots k_{m_{0}},k_{l}),\ f(k_{l+1},k_{l+2}\cdots k_{m_{0}})<1$. Thus, we see $1+f(k_{l+2}\cdots k_{m_{0}},k_{l})(f(k_{l+1},k_{l+2}\cdots k_{m_{0}})-1)\neq 0$. The Kramer formula yields the solution
\ali
{
f(k_{l+1}\cdots k_{m_{0}},k_{l})=&\frac{f(k_{l+2}\cdots k_{m_{0}},k_{l})f(k_{l+1},k_{l+2}\cdots k_{m_{0}})}{1+f(k_{l+2}\cdots k_{m_{0}},k_{l})(f(k_{l+1},k_{l+2}\cdots k_{m_{0}})-1)}\\
=&\frac{f(k_{l+2}\cdots k_{m_{0}},k_{l})f(k_{l+1},k_{l+2}\cdots k_{m_{0}})}{f(k_{l+2}\cdots k_{m_{0}},k_{l})f(k_{l+1},k_{l+2}\cdots k_{m_{0}})+f(k_{l},k_{l+2}\cdots k_{m_{0}})}.
}
By arbitrary choosing $k_{1},k_{2},\dots, k_{m_{0}}$, we get the equation
\ali
{
f(j_{1}\cdots j_{n},j_{0})=\frac{f(j_{1},j_{2}\cdots j_{n})f(j_{2}\cdots j_{n},j_{0})}{f(j_{0},j_{2}\cdots j_{n})+f(j_{1},j_{2}\cdots j_{n})f(j_{2}\cdots j_{n},j_{0})}.
}
Thus (\ref{eq:fe=0_2}) is fulfilled.
\smallskip
\par
Next we prove the equation
$B(i,j)B(j,k)=B(i,k)$ 
for disjoint elements $i,j,k\in T_{0}$. Since $f(ij,k)=f(ji,k)$ is satisfied by the definition of $f$, the equation (\ref{eq:fe=0_2}) implies
\ali
{
f(ij,k)=\frac{f(i,j)f(j,k)}{f(k,j)+f(i,j)f(j,k)}&=\frac{f(j,i)f(i,k)}{f(k,i)+f(j,i)f(i,k)}=f(ji,k)\\
f(i,j)f(j,k)(f(k,i)+f(j,i)f(i,k))&=f(j,i)f(i,k)(f(k,j)+f(i,j)f(j,k))\\
f(i,j)f(j,k)f(k,i)&=f(j,i)f(k,j)f(i,k).
}
Thus the assertion is fulfilled from $B(i,j)=f(i,j)/f(j,i)$.
\smallskip
\par
Finally we show the assertion of this lemma. For each $k\in T_{0}$, it follows from the form (\ref{eq:mue=expa}) with $k_{1}=k$ that the coefficient of $\mu^{S(k)}$ is equal to $f(j_{1}j_{2}\cdots j_{m_{0}-1},k)$, where we write $j_{1},\dots, j_{m_{0}-1}=k_{2}\cdots k_{m_{0}}$. Therefore it suffices to show that for any $1\leq n<m_{0}$ and for any disjoint elements $j_{0},j_{1},\dots, j_{n}\in T_{0}$
\alil
{
\textstyle f(j_{1}\cdots j_{n},j_{0})=(1+\sum_{i=1}^{n}B(j_{0},j_{i}))^{-1}.\label{eq:fe=3_2}
}
We prove it by induction for $n$. When $n=1$, 
\ali
{
f(j_{1},j_{0})=\frac{f(j_{1},j_{0})}{f(j_{1},j_{0})+f(j_{0},j_{1})}=\frac{1}{1+\frac{f(j_{0},j_{1})}{f(j_{1},j_{0})}}=\frac{1}{1+B(j_{0},j_{1})}
}
by $f(j_{1},j_{0})+f(j_{0},j_{1})=1$.
Therefore the assertion is valid. Assume that $(\ref{eq:fe=3_2})$ holds for some $1\leq n<m_{0}-1$. The equation (\ref{eq:fe=0_2}) yields
\ali
{
f(j_{1}\cdots j_{n+1},j_{0})=&\left(\frac{f(j_{0},j_{2}\cdots j_{n+1})+f(j_{1},j_{2}\cdots j_{n+1})f(j_{2}\cdots j_{n+1},j_{0})}{f(j_{1},j_{2}\cdots j_{n+1})f(j_{2}\cdots j_{n+1},j_{0})}\right)^{-1}\\
=&\left(1+\frac{1-f(j_{2}\cdots j_{n+1},j_{0})}{(1-f(j_{2}\cdots j_{n+1},j_{1}))f(j_{2}\cdots j_{n+1},j_{0})}\right)^{-1}\\
=&\left(1+\frac{1}{1-f(j_{2}\cdots j_{n+1},j_{1})}\frac{1-f(j_{2}\cdots j_{n+1},j_{0})}{f(j_{2}\cdots j_{n+1},j_{0})}\right)^{-1}\\
=&\left(1+\frac{1}{f(j_{2}\cdots j_{n+1},j_{1})\left(\frac{1}{f(j_{2}\cdots j_{n+1},j_{1})}-1\right)}\left(\frac{1}{f(j_{2}\cdots j_{n+1},j_{0})}-1\right)\right)^{-1}\\
=&\left(1+\frac{(1+\sum_{i=2}^{n+1}B(j_{1},j_{i}))\sum_{i=2}^{n+1}B(j_{0},j_{i})}{\sum_{i=2}^{n+1}B(j_{1},j_{i})}\right)^{-1}\\
=&\left(1+\frac{(1+\sum_{i=2}^{n+1}B(j_{1},j_{i}))B(j_{0},j_{1})\sum_{i=2}^{n+1}B(j_{1},j_{i})}{\sum_{i=2}^{n+1}B(j_{1},j_{i})}\right)^{-1}\\
=&\left(1+B(j_{0},j_{1})+\sum_{i=2}^{n+1}B(j_{0},j_{i})\right)^{-1}=\left(1+\sum_{i=1}^{n+1}B(j_{0},j_{i})\right)^{-1}
}
by using the equation $B(i,j)B(j,k)=B(i,k)$. Hence a proof is complete.
\proe
Recall the notation $\mathscr{S}_{0}$ in Section \ref{sec:intro}. Let $\mathscr{S}_{1}:=(S/\!\!\leftrightarrow)\setminus \mathscr{S}_{0}$.
We take nonempty subset $Q\subset S$. We write
\ali
{
\mathscr{S}_{0}(Q):=&\{U\in \mathscr{S}_{0}\,:\, U\cap Q\neq \emptyset\}\\
\mathscr{S}_{1}(Q):=&\{U\in \mathscr{S}_{1}\,:\, U\cap Q\neq \emptyset\}\\
\mathscr{S}(Q):=&\mathscr{S}_{0}(Q)\cup \mathscr{S}_{1}(Q).
}
\prop
{\label{prop:mueS1->0}
Assume that conditions (A.1)-(A.3) are satisfied. Take nonempty subset $Q\subset S$ with $\#\mathscr{S}(Q)<+\infty$ and $\mathscr{S}_{0}=\mathscr{S}_{0}(Q)$. Assume also that $g_{\e}^{Q\cap U}\geq c_{\adr{ubge}}\chi_{Q\cap U}$ for each $U\in \mathscr{S}(Q)$ for some $c_{\adr{ubge}}>0$. Then $\mu^{Q}_{\e}$ has a weak convergent sequence and the limit has a form of convex combination of $\{\mu^{Q\cap U}(U,\cd)\}_{U\in \mathscr{S}_{0}}$. In particular, for each $U\in \mathscr{S}_{1}$, $\mu^{Q}_{\e}(\Si_{U\cap Q})\to 0$ as $\e\to 0$. 
}
\pros
We define a $\#\mathscr{S}(Q)\times \#\mathscr{S}(Q)$ nonnegative matrix $E_{\e}^{Q}=(E_{\e}^{Q}(UV))_{UV}$ indexed by $\mathscr{S}(Q)$ with
\ali
{
E_{\e}^{Q}(UV)=\frac{\nu_{\e}^{Q\cap U}(\LR_{\e}[Q,S\setminus Q,\lam_{\e}]g_{\e}^{Q\cap V})}{\nu_{\e}^{Q\cap U}(g_{\e}^{Q\cap U})}.
}
Note that $E_{\e}^{Q}1=\lam_{\e}1$ for all $\e>0$, where $1$ means a vector whose entries are all $1$. Denote by $\gamma_{\e}=(\gamma_{\e}(U))_{U}$ the left Perron eigenvector of $E_{\e}^{Q}$. Then
\ali
{
\mu_{\e}^{Q}=\sum_{U\in \mathscr{S}(Q)}\gamma_{\e}(U)\mu_{\e}^{Q\cap U}.
}
Now we consider $\lam_{\e}-E_{\e}^{Q}(UU)$ for each $U\in \mathscr{S}(Q)$. Note that
\ali
{
1=\nu_{\e}^{Q\cap U}(1)\geq \nu_{\e}^{Q\cap U}(g_{\e}^{Q\cap U})\geq c_{\adr{ubge}}\nu_{\e}^{Q\cap U}(\Si_{Q\cap U})=c_{\adr{ubge}}>0.
}
Recall that for $U\in \mathscr{S}(Q)$ with $B_{UU}\neq [0]$, $g^{Q}(U,\cd)$ is defined by the eigenfunction of $\LR[Q\cap U,U\setminus Q,\lam]$. For $U\in \mathscr{S}(Q)$ with $B_{UU}=[0]$, we let $g^{Q\cap U}(U,\cd)=\chi_{Q\cap U}$. Then there exists a constant $c_{\adl{dgU}}>0$ such that $c_{\adr{dgU}}\chi_{Q\cap U}\leq g^{Q\cap U}(U,\cd)$. 

If $U\in \mathscr{S}_{0}(Q)$
\ali
{
\lam_{\e}-E_{\e}^{Q}(UU)=&\frac{\lam_{\e}\nu_{\e}^{Q\cap U}(g_{\e}^{Q\cap U})}{\nu_{\e}^{Q\cap U}(g_{\e}^{Q\cap U})}-\frac{\nu_{\e}^{Q\cap U}(\LR_{\e}[Q\cap U,S\setminus Q,\lam_{\e}]g_{\e}^{Q\cap U})}{\nu_{\e}^{Q\cap U}(g_{\e}^{Q\cap U})}\\
=&\frac{\nu_{\e}^{Q\cap U}((\LR_{\e}[Q\cap U,S\setminus (Q\cap U),\lam_{\e}]-\LR_{\e}[Q\cap U,S\setminus Q,\lam_{\e}])g_{\e}^{Q\cap U})}{\nu_{\e}^{Q\cap U}(g_{\e}^{Q\cap U})}\\
\leq &c_{\adr{dgU}}^{-2}\nu_{\e}^{Q\cap U}((\LR_{\e}[Q\cap U,S\setminus (Q\cap U),\lam_{\e}]-\LR_{\e}[Q,S\setminus Q,\lam_{\e}])g^{Q\cap U}(U,\cd))\\
&\qqquad (\because \LR_{\e}[Q\cap U,S\setminus (Q\cap U),\lam_{\e}]-\LR_{\e}[Q\cap U,S\setminus Q,\lam_{\e}] \text{ is positive})\\
\leq &c_{\adr{dgU}}^{-2}\nu_{\e}^{Q\cap U}((\lam_{\e}\IR-\LR_{\e}[Q,S\setminus Q,\lam_{\e}])g^{Q\cap U}(U,\cd)).
}
By using Proposition \ref{prop:conv_PC_Rue2} and the fact $\|g^{Q\cap U}(U,\cd)\|_{\infty}=1$
\ali
{
&\|(\lam_{\e}\IR-\LR_{\e}[Q\cap U,S\setminus Q,\lam_{\e}])g^{Q\cap U}(U,\cd)-(\lam\IR-\LR[Q\cap U,U\setminus Q,\lam])g^{Q\cap U}(U,\cd)\|_{\infty}\\
\leq &|\lam_{\e}-\lam|+\|\LR_{\e}[Q\cap U,S\setminus Q,\lam_{\e}]-\LR[Q\cap U,U\setminus Q,\lam]\|_{\infty}\to 0
}
as $\e\to 0$. In addition to convergence $\nu_{\e}^{Q\cap U}\to \nu^{Q\cap U}(U,\cd)$ weakly, Therefore
\ali
{
&c_{\adr{dgU}}^{-2}\nu_{\e}^{Q\cap U}((\lam_{\e}\IR-\LR_{\e}[Q\cap U,S\setminus Q,\lam_{\e}])g^{Q\cap U}(U,\cd))\\
\to &c_{\adr{dgU}}^{-2}\nu^{Q\cap U}(U,\cd)(\lam\IR-\LR[Q\cap U,U\setminus Q,\lam])g^{Q\cap U}(U,\cd)=0.
}
We get $\lam_{\e}-E_{\e}^{Q}(UU)\to 0$.

If $U\in \mathscr{S}_{1}(Q)$ with $B_{UU}\neq [0]$, we let $\lam(U)$ the spectral radius of $\LR_{UU}$. Then we see $\LR[Q\cap U,U\setminus Q,\lam(U)]g^{Q\cap U}(U,\cd)=\lam(U)g^{Q\cap U}(U,\cd)$. Moreover
\alil
{
\lam_{\e}-E_{\e}^{Q}(UU)\geq &c_{\adr{dgU}}\nu_{\e}^{Q\cap U}((\LR_{\e}[Q\cap U,S\setminus (Q\cap U),\lam_{\e}]-\LR_{\e}[Q\cap U,S\setminus Q,\lam_{\e}])g^{Q\cap U}(U,\cd))\nonumber\\
\geq &c_{\adr{dgU}}\nu_{\e}^{Q\cap U}((\lam_{\e}\IR-\LR_{\e}[Q\cap U,S\setminus Q,\lam_{\e}])g^{Q\cap U}(U,\cd)).\label{eq:le-EPeU}
}
Let $c_{\adl{lU2}}=(\lam-\lam(U))c_{\adr{dgU}}>0$. By Proposition \ref{prop:conv_PC_Rue1_v2}, the operator $(\lam_{\e}\IR-\LR_{\e}[Q\cap U,S\setminus Q,\lam_{\e}])g^{Q\cap U}(U,\cd)$ converges to $(\lam\IR-\LR[Q\cap U,U\setminus Q,\lam])g^{Q\cap U}(U,\cd)$ uniformly. Then there exists $\e_{0}>0$ such that for any $0<\e<\e_{0}$
\ali
{
\|(\lam_{\e}\IR-\LR_{\e}[Q\cap U,S\setminus Q,\lam_{\e}])g^{Q\cap U}(U,\cd)-(\lam\IR-\LR[Q\cap U,U\setminus Q,\lam])g^{Q\cap U}(U,\cd)\|_{\infty}\leq c_{\adr{lU2}}/2.
}
Now we notice that for any $\om\in \Si_{Q\cap U}$
\ali
{
&(\lam_{\e}\IR-\LR_{\e}[Q\cap U,S\setminus Q,\lam_{\e}])g^{Q\cap U}(U,\om)\\
=&(\lam_{\e}\IR-\LR_{\e}[Q\cap U,S\setminus Q,\lam_{\e}])g^{Q\cap U}(U,\om)-(\lam\IR-\LR[Q\cap U,U\setminus Q,\lam])g^{Q\cap U}(U,\om)\\
&+(\lam\IR-\LR[Q\cap U,U\setminus Q,\lam])g^{Q\cap U}(U,\om)\\
\geq&(\lam\IR-\LR[Q\cap U,U\setminus Q,\lam])g^{Q\cap U}(U,\om)-c_{\adr{lU2}}/2\\
\geq&(\lam\IR-\LR[Q\cap U,U\setminus Q,\lam(U)])g^{Q\cap U}(U,\om)-c_{\adr{lU2}}/2\\
=&(\lam\IR-\lam(U))g^{Q\cap U}(U,\cd)-c_{\adr{lU2}}/2\\
\geq&c_{\adr{lU2}}-c_{\adr{lU2}}/2=c_{\adr{lU2}}/2>0.
}
Thus (\ref{eq:le-EPeU}) yields
\ali
{
\lam_{\e}-E_{\e}^{Q}(UU)\geq c_{\adr{dgU}}c_{\adr{lU2}}/2>0
}
for all $0<\e<\e_{0}$. We get the positively $\liminf_{\e\to 0}\lam_{\e}-E_{\e}^{Q}(UU)>0$.
\smallskip
\par
If $U\in \mathscr{S}_{1}(Q)$ with $B_{UU}=[0]$, then $\LR_{\e}[Q\cap U,S\setminus Q,\lam_{\e}]$ converges to $\LR[Q\cap U,S\setminus Q,\lam]=O$ with respect to $\|\cd\|_{\infty}$. Moreover
\ali
{
\lam_{\e}-E_{\e}^{Q}(UU)\geq &c_{\adr{dgU}}\nu_{\e}^{Q\cap U}((\lam_{\e}\IR-\LR_{\e}[Q\cap U,S\setminus Q,\lam_{\e}])\chi_{Q\cap U})\\
\geq&\nu_{\e}^{Q\cap U}((\lam_{\e}-\|\LR_{\e}[Q\cap U,S\setminus Q,\lam_{\e}]\|_{\infty})\chi_{Q\cap U})\\
=&\lam_{\e}-\|\LR_{\e}[Q\cap U,S\setminus Q,\lam_{\e}]\|_{\infty}\to \lam
}
Thus $\liminf_{\e\to 0}\lam_{\e}-E_{\e}^{Q}(UU)\geq \lam>0$. As result, 
\ali
{
\lam_{\e}-E_{\e}^{Q}(UU)
\case
{
\to 0&(U\in \mathscr{S}_{0}(Q))\\
\geq c_{\adr{ulE}}&(U\in \mathscr{S}_{1}(Q))\\
}
}
for some $c_{\adl{ulE}}>0$. By $E_{\e}^{Q}1=\lam_{\e}1$, for each $U\in \mathscr{S}_{0}$
\ali
{
\sum_{V\in \mathscr{S}(Q)\,:\,V\neq U}E_{\e}^{Q}(UV)=\lam_{\e}-E_{\e}^{Q}(UU)\to 0
}
and therefore $E_{\e}^{Q}(UV)\to 0$ for any $U\in \mathscr{S}_{0}(Q)$ and $V\in \mathscr{S}(Q)$. Choose any positive sequence $(\e_{n})$ with $\inf_{n}\e_{n}=0$. Since $E_{\e}^{Q}$ is finite matrix, there exists a subsequence $(\e^{\p}_{n})$ of $(\e_{n})$ such that $E_{\e^\p_{n}}^{Q}$ converges to a matrix $M$ and $\gamma_{\e^{\p}_{n}}=(\gamma_{\e^{\p}_{n}}(U))_{U}$ converges to $\gamma=(\gamma(U))_{U}$. Then $\gamma$ is the corresponding nonnegative left eigenvector of the eigenvalue $\lam$ of $M$. The matrix $M$ has the form
\ali
{
{}^{t}M_{0}MM_{0}=\left[
\begin{array}{c|c}
M_{1}&O\\
\hline
*&M_{2}\\
\end{array}\right]\quad\text{ with }
M_{1}=\left[\mat
{
\lam& & O\\
 &\ddots& \\
O & & \lam\\
}\right]
}
for suitable permutation matrix $M_{0}$. We note that $M_{2}$ is a triangular matrix. Indeed, letting $V=S\setminus Q$, we notice $r(\LR_{VV})\leq \max\{r(\LR_{V^\p V^\p})\,:\,V^\p\in \mathscr{S}(Q)\}<\lam$. By Proposition \ref{prop:conv_PC_Rue1_v2}, for each $U_{1},U_{2}\in \mathscr{S}_{1}(Q)$
\ali
{
E_{\e}^{Q}(U_{1}U_{2})\leq &c_{\adr{dgU}}^{-1}\|\chi_{Q\cap U_{1}}\LR_{\e}[Q,S\setminus Q,\lam_{\e}]\chi_{Q\cap U_{2}}\|_{\infty}\\
\to&c_{\adr{dgU}}^{-1}\|\chi_{Q\cap U_{1}}\LR[Q,S\setminus Q,\lam]\chi_{Q\cap U_{2}}\|_{\infty}.
}
We note that $\chi_{Q\cap U_{1}}\LR[Q,S\setminus Q,\lam]\chi_{Q\cap U_{2}}\neq O$ if and only if there exists a $B$-admissible word $w=w_{1}\dots w_{n}$ such that $w_{1}\in Q\cap U_{2}$ and $w_{n}\in Q\cap U_{1}$. Therefore, choose any $U_{1},U_{2},\cdots,U_{k+1}\in \mathscr{S}_{1}(Q)$ with $U_{k+1}=U_{1}$ and $U_{1}\neq U_{i}$ for some $i$. Since there is no $B$-admissible word $w$ so that $w_{j_{i}}\in U_{i}$ $(1\leq i\leq k+1$ for some $1\leq j_{1}<j_{2}<\cdots <j_{k+1}$, we have 
\ali
{
\prod_{i=1}^{k}E_{\e}^{Q}(U_{i}U_{i+1})\to \prod_{i=1}^{k}M_{2}(U_{i}U_{i+1})=0
}
as $\e\to 0$ running through $\{\e_{n}^\p\}$. Since the spectral radius of $M_{2}$ is less than $\lam$, we obtain $\gamma(U)=0$ for all $U\in \mathscr{S}_{1}(Q)$. By Proposition \ref{prop:conv_PC_mue}, we have that for any $f\in F_{b}^{1}(X)$
\ali
{
\mu_{\e}(f)=\sum_{U\in \mathscr{S}(Q)}\gamma_{\e}(U)\mu_{\e}^{Q\cap U}(f)\to \sum_{U\in \mathscr{S}_{0}(Q)}\gamma(U)\mu^{Q\cap U}(f)
}
by using the fact $\gamma_{\e}(U)\mu^{Q\cap U}(f)\to 0$ for each $U\in \mathscr{S}_{0}(Q)$. Hence the proof is complete.
\proe
Let $Q\subset S$ be nonempty subset with $Q\cap U\neq \emptyset$ for all $U\in \mathscr{S}_{0}$, 
Recall that for $i,j\in T_{0}$ with $i\neq j$, $c_{\e}(Q,i,j)$ is the spectral radius of $\LR_{\e}[Q\cap U(i),S\setminus ((U(i)\cup U(j))\cap Q),\lam_{\e}]$, and $\delta_{\e}(Q,k)$ is defined in (\ref{eq:deQi=}). We let $b_{\e}(Q,i,j)$ be $b_{\e}(i : S\setminus \{i,j\})$.
\prop
{\label{prop:dek/tdek->1}
Assume that conditions (A.1)-(A.3) are satisfied. Let $Q\subset \bigcup\mathscr{S}_{0}$ be a nonempty finite set with $Q\cap U\neq \emptyset$ for all $U\in \mathscr{S}_{0}$. Then $c_{\e}(Q,i,j)/b_{\e}(Q,i,j)$ converges to $1$. Consequently, $\delta_{\e}(Q,i)/\ti{\delta}_{\e}(Q,i)\to 1$ as $\e\to 0$.
}
\pros
Recall the notation $\mathscr{S}_{0}=\{U(1),U(2),\dots,U(m_{0})\}$. We put $Q(i):=Q\cap U(i)$ for $i$. We write $T_{0}=\{k_{1},k_{2},\dots, k_{m_{0}}\}$. Let $\ti{Q}(k_{j})=\bigcup_{l=j+1}^{m_{0}}Q(k_{l})$ for $j=1,2,\dots, m_{0}-1$. For the sake of convenience, we set $\ti{Q}(k_{0}):=\bigcup_{l=1}^{m_{0}}Q(k_{l})=Q$. Recall the $2\times 2$ matrix (\ref{eq:Ej=...}). For each $j\in T_{0}$, we define a $2\times 2$ matrix $E^{j}_{\e}=(M^{j}_{\e}(kk^\p))$ by
\ali
{
&\matII
{
\mu_{\e}^{Q(k_{j})}(\tLR_{\e}[Q(k_{j}),V,\lam_{\e}]1)&\mu_{\e}^{Q(k_{j})}(\tLR_{\e}[\ti{Q}(k_{j-1}),V ,\lam_{\e}]_{Q(k_{j})\ti{Q}(k_{j})}1)\\
\mu_{\e}^{\ti{Q}(k_{j})}(\tLR_{\e}[\ti{Q}(k_{j-1}),V,\lam_{\e}]_{\ti{Q}(k_{j})Q(k_{j})}1)&\mu_{\e}^{\ti{Q}(k_{j})}(\tLR_{\e}[\ti{Q}(k_{j}),V,\lam_{\e}]1)\\
}\\
=&\matII
{
\di\frac{\nu_{\e}^{Q}(\LR_{\e}[Q(k_{j}),V,\lam_{\e}]g_{\e}^{Q})}{\nu_{\e}^{Q}(g_{\e}^{Q}\chi_{Q(k_{j})})}
&\di\frac{\nu_{\e}^{Q}(\chi_{Q(k_{j})}\LR_{\e}[\ti{Q}(k_{j-1}),V ,\lam_{\e}](\chi_{\ti{Q}(k_{j})}g_{\e}^{Q}))}{\nu_{\e}^{Q}(g_{\e}^{Q}\chi_{Q(k_{j})})}\\
\di\frac{\nu_{\e}^{Q}(\chi_{\ti{Q}(k_{j})}\LR_{\e}[\ti{Q}(k_{j-1}),V,\lam_{\e}](\chi_{Q(k_{j})}g_{\e}^{Q}))}{\nu_{\e}^{Q}(g_{\e}^{Q}\chi_{\ti{Q}(k_{j})})}&\di\frac{\nu_{\e}^{Q}(\LR_{\e}[\ti{Q}(k_{j}),V,\lam_{\e}]g_{\e}^{Q})}{\nu_{\e}^{Q}(g_{\e}^{Q}\chi_{\ti{Q}(k_{j})})}.
}
}
with $V=S\setminus \ti{Q}(k_{j-1})$. Denote by $\delta^{j}_{\e}=(\delta^{j}_{\e}(i))_{i=1}^{2}$ the left Perron eigenvalue with $\sum_{i=1}^{2}\delta^{j}_{\e}(i)=1$ of the matrix $M^{j}_{\e}$. 
Since the right Perron eigenvector of the matrix $E_{\e}^{j}$ is a vector whose entires are $1$, we see $E_{\e}^{j}(11)+E_{\e}^{j}(12)=\lam_{\e}$. Put $V=S\setminus Q$. Notice the equation
\ali
{
S\setminus \ti{Q}(k_{j-1})=V\cup \bigcup_{l=1}^{j-1}Q(k_{l}).
}
Therefore, by putting $k_{j}=i$ and $P=\{k_{1},\dots,k_{j-1}\}$ 
\alil
{
b_{\e}(i:P)=&\lam_{\e}-\frac{\nu_{\e}^{Q}(\LR_{\e}[Q(i),V\cup \bigcup_{k\in P}Q(k),\lam_{\e}]g_{\e}^{Q})}{\nu_{\e}^{Q}(g_{\e}^{Q}\chi_{Q(i)})}\nonumber\\
=&\lam_{\e}-E_{\e}^{j}(11)\nonumber\\
=&E_{\e}^{j}(12)=\frac{\nu_{\e}^{Q(i)}(\chi_{Q(i)}\LR_{\e}[\ti{Q}(k_{j-1}),V\cup \bigcup_{l=1}^{j-1}Q(i),\lam_{\e}](\chi_{\ti{Q}(i)}g_{\e}))}{\nu_{\e}^{Q(i)}(g_{\e}\chi_{Q(i)})}.\label{eq:lem:be/ceto1}
}
On the other hand, recall the spectral radius $c_{\e}(Q,i,j)$ of the operator $\LR_{\e}[Q(i),S\setminus (Q(i)\cup Q(j)),\lam_{\e}]$. Let $i\in T_{0}$. Note that $A_{Q(i)\times Q(i)}$ may be not irreducible. We take a finite subset $\ti{Q}(i)$ of $S$ so that $Q(i)\subset \ti{Q}(i)\subset S(i)$ and $A|_{\ti{Q}(i)\times \ti{Q}(i)}$ is irreducible (see Proposition \ref{prop:irresubmat}). Notice that $r((\LR_{\e})_{\ti{Q}(i)\ti{Q}(i)}$ is positive.
We take the Perron eigenvector $\ti{\nu}_{\e}\in M(\SiAp)$ of $\LR_{\e}[\ti{Q}(i),(V\setminus \ti{Q}(i))\cup \bigcup_{l=1}^{j-1}Q(k_{l}),\lam_{\e}]^{*}$. Then we can show that $\ti{\nu}_{\e}^{Q(i)}$ is the Perron eigenvector of $\LR_{\e}[Q(i),(V\setminus \ti{Q}(i))\cup (\ti{Q}(i)\setminus Q(i))\cup \bigcup_{l=1}^{j-1}Q(k_{l}),\lam_{\e}]^{*}=\LR_{\e}[Q(i),(V\cup \bigcup_{l=1}^{j-1}Q(k_{l}),\lam_{\e}]^{*}$ by using the proof of Proposition \ref{prop:spec_PC_Rue}. Moreover, we can check that $\ti{\nu}_{\e}^{Q(i)}$ converges to $\nu^{Q(i)}(U(i),\cd)$ (see the proof of Corollary \ref{prop:ex_PSC_PC_per}(3)). Put $U^\p=\bigcup_{l=j}^{m_{0}}Q(k_{l})$ and $U^\pp=\bigcup_{l=1}^{j-1}Q(k_{l})$. We obtain
\ali
{
c_{\e}(i:P)\ti{\nu}_{\e}^{Q(i)}(g_{\e}^{U^\p})&=\ti{\nu}_{\e}^{Q(i)}(\lam_{\e}g_{\e}^{U^\p}-r(\LR_{\e}[Q(i),V\cup U^\pp,\lam_{\e}])g_{\e}^{U^\p})\\
&=\ti{\nu}_{\e}^{Q(i)}((\LR_{\e}[U^\p,V\cup U^\pp,\lam_{\e}]-\LR_{\e}[Q(i),V\cup U^\pp,\lam_{\e}])g_{\e}^{U^\p})\\
&=\ti{\nu}_{\e}^{Q(i)}(\LR_{\e}[U^\p,V\cup U^\pp,\lam_{\e}]_{Q(i)\ti{Q}(i)}g_{\e}^{U^\p})
}
by using $U^\p\setminus Q(i)=\ti{Q}(i)$ and $\supp\,\ti{\nu}_{\e}^{Q(i)}\subset \Si_{Q(i)}$. Thus
\alil
{
c_{\e}(i:P)=\frac{\ti{\nu}_{\e}^{Q(i)}(\chi_{Q(i)}\LR_{\e}[U^\p,V\cup U^\pp,\lam_{\e}](\chi_{\ti{Q}(i)}g_{\e}))}{\ti{\nu}_{\e}^{Q(i)}(g_{\e}\chi_{Q(i)})}\label{eq:lem:be/ceto2}
}
by the form $g_{\e}^{U^\p}=g_{\e}\chi_{U^\p}/\|g_{\e}\chi_{U^\p}\|_{\infty}$.

Let $\xi_{\e}=\chi_{Q(i)}\LR_{\e}[U^\p,V\cup U^\pp,\lam_{\e}](\chi_{\ti{Q}(i)}g_{\e})/\|\chi_{Q(i)}\LR_{\e}[U^\p,V\cup U^\pp,\lam_{\e}](\chi_{\ti{Q}(i)}g_{\e})\|_{\infty}$. By Proposition \ref{prop:prop_Lamc0}(4), $\xi_{\e}$ is in $\Lambda_{c_{\adr{cL}}}$ is valid for any $\e>0$. 
Choose any positive sequence $(\e(n))$ with $\lim_{n\to \infty}\e(n)=0$ so that the number $c_{\e(n)}(i:P)/b_{\e(n)}(i:P)$ converges in $[0,\infty]$. We take a subsequence $(\e^{\p}(n))$ of $(\e(n))$ such that $\xi_{\e}$ converges to a function $\xi$ for each point in $X$. We remark that $\xi\neq 0$ on $\supp\, \nu^{Q(i)}(U(i),\cd)$. Indeed, since $Q(i)$ is finite, $\xi(\om)>0$ for some $\om\in \Si_{Q(i)}$. For any $\up\in \supp\, \nu^{Q(i)}(U(i),\cd)\cap [\om_{0}]$, we have the inequality $0<\xi(\om)\leq e^{c_{14}\theta}\xi(\up)$ from $\xi\in \Lambda_{c_{14}}$. In particular, $\nu^{Q(i)}(U(i),\xi)>0$ is satisfied. Now we recall the notation $g_{\e}^{Q(i)}=g_{\e}\chi_{Q(i)}/\|g_{\e}\chi_{Q(i)}\|_{\infty}$ and $\nu_{\e}^{Q(i)}=\nu_{\e}(\cd|\Si_{Q(i)})$.
The equations (\ref{eq:lem:be/ceto1}) and (\ref{eq:lem:be/ceto2}) yield
\ali
{
\frac{c_{\e}(i:P)}{b_{\e}(i:P)}=\frac{\ti{\nu}_{\e}^{Q(i)}(\xi_{\e})}{\nu_{\e}^{Q(i)}(\xi_{\e})}\frac{\nu_{\e}^{Q(i)}(g_{\e}^{Q(i)})}{\ti{\nu}_{\e}^{Q(i)}(g_{\e}^{Q(i)})}\to \frac{\nu^{Q(i)}(U(i),\xi)}{\nu^{Q(i)}(U(i),\xi)}\frac{\nu^{Q(i)}(U(i),g^{Q(i)}(U(i),\cd))}{\nu^{Q(i)}(U(i),g^{Q(i)}(U(i),\cd))}=1
}
as $\e\to 0$ running through $\{\e^{\p}(n)\}$. This fact does not depend on how to take a sequence $\{\e(n)\}$. Hence we obtain the assertion.
\proe
\section{Proof of main results}\label{sec:proof}
\subsection{Proof of Theorem \ref{th:main}}\label{sec:proof_main2}
\pros
Next we prove (1) $\Rightarrow$ (2). We assume that $\mu_{\e}\,:\,F_{b,0}^{1}(X)\to \R$ converges to a measure $\mu$ weakly. Choose any $Q\subset S$ so that $Q$ is nonempty finite set with $Q\cap U\neq \emptyset$ for all $U\in \mathscr{S}_{0}$. We have $\mu_{\e}(\chi_{Q})\to \mu(\chi_{Q})=:a(Q)$. If $a(Q)>0$, then for any $f\in F_{b}^{1}(X)$, $f\chi_{Q}\in F_{b,0}^{1}(X)$ and
\ali
{
\mu_{\e}^{Q}(f)=\frac{\mu_{\e}(f\chi_{Q})}{\mu_{\e}(\chi_{Q})}\to \frac{\mu(f\chi_{Q})}{\mu(\chi_{Q})}.
}
Therefore Proposition \ref{prop:mueS1->0}, $\delta_{\e}(Q,k)$ converges for all $k$. Thus the assertion is valid.
\smallskip
\\
Next we show (2) $\Rightarrow$ (3). Choose any sequence $Q_{n}\subset S$ so that $Q_{n}$ is finite, $Q_{n}\subset Q_{n+1}$, $\bigcup_{n}Q_{n}=S$, and $Q_{n}\cap U\neq \emptyset$ for all $U\in \mathscr{S}_{0}$. By (2), $\mu_{\e}(\chi_{Q_{n}})\to a(Q_{n})$ as $\e\to 0$ for each $n$. By $\mu_{\e}(\chi_{Q_{n}})\leq \mu_{\e}(\chi_{Q_{n+1}})$, we have $a(Q_{n})\leq a(Q_{n+1})$. Therefore, $\lim_{n\to \infty}a(Q_{n})=a$ exists. Moreover, assume $a>0$. Then there exists $c_{\adl{cQn}}>0$ such that $\mu_{\e}(\chi_{Q_{n}})\geq c_{\adr{cQn}}$ for any large $n$ and for any small $\e>0$. Then $\liminf_{\e\to 0}\mu_{\e}(\chi_{Q_{n}})>0$. Choose any $1\leq k\leq m_{0}$, we have $\delta_{\e}(Q_{n},k)\to \delta(Q_{n},k)$ by the assumption of (2). Thus (3) is fulfilled.
\smallskip
\\
Next we show (3) $\Rightarrow$ (1). Take a finite subset $Q_{n}\subset S$ so that $Q_{n}\cap U\neq \emptyset$ for all $U\in \mathscr{S}_{0}$, $Q_{n}\subset Q_{n+1}$ and $\bigcup_{n}Q_{n}=S$. We show $\mu_{\e}\to \mu$ weakly in the sense of $F_{b,0}(X)$ for some measure $\mu$.
Fix $f\in F_{b,0}(X)$. Choose any $0<\eta<1/2$. By the definition of $F_{b,0}(X)$, there exists a finite subset $Q\subset S$ such that $|f|\leq \eta$ on $\Si_{S\setminus Q}$. Therefore, we see $|f|\leq \eta$ on $\Si_{S\setminus Q_{n}}$ for any $n\geq n_{0}$ for some $n_{0}\geq 1$. Then
\alil
{
\mu_{\e}(f)=&\mu_{\e}(f\chi_{Q_{n}})+\mu_{\e}(f\chi_{S\setminus Q_{n}})\quad \text{ and }\quad|\mu_{\e}(f\chi_{S\setminus Q_{n}})|\leq \eta.\label{eq:muef=...}
}
Let $a:=\lim_{\e\to 0,n\to \infty}\mu_{\e}(\chi_{Q_{n}})$. If $a=0$, then $|\mu_{\e}(f\chi_{Q_{n}})|\leq \|f\|_{\infty}|\mu_{\e}(\chi_{Q_{n}})\to 0$ as $n\to \infty$ and $\e\to 0$ and thus $\limsup_{\e\to 0}|\mu_{\e}(f)|\leq \eta$ for any small $\eta>0$. This means that $\mu_{\e}(f)$ converges to zero for any $f\in F_{b,0}(X)$.

On the other hand, we assume $a>0$. There exist $n_{1}\geq 1$ and $\e_{0}>0$ such that 
\ali
{
\mu_{\e}(\chi_{Q_{n}})\geq a/2\quad\text{ and }\quad|\mu_{\e}(\chi_{Q_{n}})-a|<\eta
}
for any $0<\e<\e_{0}$ and $n\geq n_{1}$. Therefore
\ali
{
\mu_{\e}(f\chi_{Q_{n}})=\mu_{\e}(\chi_{Q_{n}})\mu_{\e}^{Q_{n}}(f)
=&a\mu_{\e}^{Q_{n}}(f)+(\mu_{\e}(\chi_{Q_{n}})-a)\mu_{\e}^{Q_{n}}(f)\\
|(\mu_{\e}(\chi_{Q_{n}})-a)\mu_{\e}^{Q_{n}}(f)|\leq &\eta\|f\|.
}
Moreover, by $\mu(U,1)=1$, there exists $n_{2}\geq 1$ such that for any $n\geq n_{2}$ and for any $U\in \mathscr{S}_{0}$
\ali
{
|\mu(U,\chi_{Q_{n}})-1|\leq \eta.
}
By the assumption, $\delta_{\e}(Q_{n},k)$ converges as $\e\to 0$ for each $1\leq k\leq m_{0}$. By Proposition \ref{prop:mueS1->0}, for each $n\geq 1$, $\mu_{\e}^{Q_{n}}(f)$ converges to $\sum_{U\in \mathscr{S}_{0}}\delta(Q_{n},U)\mu^{Q_{n}}(U,f)$ for a probability vector $(\delta(Q_{n},U))_{U}$. Now we show convergence of $\delta(Q_{n},U)$ as $n\to \infty$. Choose any $U\in \mathscr{S}_{0}$, $n\geq 1$ and $b\in Q_{n}\cap U$. Then  for any $m\geq n$ , $b\in Q_{m}$ and
\ali
{
\frac{\mu_{\e}([b])}{\mu_{\e}(\chi_{Q_{m}})}=\mu_{\e}^{Q_{m}}(\chi_{[b]})\to \delta(Q_{m},U)\mu^{Q_{m}}(U,[b])=\delta(Q_{m},U)\frac{\mu(U,[b])}{\mu(U,Q_{m})}
}
as $\e\to 0$. By $\mu_{\e}(\chi_{Q_{m}})\leq \mu_{\e}(\chi_{Q_{m+1}})$, we have $\mu_{\e}^{Q_{m}}(\chi_{[b]})\geq \mu_{\e}^{Q_{m+1}}(\chi_{[b]})$ and
\ali
{
\frac{\delta(Q_{m},U)}{\mu(U,Q_{m})}\geq \frac{\delta(Q_{m+1},U)}{\mu(U,Q_{m+1})}.
}
Then $\lim_{m\to \infty}\delta(Q_{m},U)/\mu(U,Q_{m})$ converge to a number $\delta(U)$ and $\lim_{m\to \infty}\mu(U,Q_{m})=1$, we get
\ali
{
\delta(Q_{m},U)=
\frac{\delta(Q_{m},U)}{\mu(U,Q_{m})}\mu(U,Q_{m})\to \delta(U)
}
as $m\to \infty$. Thus there exists $n_{3}\geq 1$ such that for any $n\geq n_{3}$ and $U\in \mathscr{S}_{0}$
\ali
{
|\delta(Q_{m},U)-\delta(U)|\leq \eta.
}
Fix $n=\max\{n_{0},n_{1},n_{2},n_{3}\}$. By Proposition \ref{prop:mueS1->0}, there exists $\e_{1}>0$ such that
\ali
{
|\mu_{\e}^{Q_{n}}(f)-\sum_{U\in \mathscr{S}_{0}}\delta(Q_{n},U)\mu^{Q_{n}}(U,f)|\leq \eta
}
for any $0<\e<\e_{1}$. Furthermore
\ali
{
|\mu(U,f)-\mu^{Q_{n}}(U,f)|\leq& |\mu(U,f)-\mu(U,f\chi_{Q_{n}})|+|\mu(U,f\chi_{Q_{n}})|\left|1-\frac{1}{\mu(U,\chi_{Q_{n}})}\right|\\
\leq&\|f\|_{\infty}\eta+\|f\|_{\infty}\frac{\eta}{1-\eta}\leq 3\|f\|_{\infty}\eta.
}
Therefore
\ali
{
&|\delta(Q_{n},U)\mu^{Q_{n}}(U,f)-\delta(U)\mu(U,f)|\\
\leq& |\delta(Q_{n},U)\mu^{Q_{n}}(U,f)-\delta(U)\mu^{Q_{n}}(U,f)|+|\delta(U)\mu^{Q_{n}}(U,f)-\delta(U)\mu(U,f)|\\
\leq&\|f\|_{\infty}\eta+\delta(U)3\eta\|f\|_{\infty}\leq 4\eta\|f\|_{\infty}.
}
Consequently, we obtain that for any $0<\e<\min(\e_{0},\e_{1})$
\ali
{
&|\mu_{\e}(f)-a\sum_{U\in \mathscr{S}_{0}}\delta(U)\mu(U,f)|\\
\leq&|\mu_{\e}(f\chi_{S\setminus Q_{n}})|+|(\mu_{\e}(\chi_{Q_{n}})-a)\mu_{\e}^{Q_{n}}(f)|\\
&+|a\mu_{\e}^{Q_{n}}(f)-a\sum_{U\in \mathscr{S}_{0}}\delta(Q_{n},U)\mu^{Q_{n}}(U,f)|+a\sum_{U\in \mathscr{S}_{0}}|\delta(Q_{n},U)\mu^{Q_{n}}(U,f)-\delta(U)\mu(U,f)|\\
\leq&(1+\|f\|_{\infty}+a+4a(\#\mathscr{S}_{0})\|f\|_{\infty})\eta.
}
Hence $\mu_{\e}(f)$ converges.\\
(3) $\Rightarrow$ (4). The assertion directly follows.\\
(4) $\Rightarrow$ (3). Assume the assertion (4). Let $Q:=Q_{n}$ for a fixed $n$. We may assume that $Q\cap U\neq \emptyset$ for all $U\in \mathscr{S}_{0}$. Recall the equation $\mu_{\e}^{Q}=\sum_{U\in \mathscr{S}_{0}}\ti{\delta}_{\e}(Q,U)\mu_{\e}^{Q\cap U}$ and the fact $\ti{\delta}_{\e}(Q,U)/\delta_{\e}(Q,U)\to 1$ as $\e\to 0$. Since $1=\mu_{\e}^{Q}(1)=\sum_{U\in \mathscr{S}_{0}}\ti{\delta}_{\e}(Q,U)\mu_{\e}^{Q\cap U}(1)=\sum_{U\in \mathscr{S}_{0}}\ti{\delta}_{\e}(Q,U)$, we obtain
\ali
{
\sum_{U\in \mathscr{S}_{0}}\delta_{\e}(Q,U)=\frac{\sum_{U\in \mathscr{S}_{0}}\delta_{\e}(Q,U)}{\sum_{U\in \mathscr{S}_{0}}\tilde{\delta}_{\e}(Q,U)}\to 1.
}
Therefore
\ali
{
\mu_{\e}(\chi_{Q})=\frac{\sum_{U\in \mathscr{S}_{0}}\mu_{\e}(\chi_{Q})\delta_{\e}(Q,U)}{\sum_{U\in \mathscr{S}_{0}}\delta_{\e}(Q,U)}\to \sum_{U\in \mathscr{S}_{0}}\delta(n,U).
}
That is, $\mu_{\e}(\chi_{Q_{n}})$ converges to a number $a(n)$ for $\e\to 0$ for all $n$. If $a(n)>0$ then
\ali
{
\delta_{\e}(Q_{n},U)=\frac{\mu_{\e}(\chi_{Q_{n}})\delta_{\e}(Q_{n},U)}{\mu_{\e}(\chi_{Q_{n}})}\to \frac{\delta(n,U)}{a(n)}.
}
Thus the assertion (3) holds. Hence the proof is complete.
\proe
\subsection{Proof of Theorem \ref{th:main_ac}}\label{sec:proof_main1}
\pros
Choose any positive sequence $(\e(n))$ with $\e(n)\to 0$ and a sequence of finite subsets $\{Q_{k}\}$ of $S$ with $Q_{k}\cap U(i)\neq \emptyset$ for all $i=1,\dots, m_{0}$, $Q_{k}\subset Q_{k+1}$ and $\bigcup_{k}Q_{k}=S$. Consider the double sequence $\{\mu_{\e(n)}(\chi_{Q_{k}})\delta_{\e}(Q_{k},i)\}_{n,k}$. By the diagonal argument for double sequence, there exists a subsequence $(\e^\p(n))$ of $(\e(n))$ such that for any $k=1,2,\dots$ and for any $1\leq i\leq m_{0}$, $\mu_{\e^\p(n)}(\chi_{Q_{k}})\delta_{\e^\p(n)}(Q_{k},i)$ converges as $n\to \infty$. By Theorem \ref{th:main}, $\mu_{\e^\p(n)}\,:\,F_{b,0}(X)\to \R$ converges weakly to a finite measure $\mu=\sum_{i=1}^{m_{0}}\delta(i)\mu(U(i),\cd)$ putting $\delta(i):=\lim_{k\to \infty}\lim_{n\to \infty}\mu_{\e^\p(n)}(\chi_{Q_{k}})\delta_{\e^\p(n)}(Q_{k},i)$. Hence the proof is complete.
\proe
\subsection{Proof of Corollary \ref{cor:main}}\label{sec:proof_main2_cor}
\pros
(1) $\Rightarrow$ (2). Assume that $\mu_{\e}$ converges weakly to a measure $\mu$. By Theorem \ref{th:main}(4), there exists a sequence of finite subsets $Q_{n}\subset S$ with $Q_{n}\cap U(k)\neq \emptyset$ for all $k$, $Q_{n}\subset Q_{n+1}$ and $\bigcup_{n}Q_{n}=S$ such that $\mu_{\e}(\chi_{Q_{n}})\delta_{\e}(Q_{n},k)$ converges to a number as $\e\to 0$ for any $n$ and $k$. By convergence of $\lim_{\e\to 0}\mu_{\e}(\chi_{Q_{n}})$, the number $\delta_{\e}(Q_{n},k)$ converges to a number. Thus the assertion (2) holds.
\smallskip
\\
(2) $\Rightarrow$ (1). Assume that $\delta_{\e}(Q_{n},i)$ converges to a number $\delta(Q_{n},i)$ as $\e\to 0$ for any $n,i$.
Choose any positive sequence $(\e(k))$ with $\lim_{k\to \infty}\e(k)=0$. By condition (A.4) and by the diagonal argument for $\mu_{\e(k)}(\chi_{Q_{n}})$, there exists a subsequence $(\e^\p(k))$ of $(\e(k))$ such that $\mu_{\e^\p(k)}(\chi_{Q_{n}})$ converges to a number $a(n)$ as $k\to \infty$ for any $n$, and $a(n)\geq 1-1/n$ is satisfied. By Theorem \ref{th:main}(3), $\mu_{\e^\p(k)}$ converges weakly to a measure $\mu=\sum_{i=1}^{m_{0}}\delta(i)\mu(U(i),\cd)$ with $\delta(i)=\lim_{n\to \infty}a(n)\delta(Q_{n},i)$. Since $a(n)\to 1$, we have $\delta(i)=\lim_{n\to \infty}\delta(Q_{n},i)$. Moreover, it follows from $1=\mu_{\e}^{Q_{n}}(1)=\sum_{i=1}^{m_{0}}\ti{\delta}_{\e}(Q_{n},i)$ and $\ti{\delta}_{\e}(Q_{n},i)/\delta_{\e}(Q_{n},i)\to 1$ as $\e\to 0$ that $\sum_{i=1}^{m_{0}}\delta(i)=1$.
The limiting measure $\mu$ does not depend on choosing $(\e(n))$. Thus we obtain convergence of $\mu_{\e}$. Finally, we prove that $F_{b,0}(X)$ in Theorem \ref{th:main}(1) is replaced by $F_{b}(X)$. The notion of vanishing as infinity is used only to the inequality $|\mu_{\e}(f\chi_{S\setminus Q_{n}})|\leq \eta$ in (\ref{eq:muef=...}). This is implied from condition (A.4). Indeed, by taking large $n$ if necessary, $|\mu_{\e}(\chi_{S\setminus Q_{n}})|\leq \eta/(\|f\|_{\infty}+1)$ is satisfied for any small $\e>0$. Therefore $|\mu_{\e}(f\chi_{S\setminus Q_{n}})|\leq \|f\|_{\infty}\eta/(\|f\|_{\infty}+1)\leq \eta$ holds. Hence the assertion (1) is fulfilled.
\proe
\section{Application: a perturbation of piecewise expanding Markov maps with countably infinite partitions}\label{sec:piecewiseCI}
We consider piecewise expanding Markov maps of the interval with countably infinite partitions. We study a perturbation of such a Markov map in which there is a single transitive component while the unperturbed map has a finite number of transitive components. We investigate the phenomenon of splitting of the Gibbs measure for the perturbed Markov map.

Consider a set $(G,(J_{v}),(T_{e}))$ satisfying the following conditions (G.1)-(G.4):
\begin{itemize}
\item[(G.1)] $G=(V,E,i(\cd),t(\cd))$ is a directed multigraph endowed with finite vertex set $V$, countable edge set $E$, and two maps $i(\cd)$ and $t(\cd)$ from $E$ to $V$. For each $e\in E$, $i(e)$ is called the initial vertex of $e$ and $t(e)$ called the terminal vertex of $e$.
\item[(G.2)] For each $v\in V$, $J_{v}$ is a closed interval with non-empty interior satisfying that $\bigcup_{v\in V}J_{v}=[0,1]$, and $\mathrm{int} J_{v}$ and $\mathrm{int}J_{v^\p}$ are disjoint for $v^\p\in V$ with $v\neq v^\p$.
\item[(G.3)] For each $e\in E$, $T_{e}\,:\,J_{t(e)}\to J_{i(e)}$ is a map satisfying the following:
\ite
{
\item[(i)] $T_{e}$ is injective and of class $C^{1}$. The expression $|T_{e}^\p(x)|$ for $x\in \partial J_{t(e)}$ denotes the one-sided derivative. 
\item[(ii)] (Expanding) $\sup_{e\in E}\sup_{x\in J_{t(e)}}|T_{e}^\p(x)|<1$.
\item[(iii)] (Open set condition) For $e,e^\p\in E$ with $e\neq e^\p$, $T_{e}(\mathrm{int}J_{t(e)})\cap T_{e^\p}(\mathrm{int}J_{t(e^\p)})=\emptyset$.
\item[(iv)] (Bounded distortion) There exist constants $c_{\adl{ibd}}>0$ and $0<\beta\leq 1$ such that for any $e\in E$ and $x,y\in J_{t(e)}$, $||T_{e}^\p(x)|-|T_{e}^\p(y)||\leq c_{\adr{ibd}}|T_{e}^\p(x)| |x-y|^{\beta}$.
\item[(v)] For each $v\in V$, $\overline{\bigcup_{e\in E\,:\,i(e)=v} T_{e}(J_{t(e)})}=J_{v}$.
\item[(vi)] (Summability) $\sum_{e\in E}\sup_{x\in J_{t(e)}}|T_{e}^\p(x)|<\infty$.
}
\end{itemize}
The {\it incidence matrix} $A=A_{G}$ of $G$ is a zero-one matrix indexed by $E$ such that $A(ee^\p)=1$ if $t(e)=i(e^\p)$ and $A(ee^\p)=0$ if $t(e)\neq i(e^\p)$. Consider the TMS
\alil
{
\textstyle X=X_{G}=\{\om\in \prod_{n=0}^{\infty}E\,:\,A(\om_{n}\om_{n+1})=1 \text{ for all }n\geq 0\}.\label{eq:E^inf=}
}
with countable state space $E$ and with transition matrix $A$. The coding map $\pi\,:\,X\to \R^{1}$ is well defined by
$\{\pi\om\}=\bigcap_{k=0}^{\infty}T_{\om_{0}\cdots\om_{k}}(J_{t(\om_{k})})$, where $T_{\om_{0}\cdots \om_{k}}$ means $T_{\om_{0}}\circ \cdots\circ T_{\om_{k}}$. The {\it limit set} of the system is given by the image $K:=\pi(X)$. 

We define a map $f_{0}\,:\,[0,1]\to [0,1]$ by $f_{0}(x)=T_{e}^{-1}(x)$, where $e$ is decided uniquely if $x\in \bigcup_{e\in E}\mathrm{int}(T_{e}(J_{t(e)}))$ and otherwise we arbitrarily choose $e$ so that $x\in \partial T_{e}(J_{t(e)})$. Then $f_{0}$ is a piecewise expanding map with an infinite Markov partition.
\ali
{
\graphexp{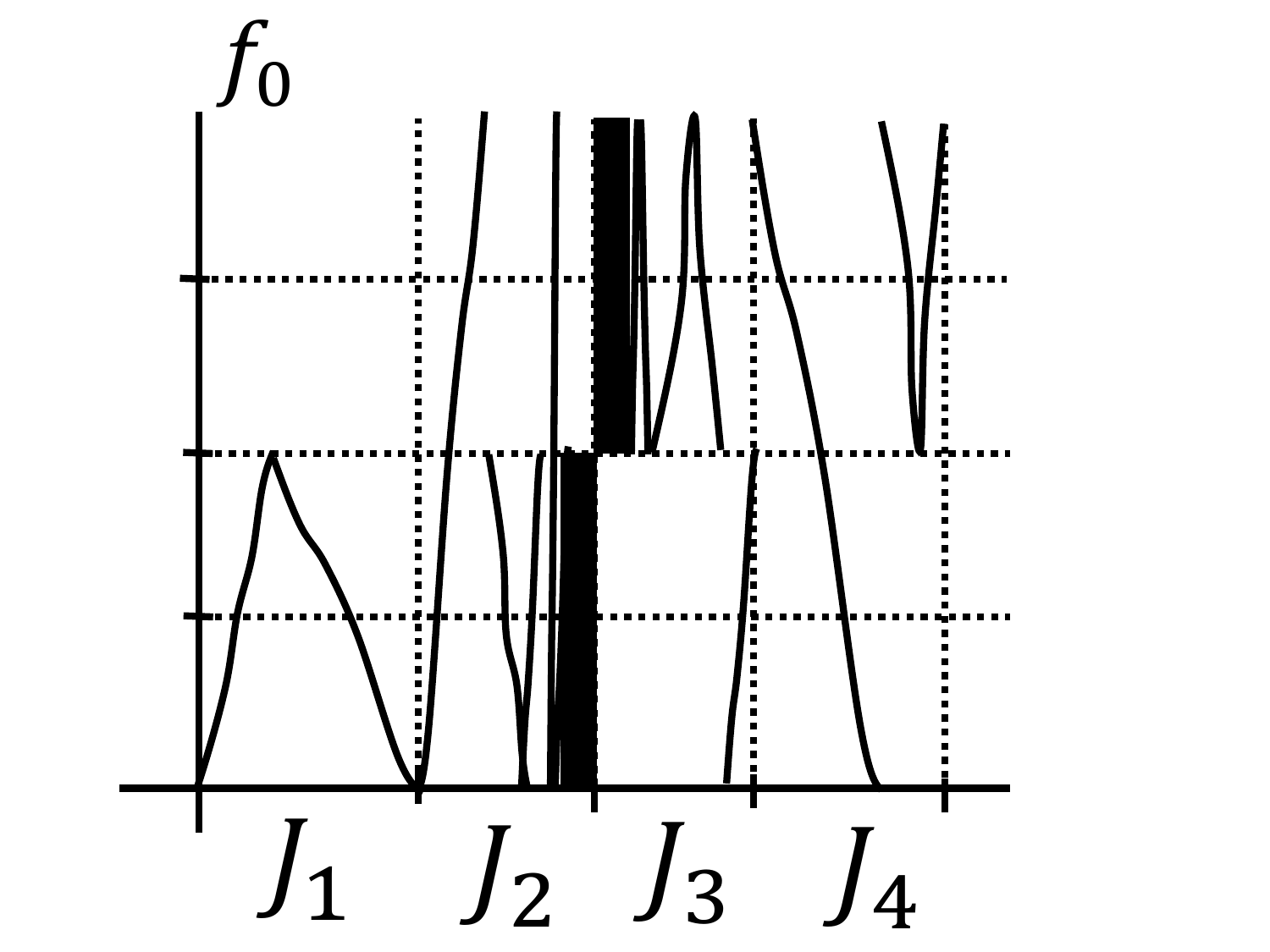}{0.16}{2.5}{1}
}
Assume that the graph $G$ is strongly connected, i.e., for any $v,u\in V$, there exists a directed path from $v$ to $u$. The set of all critical points of $f_{0}$ is given by $C_{0}=\bigcup_{e\in E}\partial T_{e}(J_{t(e)})$. It is known that the sets $D_{0}=\bigcup_{n=0}^{\infty}f_{0}^{-n}C_{0}$ and $\pi^{-1}D_{0}$ are at most countable. Moreover, for any $x\in \bigcup_{v\in V}J_{v}\setminus D_{0}$, there exist $e_{1},e_{2},\dots\in E$ and $x_{1},x_{2},\dots \in \bigcup_{v\in V}J_{v}\setminus D_{0}$ such that $x=T_{e_{1}}\circ \cdots T_{e_{n}}(x_{n})$ for any $n\geq 1$. We have $x\in \pi(X)$ and thus
\ali
{
\bigcup_{v\in V}J_{v}\setminus D_{0}\subset \pi(X)
}
In particular, we have $\dim_{H}\pi(X)=1$. We define a function $\ph\,:\,X\to \R$ by
\ali
{
\ph(\om)=\log|T_{\om_{0}}^\p(\pi\si\om)|.
}
Then it is known that $P(\ph)=0$. We take $(1,h,\nu)$ as the spectral triplet of the Ruelle operator $\LR_{\ph}$. Put $\mu=h\nu$. Then the absolutely continuous invariant probability measure (ACIM) of $f_{0}$ coincides with the measure $\mu\circ \pi^{-1}$ ($\mathrm{mod}$ $0$). In particular, $\nu\circ \pi^{-1}$ equals the Lebesgue measure. In particular, $\mu$ is the Gibbs measure for the potential $\ph$.

We formulate a perturbation of the map $f_{0}$:
\ite
{
\item[(H)] A triplet $(G,(J_{v}),(T_{e}(\e,\cd)))$ is satisfied (G.1)-(G.3) with a small parameter $0<\e<1$ s.t.
\ite
{
\item[(i)] $G$ is strongly connected.
\item[(ii)] There exist a non-empty subset $E_{0}\subset E$, maps $T_{e}\,:\,J_{t(e)}\to J_{i(e)}$ $(e\in E_{0})$, and $a_{e}\in J_{i(e)}$ $(e\in E\setminus E_{0})$ such that
\ali
{
\sup_{e\in E_{0}}\sup_{x\in J_{t(e)}}|T_{e}(\e,x)-T_{e}(x)|&\to 0\\
\sup_{e\in E\setminus E_{0}}\sup_{x\in J_{t(e)}}|T_{e}(\e,x)-T_{e}(x)|&\to 0\\
\sup_{x\in J_{t(e)}}|{\textstyle\frac{\partial}{\partial x}}T_{e}(\e,x)-T_{e}^\p(x)|&\to 0\quad \text{for any }e\in E_{0}\\
\sup_{x\in J_{t(e)}}|{\textstyle\frac{\partial}{\partial x}}T_{e}(\e,x)|&\to 0\quad \text{for any }e\in E\setminus E_{0}
}
\item [(iii)] The subsystem $(G_{0}=(V_{0},E_{0}),(J_{v})_{v\in V_{0}},(T_{e})_{e\in E_{0}})$ is satisfied conditions (G.1)-(G.3).
\item[(iv)] There exist constants $c_{\adl{pbd}}>0$ and $0<\beta\leq 1$ such that for any $e\in E$ and $x,y\in J_{t(e)}$, $||T_{e}^\p(\e,x)|-|T_{e}^\p(\e,y)||\leq c_{\adr{pbd}}|T_{e}^\p(\e,x)||x-y|^{\beta}$.
\item[(v)] For each $v\in V$, there exists a unique strongly connected subgraph $H=(V_{H},E_{H})$ of $G$ such that $v\in V_{H}$ and $\overline{\bigcup_{e\in E_{H}\,:\,i(e)=v} T_{e}(J_{t(e)})}=J_{v}$.
\item[(vi)] $\sum_{e\in E}\sup_{\e>0}\sup_{x\in J_{t(e)}}|\frac{\partial}{\partial x}T_{e}(\e,x)|<+\infty$.
}
}
For any $\e>0$, the perturbed piecewise expanding  Markov map $f_{\e}$ with respect to $(G,(J_{v}), (T_{e}(\e,\cd)))$ has only one transitive component. On the other hand, the unperturbed piecewise expanding Markov map $f_{0}$ with respect to $(G_{0},(J_{v}), (T_{e}))$ has a finite number of transitive components. Let $SC(G_{0})$ be the set of all strongly connected components of the graph $G_{0}$. For $H=(V_{H},E_{H})\in SC(G_{0})$, we let $I_{H}=\bigcup_{v\in V_{H}}J_{v}$. 
\ali
{
\graphexp{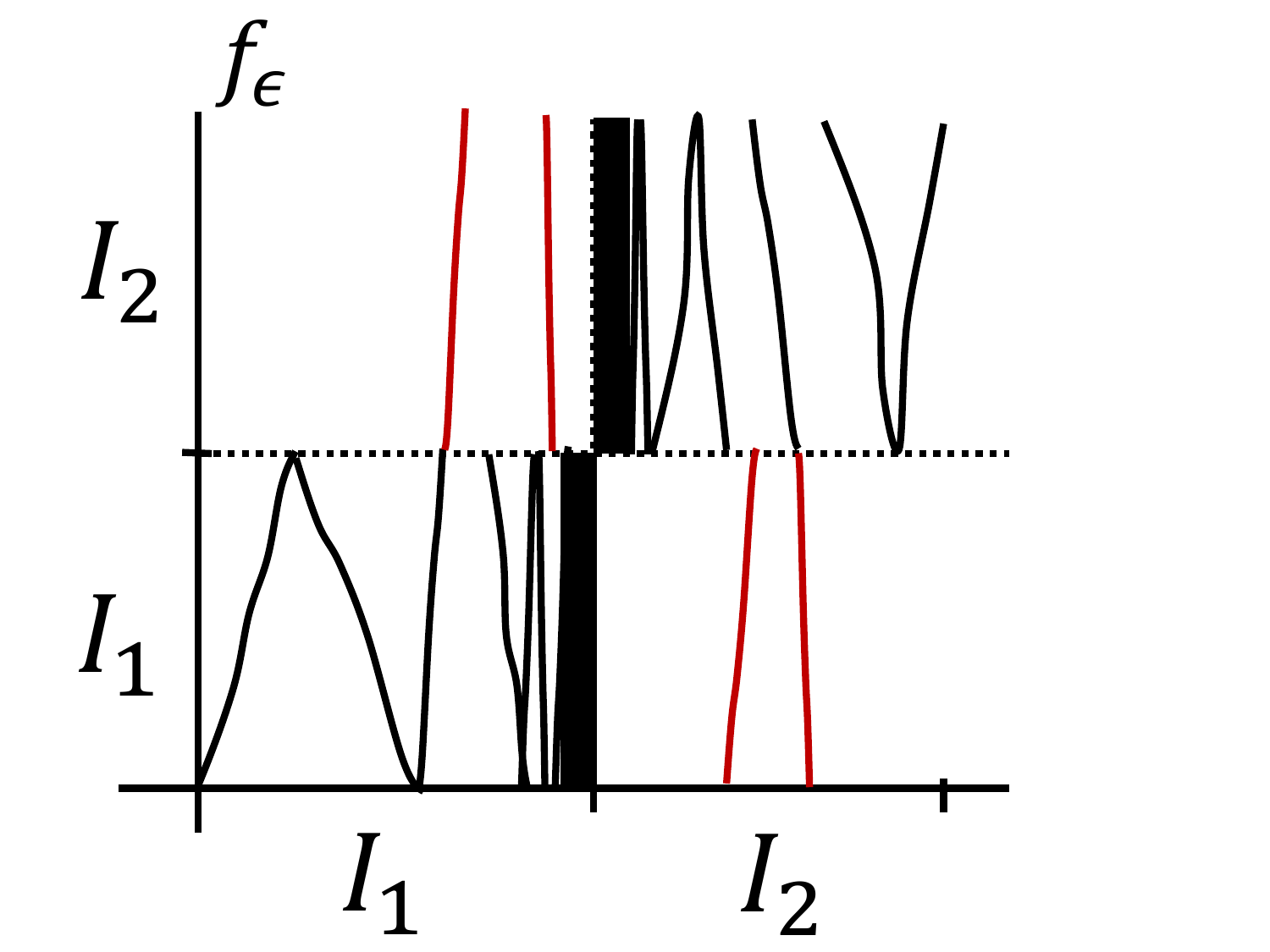}{0.18}{3.5}{1}
\underset{\e\to 0}{\rightarrow}\quad
\graphexp{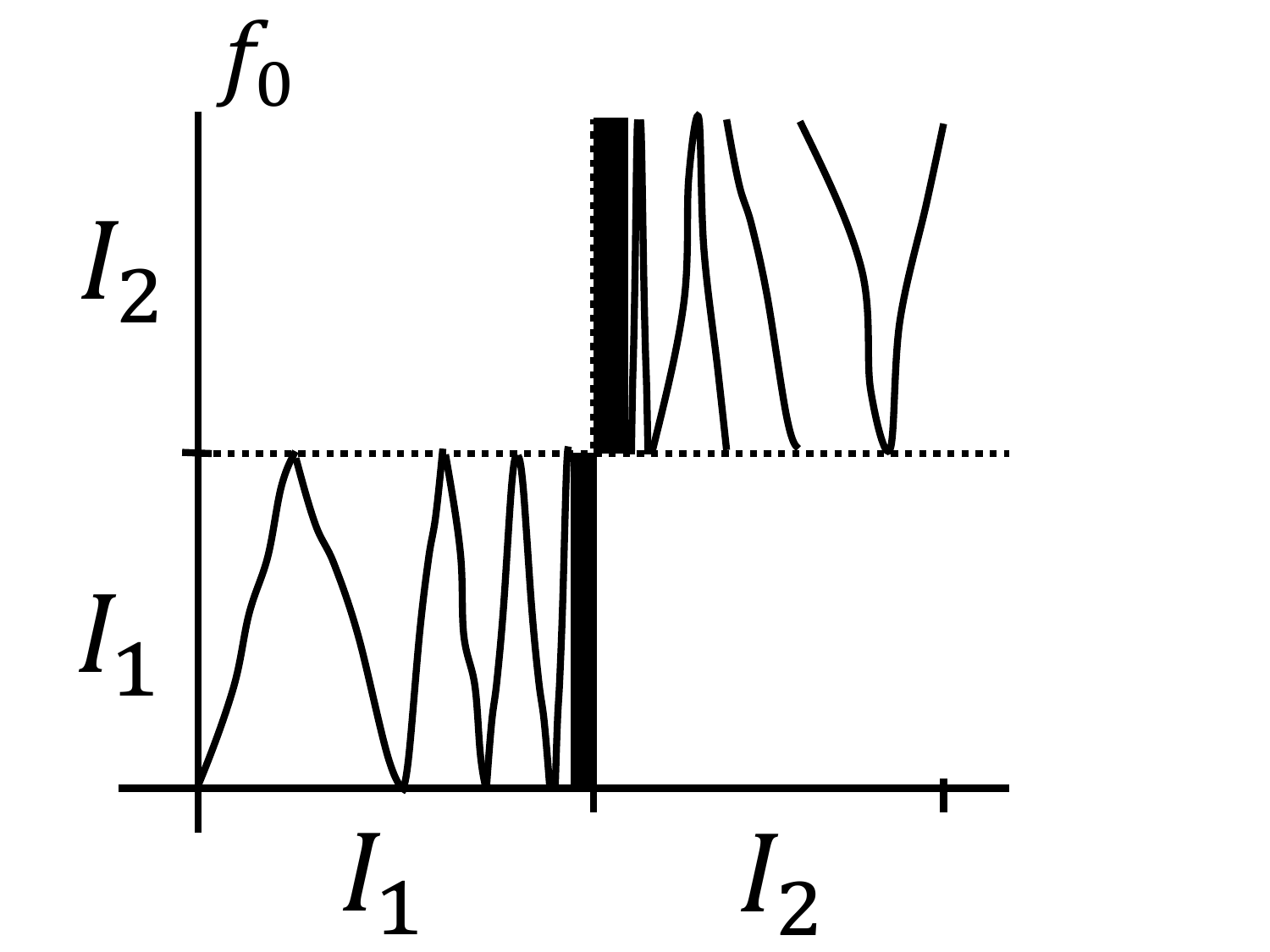}{0.18}{3.5}{1}
}
Let $\pi(\e,\cd)$ be the limit set of $f_{\e}$, and we define
\ali
{
\ph(\e,\om)=\log |{\textstyle\frac{\partial}{\partial x}}T_{\om_{0}}(\e,\pi(\e,\si\om)|.
}
To state our main result, we need more notation. We write
\ali
{
\hat{E}_{H}=\{e\in E\,:\,i(e)\in V_{H}\}
}
for $H\in SC(G_{0})$ and 
\ali
{
\{E(1),E(2),\dots, E(m_{0})\}=\{\hat{E}_{H}\,:\,H\in SC(G_{0})\}.
}
Note that $E_{H}\subset \hat{E}_{H}$. Denote by $\mu_{\e}$ the Gibbs measure of $\ph(\e,\cd)$, and by $\mu(k,\cd)$ the Gibbs measure of $\ph|_{X_{H(k)}}$ for each $1\leq k\leq m_{0}$.

For distinct integers $k,k^\p$ and a finite subset $Q\subset E$, let $c_{\e}(Q, k,k^\p)$ be the spectral radius of $\LR_{\e}[E(k)\cap Q,E\setminus ((E(k)\cup E(k^\p))\cap Q),1]$. For $1\leq k\leq m_{0}$
\ali
{
p_{\e}(Q,k):=\frac{1}{\di 1+\sum_{k^\p\,:\,k^\p\neq k}\frac{1-c_{\e}(Q, k,k^\p)}{1-c_{\e}(Q, k^\p,k)}}.
}
Now we are in a position to state our main result.
\thm
{\label{th:Gibbsiff_PEMM}
Assume that condition (H) is satisfied. Assume also that $m_{0}:=\sharp SC(G_{0})\geq 2$. Then the following are equivalent.
\ite
{
\item The Gibbs measure $\mu_{\e}$ converges to a measure $\mu$ weakly.
\item There exists a sequence of finite subsets $Q_{n}\subset S$ with $Q_{n}\cap U(k)\neq \emptyset$ for all $k$, $Q_{n}\subset Q_{n+1}$ and $\bigcup_{n}Q_{n}=S$ such that $p(n,k):=\lim_{\e\to 0}p_{\e}(Q_{n},k)$ exists for all $n\geq 1$ and for all $1\leq k\leq m_{0}$.
}
In this case, the measure $\mu$ has the convex combination form $\sum_{k=1}^{m_{0}}p(k)\mu(k,\cd)$ with $p(k)=\lim_{n\to \infty}p(n,k)$.
}
To prove above theorem, we need to show some auxiliary propositions.
For $H\in SC(G_{0})$, 
\ali
{
X_{H}=\{\om\in X\,:\,\om_{n}\in E_{H} \text{ for any }n\geq 0\}.
}
\prop
{\label{prop:finite_SC}
Assume that condition (H) is satisfied. Then the number of $SC(G_{0})$ is finite.
}
\pros
The vertex sets $V_{H}$ ($H\in SC(G_{0})$) are disjoint and are subsets of the finite set $V$. Hence $SC(G_{0})$ must be finite.
\proe
\prop
{\label{prop:SC_dimone}
Assume that condition (H) is satisfied. Then $\dim_{H}\pi(\e,X)=1$ for all $\e>0$, and $\dim_{H}\pi(X_{H})=1$ for any $H\in SC(G_{0})$.
}
\pros
The subsystem $(H,(J_{v})_{v\in V_{H}},(T_{e})_{e\in E_{H}})$ satisfies conditions (G.1)-(G.3) replacing $G$ by $H$. Hence each limit set $\pi(X_{H})$ has one Hausdorff dimension.
\proe
For $H\in SC(G_{0})$, we set
\ali
{
\hat{E}_{H}=\{e\in E\,:\,i(e)\in V_{H}\}
}
\prop
{\label{prop:partition_E}
Assume that condition (H) is satisfied. Then $E$ is decomposed into $\hat{E}_{H}$ ($H\in SC(G_{0})$).
}
\pros
By condition (H)(v), any $e\in E$ belongs to $V_{H}$ for some $H\in SC(G_{0})$. Moreover, since for any distinct subgraphs $H,H^\p\in SC(G_{0})$, $V_{H}$ and $V_{H^\p}$ are disjoint. Then $\hat{E}_{H}\cap \hat{E}_{H^\p}=\emptyset$. Hence the assertion is valid.
\proe
For $e\in E\setminus E_{0}$, we define $T_{e}\,:\,J_{t(e)}\to J_{i(e)}$ by $T_{e}\equiv a_{e}$. Then the map $\pi\,:\,X\to \R$ is given as the defition of the coding map, i.e. $\pi\om=\bigcap_{n=0}^{\infty}T_{\om_{0}}\circ\cdots\circ T_{\om_{n}}(J_{t(\om_{n})})$.
\prop
{\label{prop:conv_pie}
Assume that condition (H) is satisfied. Then $\sup_{\om\in X}|\pi(\e,\om)-\pi\om|\to 0$ as $\e\to 0$. Moreover, there exists a constant $c_{\adl{pi0}}>0$ such that for any $\om,\up\in X$ with $i(\om)=i(\up)$, $|\pi(\e,\om)-\pi(\e,\up)|\leq c_{\adr{pi0}}d_{r}(\om,\up)$, where $r=\sup_{e\in E}\sup_{x\in J_{t(e)}}|T^\p_{e}(x)|$.
}
\pros
See \cite[Lemma 2.12]{T2025_IGIFSdege}.
\proe
\prop
{\label{prop:summ_PEMM}
Assume that condition (H) is satisfied. Then the uniform summability $\sum_{e\in E}\sup_{\e>0}\sup_{\om\in [e]}\exp(\ph(\e,\om))<+\infty$ holds. Therefore condition (A.2) holds.
}
\pros
See \cite[Lemma 2.13]{T2025_IGIFSdege}.
\proe
\prop
{\label{prop:1-vari_PEMM}
Assume that condition (H) is satisfied. Then $\sup_{\e>0}[\phe]_{1}<+\infty$. In particular, condition (A.1) holds.
}
\pros
See \cite[Lemma 2.14]{T2025_IGIFSdege}.
\proe
Put
\ali
{
\ph_{0}(\om):=
\case
{
\log|T_{\om_{0}}^\p(\pi\si\om)|,&\om_{0}\in E_{0}\\
0,&\om_{0}\in E\setminus E_{0}\\
}\qquad
\psi(\om):=
\case
{
e^{\ph_{0}(\om)},&\om_{0}\in E_{0}\\
0,&\om_{0}\in E\setminus E_{0}\\
}
}
\prop
{\label{prop:conv_expphie}
Assume that condition (H) is satisfied. Then $\sup_{\om\in [e]}|e^{\ph(\e,\om)}-\psi(\om)|\to 0$ as $\e\to 0$ for each $e\in E$.  Therefore condition (A.3) holds.
}
\pros
See \cite[Lemma 2.15]{T2025_IGIFSdege}.
\proe
Recall the spectral triplet $(1,h_{\e},\nu_{\e})$ and $g_{\e}=h_{\e}/\|h_{\e}\|_{\infty}$.
\prop
{\label{prop:positivebdd_ge^E}
Assume that condition (H) is satisfied. Then there exists a constant $c_{\adl{bgH}}>0$ such that for any $H\in SC(G_{0})$, $g_{\e}^{\hat{E}_{H}}\geq c_{\adr{bgH}}\chi_{\hat{E}_{H}}$ for any small $\e>0$.
}
\pros
We will check the additional conditions (i)-(iii) of Proposition \ref{prop:gP_positive_inf}.
For each $v\in V_{H}$, we take $e(v)\in \hat{E}_{H}$ with $i(e(v))=v$. Then for any $e\in E$ with $A(ee^\p)=1$ for some $e^\p\in \hat{E}_{H}$, we have $A(ee(v))=1$ letting $v:=t(e)\in V_{H}$. Therefore the additional condition of Proposition \ref{prop:gP_nonzero} by putting $P_{0}:=\{e(v)\,:\,v\in V_{H}\}$. By this proposition, we have $\liminf_{\e\to 0}g_{\e}^{\hat{E}_{H}}\neq 0$ and condition (i) holds. Moreover, for any $b\in \hat{E}_{H}$, $B(ab)=1$ for some $a\in Q_{0}:=P_{0}$. Thus condition (ii) is satisfied. finally, condition (iii) follows from Proposition \ref{prop:1-vari_PEMM}. Hence the assertion holds from Proposition \ref{prop:gP_positive_inf}.
\proe
Write $\LR_{\e}$ as the Ruelle operator for the potential $\phe$, and $\LR$ as the Ruelle operator for the potential $\ph$.
\prop
{\label{prop:conv_mueEH}
Assume that condition (H) is satisfied. Then for any $H\in SC(G_{0})$, $\mu_{\e}^{\hat{E}_{H}}$ converges to $\mu(H,\cd)$ weakly.
}
\pros
Note that $\supp\, \nu_{\e}^{\hat{E}_{H}}=X\cap (\hat{E}_{H})^{\Z_{+}}=X\cap (E_{H})^{\Z_{+}}=\supp\, \nu_{\e}^{E_{H}}$, we have the equation $\mu_{\e}^{\hat{E}_{H}}=\mu_{\e}^{E_{H}}$.
Together with Proposition \ref{prop:positivebdd_ge^E}, the assertion follows from Proposition \ref{prop:conv_PC_mue}.
\proe
\prop
{\label{prop:vani_mu}
Assume that condition (H) is satisfied. Then condition (A.4) holds.
}
\pros
Choose any $\eta\in (0,1)$. 
Consider the form $\mu_{\e}=\sum_{H\in SC(G_{0})}\mu(\Si_{\hat{E}_{H}})\mu_{\e}^{\hat{E}_{H}}$. Since $\mu(H,\cd)$ is a probability measure, there exists a finite set $Q\subset E$ such that $\mu(H,\Si_{E\setminus Q})\leq \eta/2$. 
By Proposition \ref{prop:conv_mueEH}, $\mu_{\e}^{\hat{E}_{H}}(\Si_{E\setminus Q})\to \mu(H,\Si_{E\setminus Q})$. Therefore, for any small $\e>0$, $\mu_{\e}^{\hat{E}_{H}}(\Si_{E\setminus Q})\leq \eta$. We have
\ali
{
\mu_{\e}(\Si_{E\setminus Q})\leq \sum_{H\in SC(G_{0})}\mu(\Si_{\hat{E}_{H}})\eta=\eta.
}
Hence condition (A.4) is fulfilled.
\proe
\pros[Proof of Theorem \ref{th:Gibbsiff_PEMM}.] 
Together with Proposition \ref{prop:summ_PEMM}, Proposition \ref{prop:1-vari_PEMM}, Proposition \ref{prop:conv_expphie}, and Proposition \ref{prop:vani_mu}, the assertion follows from Corollary \ref{cor:main}.
\proe

\endthebibliography
\end{document}